\documentclass{amsart}

\usepackage{graphicx}
\usepackage{enumerate}
\usepackage{amssymb}

\newtheorem{theorem}{Theorem}[section]
\newtheorem{lemma}[theorem]{Lemma}
\newtheorem{proposition}[theorem]{Proposition}

\newtheorem{mtheorem}{Theorem}

\theoremstyle{definition}
\newtheorem{definition}[theorem]{Definition}
\newtheorem{remark}[theorem]{Remark}
\newtheorem{setting}[theorem]{Setting}

\numberwithin{equation}{section}
\numberwithin{figure}{section}

\begin{document}

\title[Cubic homoclinic tangencies]{ 
Persistent antimonotonic bifurcations and strange attractors for cubic homoclinic tangencies}

\author{Shin Kiriki} 
\address{Department of Mathematics, Kyoto University of Education, 
1 Fukakusa-Fujinomori, Fushimi-ku, Kyoto, 612-8522, Japan}
\email{skiriki@kyokyo-u.ac.jp}

\author{Teruhiko Soma}
\address{Department of Mathematics and Information Sciences, Tokyo Metropolitan
University, Minami-Ohsawa 1-1, Hachioji, Tokyo 192-0397, Japan}
\email{tsoma@tmu.ac.jp}

\date{\today}

\subjclass[2000]{Primary:  37C29, 37D45; Secondary: 37D20}

\keywords{Cubic homoclinic tangency, contact-making, contact-breaking tangencies, strange attractors}

\maketitle

\begin{abstract}
In this paper, we study a two-parameter family $\{\varphi_{\mu,\nu}\}$ of two-dimensional diffeomorphisms such that $\varphi_{0,0}=\varphi$ has a cubic homoclinic tangency unfolding generically which is associated with a dissipative saddle point.
Our first theorem presents an open set $\mathcal{O}$ in the $\mu\nu$-plane with $\mathrm{Cl}(\mathcal{O})\ni (0,0)$ such that, for any $(\mu_0,\nu_0)\in \mathcal{O}$, there exists a one-parameter subfamily of $\{\varphi_{\mu,\nu}\}$ passing through $\varphi_{\mu_0,\nu_0}$ and exhibiting cubically related persistent contact-making and 
contact-breaking quadratic tangencies.
Moreover, the second theorem shows that any such two-parameter family satisfies Wang-Young's conditions 
which guarantee that some $\varphi_{\mu,\nu}$ arbitrarily near $\varphi$ exhibits a cubic polynomial-like strange attractor with an SRB measure.
\end{abstract}

\setcounter{section}{-1}

\section{Introduction}
One of main motivations behind studies in monotonic or nonmonotonic bifurcation phenomena is attributed to results given in Milnor-Thurston \cite{MT} for the logistic family, 
where there are only orbit creation parameters but no orbit annihilation parameters.
In contrast with the  logistic family, 
some family of one-dimensional cubic maps has nonmonotonic bifurcations, for example see Milnor-Tresser \cite{MTr} and Li \cite{Li}. 
The main result  in this paper is that the unfolding of cubic homoclinic tangencies for two-dimensional diffeomorphisms generates the simultaneous creation and destruction of quadratic tangencies in a persistent way. 
Nonmonotonic phenomena near a quadratic homoclinic tangency have been already obtained 
 by Kan-Ko\c{c}ak-Yorke \cite{KKY}.
Also, the nonmonotonic 
phenomena and  essential mechanism of their generations
associated with any order homoclinic tangencies have been studied extensively by 
Gonchenko, Shilnikov and Turaev\ in  \cite{G93}-\cite{G07}. 
In this paper, we will give generic conditions which guarantee that contact-making and contact-breaking 
tangencies occur simultaneously as well as  persistently for two-parameter families of planar diffeomorphisms, 
which are generated by a mechanism different from that in \cite{KKY}.
D\'{\i}az and Rocha \cite{DR} gave another example of bifurcations  closely related to ours which are not monotonic but have heterodimensional cycles in dimension equal to or greater than three.

As for phenomena of the cubic tangency, 
Bonatti,  D\'{\i}az and  Vuillemin \cite{BDV} showed that there exists a three-parameter family of diffeomorphisms  in the boundary of the  Anosov diffeomorphisms which has 
the first cubic heteroclinic tangency, see also Enrich \cite{E}.
Carvalho \cite{Ca} proved the corresponding result where the cubic tangency is homoclinic.
In \cite{Ca}, she also presented some supporting evidence for the existence of cubic tangencies in the H\'enon family.

In order to set main results in our context, we first recall some previous results.
One of the key facts, which was presented by Newhouse \cite{N1}, is that the unfolding of 
a quadratic homoclinic tangency leads to persistence of tangencies, see also \cite{Dav, Kal, LS, N2, R, TY} for related results.
Palis and Takens \cite{PT} gave a new proof of this result which has three major steps: 
(i) there is a renormalization at
quadratic homoclinic tangencies whose limit dynamics is controlled by the logistic map, 
(ii) this map is conjugate to the tent map, and (iii) the tent map has Cantor invariant sets
with arbitrarily large thickness.
Here, the thickness is one of main tools with Gap Lemma (see Subsection \ref{thickness}) to get the persistent non-transverse intersections.

\begin{figure}[htb]
\begin{center}
\includegraphics[clip]{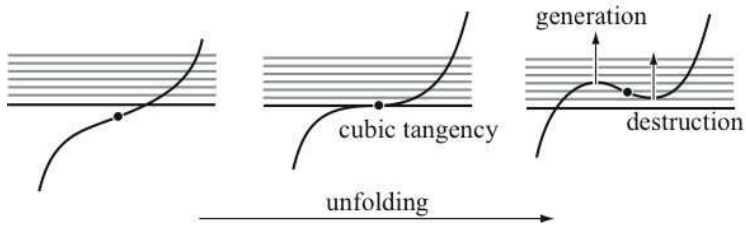}
\caption{ Unfolding of a cubic tangency.}
\label{fg_0_1}
\end{center}
\end{figure}

In this paper, we obtain counterparts of these three objects in the cubic tangency case. 
The logistic map is replaced by a cubic map and the tent map is by the N-map defined in Subsection \ref{affine}.
Then, one can adapt the techniques of Palis-Takens to our context. 
The renormalization rescaling for cubic tangencies in Subsection \ref{renorm} closely 
follows arguments in \cite[pp.\ 47-51]{PT} where the corresponding results for
the quadratic case are studied. 
From these adaptations and reformulations, we know that the unfolding of a cubic homoclinic tangency provides the unfolding of
a pair of persistent quadratic tangencies one of which generates  intersections and the other destructs intersections simultaneously, see Fig.\ \ref{fg_0_1}.

We first prove the following theorem.

\begin{mtheorem}\label{main_a}
Suppose that $\varphi$ is any $C^\infty$-diffeomorphism on the plane $\mathbb{R}^2$ with a dissipative 
saddle fixed pint $p$ at which $\varphi$ satisfies the open $C^4$-linearizing condition and 
such that the stable and unstable manifolds of $p$ have a cubic homoclinic tangency $q$.
Let $\{\varphi_{\mu,\nu}\}$ be any two-parameter family in $\mathrm{Diff}^\infty(\mathbb{R}^2)$ with 
$\varphi_{0,0}=\varphi$ such that the cubic tangency unfolds generically.
Then, there exists an open subset $\mathcal{O}$ in the $\mu\nu$-plane with $\mathrm{Cl}(\mathcal{O})\ni (0,0)$ 
and such that, for any $(\mu,\nu)\in \mathcal{O}$, there is a regular curve $c:(-\varepsilon,\varepsilon) 
\longrightarrow \mathcal{O}$ such that the one-parameter family $\{\varphi_{c(t)}\}$ exhibits cubically 
related persistent contact-making and contact-breaking tangencies.
\end{mtheorem}

Here, we say that the saddle fixed point $p$ is \emph{dissipative} if 
$\vert\det (d\varphi)_p\vert <1$ and the cubic tangency \emph{unfolds generically} with respect to $\{\varphi_{\mu,\nu}\}$ if it satisfies the generic condition (\ref{generic condition}) given in Subsection \ref{initial_set}.
The $\varphi$ is said to satisfy the \emph{open $C^4$-linearizing condition} at 
$p$ if, for any $\psi\in \mathrm{Diff}^\infty(\mathbb{R}^2)$ sufficiently close to $\varphi$, 
there exists a $C^4$-coordinate on a small neighborhood $U$ of $p$ with respect to which 
the restriction $\psi|_U$ is a linear map.
This linearizability condition is not so strong.
In fact, the \emph{$C^4$-Sternberg condition} given in \cite{S, T0} is a generic sufficient condition for it.
The curve $c$ is \emph{regular} is $dc/dt(t)\neq (0,0)$ for any $t\in (-\varepsilon,\varepsilon)$ 
and the family $\{\varphi_{c(t)}\}$ has persistent \emph{contact-making tangencies} (resp.\ \emph{contact-breaking tangencies}) if, for any $t\in (-\varepsilon,\varepsilon)$, $\varphi_{c(t)}$ has a tangency $r_t$ in a neighborhood $\mathcal{N}(q)$ of $q$ associated to basic sets and such that $r_t$ generates (resp.\ annihilates) transverse intersections in $\mathcal{N}(q)$ as $t$ increases, see Subsection \ref{p_tangency} for the precise 
definition.
The situation such that a pair of persistent contact-making and contact-breaking tangencies are 
\emph{cubically related} is also explained there.

\begin{remark}
Theorem \ref{main_a} is not a new result in the sense that it detects some antimonotonic persistent tangencies 
near cubic homoclinic bifurcations.
In fact, we know that some applications of results in Kan-Ko\c{c}ak-Yorke \cite{KKY} and Newhouse \cite{N1} 
imply the existence of such persistent tangencies. 
If an unstable segment $l^u_t$ and a stable segment $l^s_t$ have a contact-making homoclinic tangency, 
then by using Bubble Lemma in \cite{KKY} we have another unstable segment in a small neighborhood of $l^u_t$ 
which has a contact-breaking tangency with $l^s_t$. 
Moreover, by invoking \cite{N1}, one can find a connected open set of diffeomorphisms having
persistent contact-making and contact-breaking tangencies.
These results are stronger than that given in Theorem \ref{main_a} in some sense, for example 
our open set $\mathcal{O}$ given in Theorem \ref{main_a} is not connected. 
However, in our theorem, we detect antimonotonicity phenomena of persistent bifurcations  
arbitrarily near the cubic tangency which are generated by a mechanism different from that by Bubble Lemma. 
\end{remark}

As for another application of our renormalization around cubic tangency, 
we show in Theorem \ref{main_b} that, for any family $\{\varphi_{\mu,\nu}\}$ of 2-dimensional diffeomorphisms with a cubic tangency unfolding generically and satisfying certain generic conditions, there exists a positive measure subset of the parameter space for each element $(\mu,\nu)$ of which $\varphi_{\mu,\nu}$ has a strange attractor with an SRB measure.
Our proof is accomplished by checking that the limit dynamics associated to the renormalization verify the hypothesis of the models in Wang-Young \cite{WY}.

A \emph{strange attractor} $\Omega$ of a diffeomorphism $\varphi$ of $\mathbb{R}^2$ is a non-hyperbolic invariant set of special type, see Subsection \ref{strange_attractor} for its precise definition.
In particular, $\Omega$ has a point $z_0$ whose positive orbit is dense in $\Omega$ and such that $\varphi$ has a positive Lyapunov exponent at $z_0$.
A $\varphi$-invariant, Borel probability measure $\boldsymbol{\mu}$ is called an \emph{SRB measure} if $\varphi$ has a positive Lyapunov exponent $\boldsymbol{\mu}$-almost everywhere and the conditional measures of $\boldsymbol{\mu}$ on unstable manifolds are absolutely continuous with respect to the Lebesgue measure on these manifolds.

\begin{mtheorem}\label{main_b}
Let $\{\varphi_{\mu,\nu}\}$ be a two-parameter family in $\mathrm{Diff}^{\infty}(\mathbb{R}^2)$ satisfying the same conditions as in Theorem \ref{main_a}.
Then, there exists a subset $\mathcal{Z}$ of the $\mu\nu$-plane with positive $2$-dimensional Lebesgue measure on any neighborhood of $(0,0)$ and such that, for any $(\mu,\nu)\in \mathcal{Z}$, $\varphi_{\mu,\nu}$ exhibits a strange attractor $\Omega_{\mu,\nu}$ supported by an SRB measure.
Moreover, the strange attractor $\Omega_{\mu,\nu}$ is different from quadratic H\'enon-like ones.
\end{mtheorem}

\begin{figure}[hbt]
\begin{center}
\scalebox{0.85}{\includegraphics[clip]{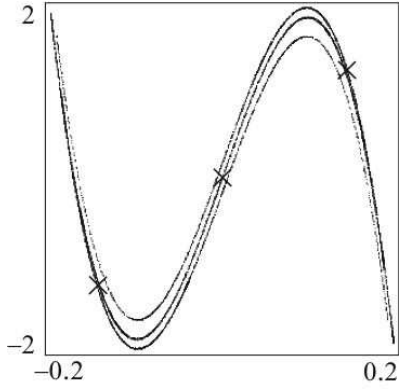}}
\caption{ The numerical results suggest the cubic polynomial-like family 
$(x,y)\mapsto(0.1y, -y^{3}+ay+x)$ would have a matured strange attractor with 
three saddle fixed points at $\times$-indications for $a\approx 2.8$.}
\label{fg_0_2}
\end{center}
\end{figure}

\begin{remark}
From the property of $\mathcal{Z}$, $\mathrm{Cl}(\mathcal{Z})$ contains $(0,0)$. 
Hence, $\varphi=\varphi_{0,0}$ is well approximated by $\varphi_{\mu,\nu}$ such that $\varphi_{\mu,\nu}^{N+n}$ 
exhibits a cubic polynomial-like strange attractor $\Omega_{\mu,\nu}$ containing three fixed points for some integer $n+N>0$, 
see Fig.\ \ref{fg_0_2}. 
On the other hand, quadratic H\'enon-like strange attractors can not contain  
three saddle fixed points.
This means that our strange attractor $\Omega_{\mu,\nu}$ is different from quadratic H\'enon-like ones which have been studied by Benedicks-Carleson \cite{BC}, Mora-Viana \cite{MV} and so on.
\end{remark}

\noindent\textbf{Acknowledgment.} 
We would like to thank M. Asaoka,  N. Sumi and M. Tsujii for helpful discussions and 
to the referees for their useful comments and suggestions.
We also thank S.\,V. Gonchenko 
for his detailed explanation of his works 
on homoclinic dynamics closely related to ours.

\section{Initial settings}\label{Initial}

\subsection{Transverse points and tangencies}\label{tangency}
We say that  a $\varphi$-invariant set 
$\Lambda$ is  \textit{basic}  if 
it  is  
compact, transitive, hyperbolic and  has  a dense subset of periodic orbits.
A diffeomorphism $\varphi$ on $\mathbb{R}^2$ has a \textit{transverse point} $r$ associated with 
non-trivial basic sets $\Lambda_{1}$ and $\Lambda_{2}$ if, for some $p_1\in \Lambda_1$ and $p_2\in \Lambda_2$,
\begin{itemize}
\item $r\in W^{u}(p_{1})\cap W^{s}(p_{2})\setminus\{p_{1},p_{2}\}$,
\item $\dim(T_{r}W^{u}(p_{1})+T_{r}W^{s}(p_{2}))=2$.
\end{itemize}
Here, $W^{u}(p_1)$, $W^s(p_2)$ denote the maximal smooth curves in $W^u(\Lambda_1)$, $W^s(\Lambda_2)$ passing through $p_1$ and $p_2$, respectively.
We also say that $\varphi$ has a \textit{tangency}  $r$ of order $n$ 
associated with basic sets $\Lambda_{1}$ and $\Lambda_{2}$ if, for some $p_1\in \Lambda_1$ and $p_2\in \Lambda_2$,
\begin{itemize}
\item $r\in W^{u}(p_{1})\cap W^{s}(p_{2})\setminus\{p_{1},p_{2} \}$, 
\item $\dim(T_{r}W^{u}(p_{1})+T_{r}W^{s}(p_{2}))=1$,
\item there exists a local coordinate $(x, y)$ in a neighborhood of $r$ such that 
$r = (0,0)$, 
$W^{s}(p_{1}) = \{(x, y) : y=0\}$ and 
$W^{u}(p_{2}) = \{(x, y) :  y=u(x)\}$, 
where $u$ is a $C^{n+1}$-function satisfying
\begin{equation}\label{eqn_tangency}
u(0)=u^{\prime}(0)=\cdots=u^{(n)}(0)=0\quad\mbox{and}\quad
u^{(n+1)}(0)\neq 0.
\end{equation}
\end{itemize}
The transverse point or the tangency is called to be \textit{homoclinic} if $\Lambda_{1}=\Lambda_{2}$. 
Usually, the first order tangency is called to \textit{quadratic}, and the second order is \textit{cubic}. 
In particular, the tangency $r$ is quadratic if and only if $W^u(p_1)$ and $W^s(p_2)$ have distinct curvatures at $r$.

\subsection{Generic cubic tangencies}\label{initial_set}
We say that a saddle fixed point $p$ of a $C^\infty$-diffeomorphism $\varphi$ on $\mathbb{R}^2$ is {\it dissipative} if the eigenvalues $\lambda$, $\sigma$ of the differential $(d\varphi)_p$ satisfy $|\lambda\sigma|<1$.
If necessary replacing $\varphi$ by $\varphi^2$, we may assume that such eigenvalues satisfy
\begin{equation}\label{lamda_sigma}
0<\lambda<1<\sigma\quad\mbox{and}\quad\lambda\sigma<1.
\end{equation}

Suppose that $q$ is a cubic homoclinic tangency of $\varphi$ associated with $p$.
The cubic tangency is said to \textit{unfold generically} or to be \textit{generic} with respect to a two-parameter family $\{\varphi_{\mu,\nu}\}$ in $\mathrm{Diff}^{\infty}(\mathbb{R}^2)$ with $\varphi_{0,0}=\varphi$ 
if there exist $(\mu,\nu)$-dependent local coordinates $(x, y)$ on a neighborhood of $q$ with $q = (0,0)$ such that $W^{s}(p_{\mu,\nu}) = \{(x, y) : y=0\}$ and 
$W^{u}(p_{\mu,\nu}) = \{(x, y) :  y=u_{\mu,\nu}(x)\}$, where $\{p_{\mu,\nu}\}$ is the continuation of $p$ consisting of saddle fixed points of $\varphi_{\mu,\nu}$ and $u_{\mu,\nu}(x)=u(\mu,\nu,x)$ is a $C^\infty$-function satisfying
\begin{equation}\label{generic condition}
(\partial_{\mu}u \cdot \partial_{\nu x}u-\partial_{\nu}u\cdot  \partial_{\mu x}u)(0,0,0)\neq 0.
\end{equation}

In particular, $W^u(p_{\mu,\nu})$ is locally parametrized by the curve $c(t)=(t,u_{\mu,\nu}(t))$ in a small neighborhood of $q$.
Since its tangent vector $c'(t)=(1,u'_{\mu,\nu}(t))$ is non-zero for any parameter values $t$, $c(t)$ is a regular curve.

\begin{figure}[hbt]
\begin{center}
\scalebox{0.85}{\includegraphics[clip]{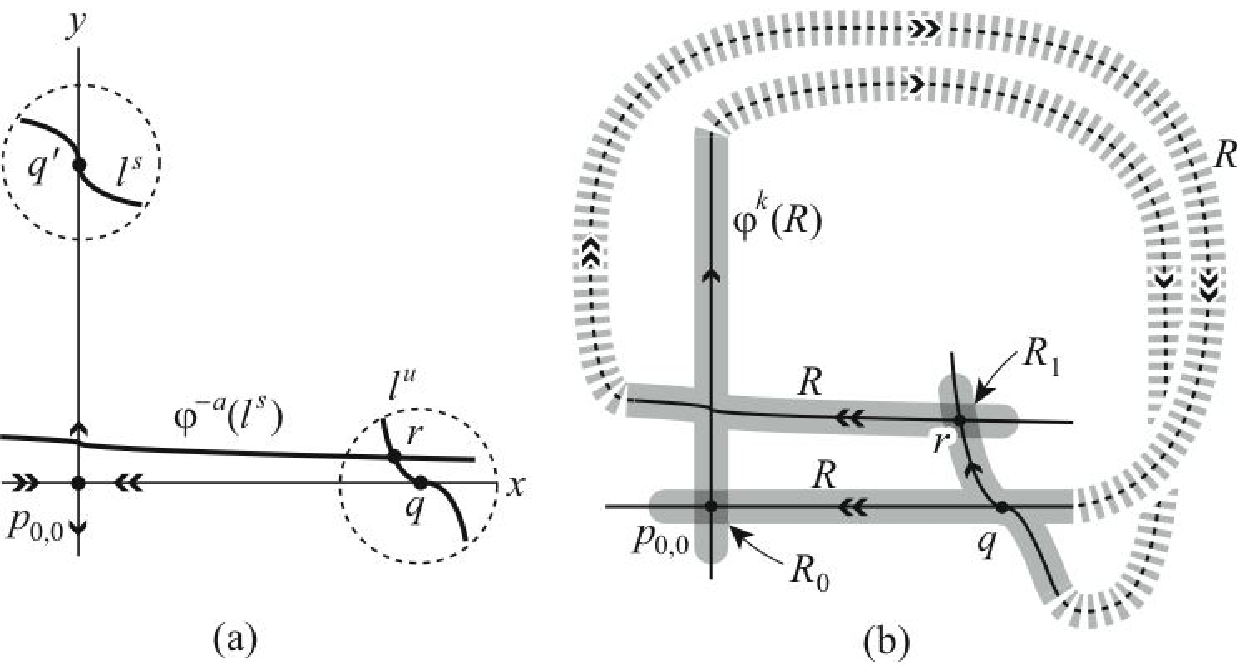}}
\caption{}
\label{fg_1_1}
\end{center}
\end{figure}

\begin{setting}\label{set}
Throughout the remainder of this paper, let $\{\varphi_{\mu,\nu}\}$ represent a two-parameter family in $\mathrm{Diff}^\infty(\mathbb{R}^2)$ given in Theorem \ref{main_a}.
We may assume that $\varphi_{0,0}=\varphi$ has the origin $p_{0,0}=(0,0)$ 
as a dissipative saddle fixed point 
and admitting a cubic homoclinic tangency $q$ which unfolds generically.
Since $\varphi_{0,0}$ satisfies the open $C^4$-linearizing condition at $(0,0)$, there exists a small open neighborhood $\mathcal{O}_0$ of $(0,0)$ in the $\mu\nu$-plane, such that, for any $(\mu,\nu)\in \mathcal{O}_0$, 
$\varphi_{\mu,\nu}$ is $C^4$-linearizable in a small neighborhood $U$ of $(0,0)$ in $\mathbb{R}^2$.
If necessary retaking the $C^4$-linearizing coordinate on $U$, 
one can suppose without loss of generality that it satisfies the following conditions:
\begin{itemize}
\item
$p_{\mu,\nu}=(0, 0)$ and the restriction 
$\varphi_{\mu,\nu}|_{U}$ has the form 
$$\varphi_{\mu,\nu}(x,y)=(\lambda_{\mu,\nu} x, \sigma_{\mu,\nu} y)$$
for any $(x,y)\in U$ where 
$\lambda_{\mu,\nu}$ and $\sigma_{\mu,\nu}$ are the contracting and expanding eigenvalues of $d \varphi_{\mu,\nu}(p_{\mu,\nu})$, respectively.
\item
$\{ (x, y) : |x|\leq 2,  |y|\leq 2\}\subset U$ and  $q=(1, 0)$.
\item
$q^{\prime}=(0, 1)$ is a point with $\varphi^{N}(q^{\prime})=q$ for some integer $N>0$, see Fig.\ \ref{fg_1_1}\,(a).
\end{itemize}
\end{setting}

\subsection{Thickness of Cantor sets}\label{thickness}
We here recall the definition of the thickness given in Newhouse \cite{N1} and Palis-Takens \cite{PT} for a Cantor set $K$ contained in an interval $I$.
A \textit{gap} of $K$ is a connected component of $I\setminus K$ which does not contain a boundary point of $I$. 
Let $G$ be a gap and $p$ a boundary point of $G$.
A closed interval $B\subset  I$ is called the \textit{bridge} at $p$ if $B$ is the maximal interval with $G\cap B=\{p\}$ such that $B$ does not intersect any gap whose length is at least that of $G$. 
The  \textit{thickness}  of  $K$ at $p$ is defined by  $\tau(K,p)= \mathrm{Length}(B)/\mathrm{Length}(G)$.
The \textit{thickness} $\tau(K)$ of $K$ is the infimum  over these $\tau(K,p)$ for all boundary points $p$ of gaps of $K$.
Let $K_1,K_2$ be two Cantor sets in $I$ with thickness $\tau_1$ and $\tau_2$ respectively.
Then, Gap Lemma (see \cite[\S 4]{N1}, \cite[\S 4.2]{PT}) shows that, if $\tau_1\cdot \tau_2>1$, then either $K_1$ is contained in a gap of $K_2$, or $K_2$ is contained in a gap of $K_1$, or $K_1\cap K_2\neq \emptyset$.

The thickness of a Cantor subset $K$ of a $C^1$-curve $\alpha$ in $\mathbb{R}^2$ is defined similarly by supposing that $\alpha$ is parametrized by arc-length.
Properties of the thickness are explained in following subsections.

\subsection{Dynamically defined Cantor sets}\label{dynam_Cantor}

Let $K$ be a Cantor set in $\mathbb{R}$ represented as $K=\bigcap_{i=0}^\infty \psi^{-i}(K_1\cup\cdots\cup K_n)$ for a Markov partition $\{K_1,\dots,K_n\}$ in the sense of \cite[\S 4.1]{PT} with respect to an expanding map $\psi:K_1\cup\cdots\cup K_n\longrightarrow \mathbb{R}$.
We say that such a $K$ is a \textit{dynamically defined Cantor set} if $K$ and $\psi$ satisfy the following conditions.
\begin{itemize}
\item
For some real number $\lambda$ with $0<\vert\lambda\vert <1$, $\lambda K$ is a neighborhood of $0$ in $K$.
\item
The map $\psi$ has a $C^{1+\varepsilon}$-expansive extension to a neighborhood of $K$.
\end{itemize}
By \cite[\S 4.2, Proposition 7]{PT}, the thickness $\tau(K)$ of such a Cantor set is positive.
Moreover, by \cite[\S 4.3, Theorem 2]{PT}, $\tau(K)$ depends continuously on $K$.

Let $\varphi$ be a $C^3$-diffeomorphism of $\mathbb{R}^2$ and $\Lambda$ a non-trivial basic set of $\varphi$.
If $\varphi$ has a 
saddle fixed point $p$ in $\Lambda$ which is linearizable, 
there exists a smooth identification $\alpha:\mathbb{R}\longrightarrow W^s(p)$ such that $\alpha^{-1}\circ (\varphi|W^s(p))\circ \alpha$ is a linear contraction.
Then, there exists a dynamically defined Cantor set $K$ in $\mathbb{R}$ such that $\alpha(K)=W^u_{\mathrm{loc}}(\Lambda)\cap W_{\mathrm{loc}}^s(p)$.
The image $W^u_{\mathrm{loc}}(\Lambda)\cap W_{\mathrm{loc}}^s(p)$ is also called a dynamically defined Cantor set.
Here, we suppose that any local unstable manifold $W^u_{\mathrm{loc}}(\Lambda)$ of $\Lambda$ considered in this paper consists of curves, called {\it leaves} of $W^u_{\mathrm{loc}}(\Lambda)$, meeting $W^s_{\mathrm{loc}}(p)$ transversely in a single point.

We refer to \cite[Chapter 4]{PT} for details on dynamically defined Cantor sets.

\subsection{Outside basic set and  its thickness}\label{outside}

From Setting \ref{set}, $q'$ is also a cubic tangency of $\varphi$ unfolded generically.
The unstable manifold $W^{u}(p_{0,0})$ contains a regular curve $l^u$ cubically tangent to $W^{s}(p_{0,0})$ at $q$.
Similarly, $W^s(p_{0,0})$ contains a regular curve $l^s$ cubically tangent to $W^{u}(p_{0,0})$ at $q'$.


\begin{lemma}\label{lem_2_2}
For a sufficiently large integer $a>0$, $\varphi^{-a}(l^{s})$ and $l^{u}$ have a transverse intersection $r$.
\end{lemma}
\begin{proof}
By the intermediate value theorem, for any sufficiently large integer $a>0$, $\varphi^{-a}(l^{s})\cap l^{u}$ has a point $r=r(a)$ in a small neighborhood of $q$.
Let $s_1,s_2$ be the slopes of $\varphi^{-a}(l^s)$ and $l^u$ at $r$ respectively.
Then, it is not hard to see that
$$|s_1|\leq c_1 \sigma^{-a}\quad\mbox{and}\quad |s_2|\geq c_2\sigma^{-2a/3}$$
for some constants $c_1,c_2>0$ independent of $a$.
Thus, one can take $a$ so that $|s_1|/|s_2|<1$ holds, and hence in particular $\varphi^{-a}(l^{s})$ meets $l^{u}$ transversely at $r$.
\end{proof}

Let $\rho$ be the arc in $W^s(p_{0,0})$ with end points $p_{0,0}$ and $r$.
Suppose that $R$ is a thin tubular neighborhood of $\rho$ in $\mathbb{R}^2$ as illustrated in Figure\ \ref{fg_1_1}\,(b) such that two of the components of $R\cap \varphi^k(R)$ are rectangles $R_0,R_1$ with $\mathrm{Int}(R_0)\ni p_{0,0}$, $\mathrm{Int}(R_1)\ni r$ for some integer $k>0$.
One can choose such an $R$ so that
$${\Lambda_{0,0}^{\mathrm{out}}}=\bigcap_{n=-\infty}^\infty \varphi^{kn}(R_0\cup R_1)$$
is a basic set of $\varphi^k$.
Then, for all $(\mu,\nu)$ sufficiently near $(0,0)$, there exists the continuation $\{\Lambda_{\mu,\nu}^{\mathrm{out}}\}$ of $\Lambda_{0,0}^{\mathrm{out}}$ consisting of basic sets of $\varphi^k_{\mu,\nu}$.
Here, the superscript ``out'' implicitly suggests that $\Lambda_{\mu,\nu}^{\mathrm{out}}$ is in the outside of a small neighborhood of $q$.
Replace  $\varphi_{\mu,\nu}^k$ by $\varphi_{\mu,\nu}$ and suppose from now on that $\Lambda^{\mathrm{out}}_{\mu,\nu}$ is a basic set of $\varphi_{\mu,\nu}$.

Since $\Lambda_{\mu,\nu}^{\mathrm{out}}\cap W^s(p_{\mu,\nu})$ is a dynamically defined Cantor set in the sense of 
Subsection \ref{dynam_Cantor}, the following lemma is derived directly from some results in \cite{PT}.

\begin{lemma}\label{lem_2_3}
There exists an open neighborhood $\mathcal{O}_1$ of $(0,0)$ in $\mathcal{O}_0$ and $\alpha>0$ such that, for any $(\mu,\nu)\in \mathcal{O}_1$, the unstable thickness $\tau(\Lambda_{\mu,\nu}^{\mathrm{out}}\cap W^s(p_{\mu,\nu}))$ is greater than $\alpha$, where $\mathcal{O}_0$ is the open neighborhood given in Setting \ref{set}.
\end{lemma}
\begin{proof}
By \cite[\S 4.2, Proposition 7]{PT}, $\tau(\Lambda_{0,0}^{\mathrm{out}}\cap W^s(p_{0,0}))$ is positive.
Since $\tau(\Lambda_{\mu,\nu}^{\mathrm{out}}\cap W^s(p_{\mu,\nu}))$ is continuous on $(\mu,\nu)$ (see \cite[\S 4.3, Theorem 2]{PT}), one can have $\alpha>0$ satisfying our desired properties if the neighborhood $\mathcal{O}_1$ is taken sufficiently small.
\end{proof}

\section{Local representations near cubic tangencies}

In \cite{PT}, Palis and Takens showed that there is a renormalization at generically unfolding homoclinic quadratic tangencies 
whose limit dynamics is the logistic map.
Similar results for cubic homoclinic tangencies have been already obtained by several authors, see \cite{G07} and references 
therein.
Especially,  Gonchenko et al \cite{G07} present such a renormalization without relying on the 
linearizability assumption.
However, since they work in general settings, it might be hard for the reader to detect in \cite{G07} 
explicit formulas needed to prove Theorems \ref{main_a} and \ref{main_b}.
So, in this section, we will present formulas of a renormalization at a 
generically unfolding homoclinic cubic tangency which are just suitable for our purpose.
Among others, the equations (\ref{renormal})-(\ref{norm}) below play a crucial role in the proof of Lemma 
\ref{lem_6_5}, which is a key lemma for Theorem \ref{main_b}.
Our proofs are done under the linearizability assumption.
The authors do not know whether one can obtain the same results without the assumption.

\subsection{Renormalizations near cubic tangencies}\label{renorm}
Roughly speaking, the renormalization is a sequence of local representations  of $\varphi_{\mu, \nu}$ near the generic cubic tangency $q$ such that the representations converge to a two-parameter family of cubic endomorphisms.

For any $(x,y)\in U$, set $\varphi_{\mu,\nu}^N(x,y)=(h_1(\mu,\nu,x,y-1),h_2(\mu,\nu,x,y-1))$, where $U$ is the neighborhood of $(0,0)$ and $N$ is the positive integer given in Setting \ref{set}.
Since $\varphi^N(0,1)=(1,0)$ and the map $c:(-\delta,\delta)\longrightarrow \mathbb{R}^2$ defined by $c(y)=\varphi^N(0,y)$ is a regular curve tangent to the $x$-axis at $q=(1,0)$ for a small $\delta >0$,
\begin{equation}\label{h_1h_2}
h_1(\mathbf{0})=1\quad\mbox{and}\quad \partial_y h_1(\mathbf{0})\neq 0,
\end{equation}
where $\mathbf{0}=(0,0,0,0)$ is the origin of $\mathbb{R}^4$.
From the latter condition of (\ref{h_1h_2}), there exists an inverse function $y=\eta(t)$ of $t=h_1(0,0,0,y-1)-1$ for any $t$ near $0$ with $\eta(0)=1$.
Set $u(t)=h_2(0,0,0,\eta(t)-1)$.
Then, any point of the image $c(1-\delta,1+\delta)$ of $c$ is represented as $(t+1,u(t))$.
From the form of $u(t)$,
\begin{align*}
u'(0)&=\partial_y h_2(\mathbf{0})\eta'(0),\\
u''(0)&=\partial_{yy}h_2(\mathbf{0})(\eta'(0))^2+\partial_y h_2(\mathbf{0})\eta''(0),\\
u'''(0)&= \partial_{yyy}h_2(\mathbf{0})(\eta'(0))^3+3\partial_{yy}h_2(\mathbf{0})\eta'(0)\eta''(0)+\partial_y h_2(\mathbf{0})\eta'''(0).
\end{align*}
Since $\eta'(0)=1/\partial_y h_1(\mathbf{0})\neq 0$, the cubically tangential condition $\mbox{(\ref{eqn_tangency})}_{n=2}$ at $q=(1,0)$ is rewritten as
\begin{equation}\label{tangent conditions}
h_2(\mathbf{0})=\partial_y h_2(\mathbf{0})=\partial_{yy}h_2(\mathbf{0})=0\quad\mbox{and}\quad \partial_{yyy}h_2(\mathbf{0})\neq 0.
\end{equation}
Consider the $C^4$-functions $H_1,H_2$ defined by
\begin{align*}
H_1(\mu,\nu,x,y-1)&=h_1(\mu,\nu,x,y-1)-1-a_0(y-1),\\
H_2(\mu,\nu,x,y-1)&=h_2(\mu,\nu,x,y-1)-a_1x-a_2\mu-a_3\nu\\
&\qquad\qquad\qquad -(a_4\mu+a_5\nu)(y-1)-a_6(y-1)^3,
\end{align*}
where $a_i$ $(i=0,1,\ldots,6)$ are the constants defined as
\begin{align*}
a_0&=\partial_y h_1(\mathbf{0}),\ 
a_1=\partial_x h_2(\mathbf{0}),\ 
a_2=\partial_\mu h_2(\mathbf{0}),\ 
a_3=\partial_\nu h_2(\mathbf{0}),\\
a_4&=\partial_{\mu y} h_2(\mathbf{0}),\ 
a_5=\partial_{\nu y} h_2(\mathbf{0}),\ 
a_6=\partial_{yyy} h_2(\mathbf{0})/6.
\end{align*}
Then, the conditions (\ref{h_1h_2}) and (\ref{tangent conditions}) imply that, for $(\mu,\nu,x,y-1)=\mathbf{0}$,
\begin{equation}\label{conditions of hot}
\left\{
\begin{aligned}
H_{1}&=\partial_{y}{H}_{1}=0,\\
H_{2}&=\partial_{\mu}{H}_{2}=\partial_{\nu}{H}_{2}=\partial_{x}{H}_{2}=\partial_{y}{H}_{2}=\partial_{\mu y}{H}_{2}=\partial_{\nu y}{H}_{2}=\partial_{y y}{H}_{2}\\
&=\partial_{yyy}H_{2}=0.
\end{aligned}
\right.
\end {equation}
From the genericity condition (\ref{generic condition}) for $\{\varphi_{\mu,\nu}\}$, we have
$$a_2a_5-a_3a_4\neq 0.$$
Thus, one can consider the new parameters $(\hat \mu,\hat \nu)=(a_4\mu+a_5\nu,a_2\mu+a_3\nu)$.
For saving symbols, we denote $(\hat \mu,\hat \nu)$ again by $(\mu,\nu)$.
Since the coordinate change is linear, the equalities in (\ref{conditions of hot}) still hold with respect to the new coordinate.
Thus, $\varphi_{\mu,\nu}^N$ is represented in a neighborhood of $q'=(0,1)$ as
\begin{equation}\label{form}
\begin{aligned}
\varphi_{\mu,\nu}^{N}(x, y)&=\big(1+a(y-1)+H_{1}(\mu,\nu,x,y-1),\\ 
& \qquad -b(y-1)^{3}+\mu (y-1)+\nu+ cx+H_{2}(\mu,\nu,x,y-1) \big),
\end{aligned}
\end{equation}
where $a=a_0$, $b=-a_6$ and $c=a_1$.

The latter condition of (\ref{h_1h_2}) implies that $a\neq 0$ and that of (\ref{tangent conditions}) does $b\neq 0$.
From now on, we only consider the case of $b>0$, which is not an essential restriction in our arguments.

\begin{lemma}\label{lem_3_1}
Under the conditions in Setting \ref{set}, for any sufficiently large $n>0$, there exists a reparametrization $\Theta_{n}$ on $\bar{\Sigma}=[0,4]\times [-1,1]$ and a coordinate change $\Phi_{n}$ on $\mathbb{R}^2$ satisfying following conditions.
\begin{enumerate}[\rm (i)]
\item
The domains of $\Phi_n$ are open sets of $\mathbb{R}^2$ eventually containing any compact subset of $\mathbb{R}^2$.
\item  
For any $(\bar{\mu}, \bar{\nu},\bar x,\bar y) \in \bar{\Sigma}\times \mathbb{R}^2$, 
$$\lim_{n\rightarrow\infty} \Theta_{n}(\bar{\mu}, \bar{\nu})=(0, 0),\quad 
 \lim_{n\rightarrow\infty} \Phi_{n} (\bar{x}, \bar{y})=q.$$
\item 
For any $(\bar{\mu}, \bar{\nu}) \in \bar{\Sigma}$, the diffeomorphisms on $\mathbb{R}^2$ defined by
$$(\bar{x},\bar{y})\mapsto
  \Phi_{n}^{-1}\circ\varphi_{\Theta_{n}(\bar\mu,\bar\nu)}^{N}\circ  \varphi_{\Theta_{n}(\bar\mu,\bar\nu)}^{n}\circ \Phi_{n}(\bar{x},\bar{y})
 $$ 
$C^3$-converge as $n\rightarrow \infty$ locally uniformly to the cubic endomorphism
$$
(\bar{\mu},  \bar{\nu},  \bar{x}, \bar{y})\mapsto
(\bar{\mu},  \bar{\nu},  \psi_{\bar{\mu},  \bar{\nu}}( \bar{x},   \bar{y})),
$$
where
$\psi_{\bar{\mu},  \bar{\nu}}( \bar{x},   \bar{y})=( \bar{y},   - \bar{y}^{3}+\bar{\mu}  \bar{y}+\bar{ \nu}).$
\end{enumerate} 
\end{lemma}
\begin{figure}[hbt]
\begin{center}
\includegraphics{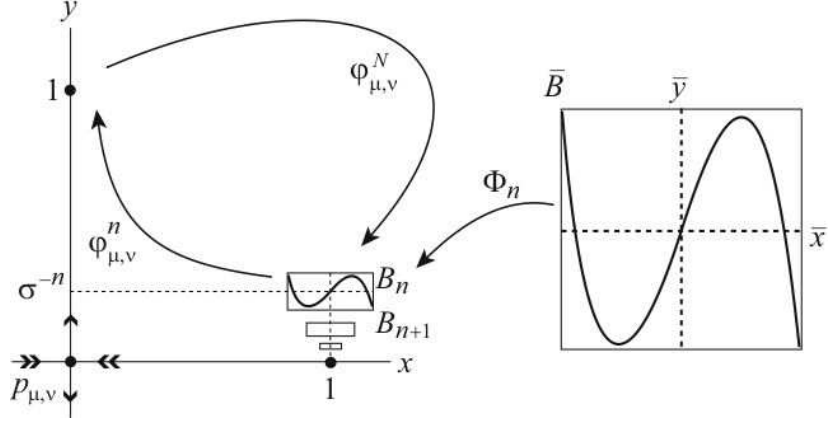}
\caption{ $\bar B=[-2,2]^2$ and $B_n=\Phi_n(\bar B)$.}
\label{fg_2_1}
\end{center}
\end{figure}
\begin{proof}
As in Subsection \ref{initial_set}, the contracting and expanding eigenvalues of 
 $d \varphi_{\mu,\nu}(p_{\mu,\nu})$ for $(\mu,\nu)$ near $(0, 0)$ are denoted by $\lambda=\lambda_{\mu,\nu}$ and  $\sigma=\sigma_{\mu,\nu}$.
For any  $n>0$, define the $\mu,\nu$-dependent coordinate changes $\Phi_{n}$ on $\mathbb{R}^2$ by
\begin{equation}\label{coordinates change}
(x,y)=\Phi_{n}(\bar{x}, \bar{y})=
(1+ag^{-1}\sigma^{-n/2}\bar{x},\ \sigma^{-n}+g^{-1}\sigma^{-3n/2}\bar{y}),
\end{equation}
where $g=\sqrt{b}$.
From the definition, we have obviously $\lim_{n\rightarrow \infty}\Phi_n(\bar x,\bar y)=q$, see Fig.\ \ref{fg_2_1}.
We define the reparametrization of $\bar \Sigma$ by
\begin{equation}\label{reparameterization}
\bar{\mu}=\sigma^{n}\mu,\quad 
\bar{\nu}=g\sigma^{n/2}(\sigma^{n}\nu+ c(\lambda\sigma)^{n} -1),
\end{equation}
or equivalently by
$$\mu=\sigma^{-n}\bar \mu,\quad \nu=g^{-1}\sigma^{-3n/2}\bar \nu-c\lambda^n+\sigma^{-n}.$$
Let $\Theta$ be the diffeomorphisms on $\bar\Sigma$ given by $\Theta_n(\bar \mu,\bar \nu)=(\mu,\nu)$.
Then, the image $\Theta_{n}(\bar{\Sigma})$ converges to $(0, 0)$ as $n\rightarrow \infty$.

The form of $\varphi_{\mu, \nu}^n$ in the $C^4$-linearizing neighborhood of the origin of $\mathbb{R}^2$ is represented as
$$\varphi_{\mu, \nu}^{n}(x,y)=(\lambda^{n}x, \sigma^{n}y).$$
From the form (\ref{form}) of $\varphi_{\mu, \nu}^{N}$ near  $q^{\prime}$, one can show that $(\bar{x}^{\prime}, \bar{y}^{\prime})=   
\Phi_{n}^{-1}\circ \varphi_{\mu, \nu}^{N}\circ \varphi_{\mu, \nu}^{n} \circ \Phi_{n}(\bar{x}, \bar{y})$ is represented as
\begin{align*}
\bar{x}^{\prime} &= \bar{y}+a^{-1}g\sigma^{n/2}H_{1}(\mu, \nu, x, y-1),\\
 \bar{y}^{\prime} &= - \bar{y}^{3}+\sigma^{n}\mu\bar{y}
+g\sigma^{n/2}(\sigma^{n}\nu+  c(\lambda\sigma)^{n} -1)\\
& \qquad\qquad\qquad\qquad +ac(\lambda\sigma)^{n}\bar{x}
+g\sigma^{3n/2} H_{2}(\mu, \nu, x, y-1)\\
&= -\bar{y}^{3}+\bar{\mu} \bar{y}+\bar{\nu}+ac(\lambda\sigma)^{n}\bar{x}+g\sigma^{3n/2} H_{2}(\mu, \nu,  x, y-1 ),
\end{align*}
where 
\begin{equation}\label{coordinates change2}
( x, y)=\varphi^{n}_{\mu,\nu}\circ \Phi_{n}(\bar{x}, \bar{y})=(\lambda^{n}+ag^{-1}\sigma^{-n/2}\lambda^{n} \bar{x} , 1+g^{-1}\sigma^{-n/2}\bar{y}).
\end{equation}
Set
\begin{align*}
\bar H_{n;1}(\bar\mu,\bar\nu,\bar x,\bar y)&=a^{-1}g\sigma^{n/2}H_{1}(\mu, \nu, x, y-1),\\
\bar H_{n;2}(\bar\mu,\bar\nu,\bar x,\bar y)&=ac(\lambda\sigma)^{n}\bar{x}+g\sigma^{3n/2} H_{2}(\mu, \nu,  x, y-1 ),
\end{align*}
where each of $\mu,\nu$ (resp.\ $x,y$) is considered as a function of $\bar \mu,\bar \nu$ (resp.\ $\bar \mu,\bar \nu,\bar x,\bar y$).
Then, we have
\begin{equation}\label{renormal}
\begin{aligned}
\Phi_{n}^{-1}\circ\varphi_{\Theta_{n}(\bar\mu,\bar\nu)}^{N+n}\circ \Phi_{n}(\bar{x},\bar{y})= \bigl(\bar y+&\bar H_{n;1}(\bar \mu,\bar \nu,\bar x, \bar y),\\
&-\bar y^3+\bar \mu \bar y +\bar \nu +\bar H_{n;2}(\bar \mu,\bar \nu,\bar x, \bar y)\bigr).
\end{aligned}
\end{equation}

It remains to show that $\bar H_{n;1}$ and $\bar H_{n:2}$ $C^3$-converge uniformly to zero on $\bar \Sigma\times K$ for any compact subset $K$ in $\mathbb{R}^2$.
It is not hard to show that the following holds on $\bar \Sigma\times K$ from (\ref{reparameterization}) and (\ref{coordinates change2}),
\begin{equation}\label{eight_eqn}
\begin{split}
&x=O(\lambda^n),\ y-1=O(\sigma^{-n/2}),\ \mu=O(\sigma^{-n}),\ \nu=O(\sigma^{-n}),\\
&\frac{\partial \mu}{\partial \bar \alpha}=O(\sigma^{-n}),\ 
\frac{\partial \nu}{\partial \bar \alpha}=O(\sigma^{-3n/2}),\ 
\frac{\partial x}{\partial \bar \beta}=O(\lambda^n\sigma^{-n/2}),\ 
\frac{\partial y}{\partial \bar \beta}=O(\sigma^{-n/2}),
\end{split}
\end{equation}
where $\bar \alpha$ is either $\bar \mu$ or $\bar \nu$, and $\bar \beta$ is either $\bar \mu,\bar \nu,\bar x$ or $\bar y$.
Using these equalities, one can prove our desired uniform $C^3$-convergence.

To state it more precisely, we note first that the closure $L$ of $\bigcup_{n\geq 0}\Theta_n(\bar \Sigma)$ in $\mathbb{R}^2$ is compact.
In fact, it is shown by the condition (ii) of Lemma \ref{lem_3_1}.
Thus, we have the constants
$$\sigma_{\bar \Sigma}^{-1}=\max\{\sigma_{\mu,\nu}^{-1};\,(\mu,\nu)\in L\}<1\quad \mbox{and}\quad (\sigma\lambda)_{\bar \Sigma}=\max\{\sigma_{\mu,\nu}\lambda_{\mu,\nu};\,(\mu,\nu)\in L\}<1$$
independent of $\mu,\nu$.
Then, the equations in (\ref{eight_eqn}) imply
\begin{equation}\label{norm}
\Vert \bar H_{n;1}\Vert_{C^3,\bar \Sigma\times K}=O(\sigma_{\bar \Sigma}^{-n/2}),\quad 
\Vert \bar H_{n;2}\Vert_{C^3,\bar \Sigma\times K}=O(\sigma_{\bar \Sigma}^{-n/2})+O((\sigma\lambda)_{\bar \Sigma}^n).
\end{equation}
For example, Taylor's expansion of $H_2$ up to the 4th order derivatives together with (\ref{conditions of hot}) implies
\begin{align*}
\bar H_{n;2}&=O(\lambda^n\sigma^n)+\sigma^{3n/2}\bigl(O(|\mu|+|\nu|+|x|)^2\\
&\qquad\qquad\qquad+O(x(y-1))+O(|\mu|+|\nu|)O(y-1)^2\bigr)\\
&=O(\lambda^n\sigma^n)+\sigma^{3n/2}\bigl(O(\sigma^{-2n})+O(\sigma^{-n/2}\lambda^n)+O(\sigma^{-2n})\bigr)\\
&=O(\lambda^n\sigma^n)+O(\sigma^{-n/2}).
\end{align*}
Since $\partial_x H_2=O(\sigma^{-n/2})$ and $\partial_y H_2=O(\lambda^n)+O(\sigma^{-3n/2})$, it follows that
$$
\partial_{\bar y}\bar H_{n;2}=\sigma^{3n/2}\left(\partial_x H_2\frac{\partial x}{\partial \bar y}+\partial_y H_2\frac{\partial y}{\partial \bar y}\right)=O(\lambda^n\sigma^n)+O(\sigma^{-n/2}).
$$
The remaining cases for derivatives on $\bar H_{n;1}$, $\bar H_{n;2}$ up to the third order are verified similarly.
\end{proof}

\section{Dynamically defined Cantor sets for one-dimensional cubic maps}\label{one-dim_cubic_maps}

\subsection{Affine Cantor sets for the N-map}\label{affine}
We consider the \textit{N-map} $S$ on $[-3/2,3/2]$ defined as
 $$
S(x)=\left\{
	\begin{array}{ll}
		-3x-3 & \mathrm{for}\   -{3}/{2}\leq x<  -{1}/{2} \\
		3x &   \mathrm{for}\  -{1}/{2}\leq x\leq {1}/{2}\\
		-3x+3 &  \mathrm{for}\  {1}/{2}< x\leq {3}/{2},
	\end{array}\right.
$$
and estimate the thickness of subsets $K_m$ of $[-3/2,3/2]$ given as follows, 
which are dynamically defined Cantor sets (Subsection \ref{dynam_Cantor}) and do not contain the turning points $\pm 1/2$ of $S$.

The following lemma is an analogous to the estimate for the tent map in \cite[\S 6.2]{PT}.

\begin{lemma}\label{lem_4_1}
For any even number $m\geq 6$, 
there exists a dynamically defined  affine Cantor set $K_{m}$
and an $m$-periodic point $q=q(m)$ for  $S$ satisfying the following conditions.
\begin{enumerate}[\rm (i)]
\item $K_{m}$ does not have the turning points $\pm 1/2$ of $S$.
\item The end points of any gap of $K_{m}$ are eventually mapped to $q$.
\item There exists a compact neighborhood $U'=U_m'$ of $K_m$ in $(-3/2,3/2)$ with $U'\cap \{-1/2,1/2\}=\emptyset$ and $\Vert dS_x\Vert =3$ for any $x\in U'$.
\item $\lim_{m(\mathrm{even})\rightarrow \infty} \tau(K_{m})=\infty$.
\end{enumerate}
\end{lemma}
\begin{proof}
Let
$S_{1}=S\vert_{[-{3}/{2}, -{1}/{2})}$, 
$S_{2}=S\vert_{[-{1}/{2}, {1}/{2}]}$ 
and
$S_{3}=S\vert_{({1}/{2}, {3}/{2}]}$. 
For any even number $m\geq 6$, we will show that there exists an $m$-periodic point $q_{0}=q_{0}(m)$ such that $0<q_{0}<1/2$ and 
$$
S_{2}\circ S_{1}\circ (S_{3}\circ S_{1})^{m/2-2}\circ S_{3}\circ S_{2}
(q_{0})=q_{0}.
$$
Set $S^{(m)}=S_{2}\circ S_{1}\circ (S_{3}\circ S_{1})^{m/2-2}\circ S_{3}\circ S_{2}$.
Then, the number $x_m=(1-1/3^{m-2})/2$ is a zero-point of $S^{(m)}$.
Since $S^{(m)}(1/2)=3/2$, by the Intermediate Value Theorem, there exists $q_0\in (x_m,1/2)$ with $S^{(m)}(q_0)=q_0$.

The orbit of $q_0$ is denoted by $\{q_{i}\}_{i=0,\cdots,m-1}$, that is, $q_{i+1}=S(q_i)$, see Fig.\ \ref{fg_3_1}. 
\begin{figure}[hbt]
\begin{center}
\includegraphics{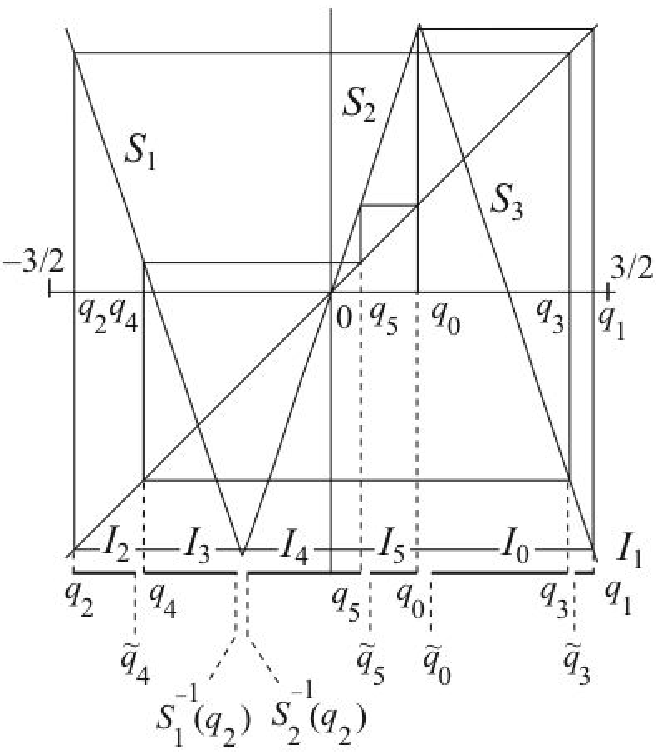}
\caption{The case of $m=6$.}
\label{fg_3_1}
\end{center}
\end{figure}

To construct our desired Cantor set, we need some backward points
defined as 
\begin{align*}
& \tilde{q}_{0}=S_{3}^{-1}(q_{1}),\quad \tilde{q}_{m-1}=S_{2}^{-1}(\tilde{q}_{0}),\\
& \tilde{q}_{m-j-2}=  \left\{
 \begin{array}{l}  S_{1}^{-1}(  \tilde{q}_{m-j-1}) \quad  \mathrm{for}\  j=0, 2,\cdots, m-6,  \\ 
 S_{3}^{-1}(  \tilde{q}_{m-j-1}) \quad  \mathrm{for}\  j=1,3,\cdots, m-5.
\end{array}\right.\\
\end{align*}
These points are mapped eventually to the $m$-periodic orbit of $q_0$ by the forward iterations.
Observe that these (eventually or proper) periodic  points are ordered in $[-3/2, 3/2]$ as follows: 
\begin{align*}
-3/2<q_{2}<\tilde{q}_{4}&<q_{4}<\cdots<\tilde{q}_{m-2}<q_{m-2}<0<q_{m-1}<\tilde{q}_{m-1}<q_{0}\\
&<1/2<\tilde{q}_{0}<q_{m-3}<\tilde{q}_{m-3}<\cdots<q_{3}<\tilde{q}_{3}<q_{1}<3/2.
\end{align*}
Using these points, we define the intervals:
\begin{itemize}
\item $I_{0}=[\tilde{q}_{0}, q_{m-3}], $
\item $ I_{2i-1}= [\tilde{q}_{2i+1}, q_{2i-1}]$ and $I_{2i}= [q_{2i}, \tilde{q}_{2i+2}]$ for $i= 1,\cdots, m/2-2$,
\item $I_{m-3}=[q_{m-2}, S^{-1}_{1}(q_{2})]$ and $I_{m-2}=[S^{-1}_{2}(q_{2}), q_{m-1}]$,
\item  $I_{m-1}=[\tilde{q}_{m-1}, q_{0}]$.
\end{itemize}
For any integer $n>0$, let
$$K_{m}^{(1)}:=\bigcup_{i=0}^{m-1}I_{i},\quad 
K_{m}^{(n)}:=\Bigl(\bigcap_{i=0}^{n-1} S^{-i}(K_{m}^{(1)})\Bigr)\cap [q_2,q_1].
$$
Since $S|_{K_{m}^{(n)}}$ 
maps $K_{m}^{(n)}$ affinely onto $K_{m}^{(n-1)}$,  
the limit 
$K_{m}:=\lim_{n\rightarrow\infty} K_{m}^{(n)}$ is  a 
dynamically defined affine Cantor set  in the sense of Subsection \ref{dynam_Cantor}.
From our definition of $K_m$, one can check the assertions (i)-(iii) easily.

Due to the affine construction of $\{K_{m}^{(n)}\}_{n\geq 1}$ as above, for any $n>0$, the relation between the bridges and the gaps in the $n$-generation $K^{(n)}_m$ is the same as that in the first generation $K^{(1)}_m$ up to scale.
Hence, for estimating the thickness of $K_{m}$, 
it suffices to consider the ratios of the lengths of the bridges and the gaps 
in the first generation. 
We denote here the length of the gap between $\tilde{q}_0$ and $q_{0}$ by $\delta$. 
Then, the length of the bridge at  $\tilde{q}_{0}$ is $1-2\delta$, and that at $q_{0}$ is $2-5\delta$.
This implies that $\tau(K_{m})=(1-2\delta)/\delta$.  
Since the points $\tilde{q}_0$ and $q_{0}$ depends on the even number $m$, $\delta$ also does.
In fact, since $q_{m-2}=(S_3\circ S_1)^{m/2-2}\circ S_3\circ S_2(q_0)$ and $q_0=1/2-\delta$, the direct calculation shows
$$
q_{m-2}=-3^{m-2}\Bigl(\frac{1}{2}-\delta\Bigr)+\frac{3}{2}(3^{m-3}-1).
$$ 
Moreover,  since $S_2\circ S_1(q_{m-2})=q_0$, 
  $$
  \delta= \frac{22}{3^m-1}.
  $$
So,  for any $m\geq 6$,  we obtain $$
\tau(K_m)\geq \frac{3^m-45}{22}.
 $$
This proves the assertion (iv).
\end{proof}

\subsection{Dynamically defined Cantor sets for  cubic maps} 

In this subsection, we show that there is a conjugation between the N-map $S$ given in the previous subsection and the limit map given by the renormalization procedure from $\{\varphi_{\mu,\nu}\}$.
Our argument is an analogy of that in \cite[Chapter 6.2]{PT}, where the tent map is conjugated to the logistic map obtained  by the renormalization associated to a quadratic tangency.

Let us consider  the family of the one-dimensional cubic maps defined  as  
$$
F_{\bar{\mu},\bar{\nu}}(\bar{y})= -\bar{y}^{3}+\bar{\mu}  \bar{y}+ \bar{\nu}.
$$

\begin{proposition}\label{prop_4_2}
For $(\bar{\mu},\bar{\nu})=(3,0)$, the cubic map 
$F=F_{3,0}$ on  $[-2, 2]$ 
is  conjugate to 
 the N-map $S$ on $[-3/2, 3/2]$.
 \end{proposition}
\begin{proof}
Consider the homeomorphism $h:[-{3}/{2}, {3}/{2}]\rightarrow[-2,2]$ given by $h(x)=2\sin ({\pi x}/{3})$.
From the formula for triple angle, we have
$$F\circ h(x)=-8\sin^3(\pi x/3)+6\sin (\pi x/3)=2\sin(\pi x).$$
On the other hand, it is not hard to check that $h\circ S(x)=2\sin (\pi x)$ for any $x\in [-3/2,3/2]$.
This shows $h\circ S=F\circ h$.
\end{proof}

\begin{proposition}\label{prop_4_3}
For any even number $m\geq 6$, 
there exists a dynamically defined Cantor set $\bar{K}_{m}\subset [-2, 2]$ and an $m$-periodic point $\bar{q}=\bar{q}(m)$ for $F$ such that 
\begin{enumerate}[\rm (i)]
\item $\bar{K}_{m}=\mathrm{Cl}( \bigcup_{i\geq 0} F^{-i}(\bar{q}))$ has no critical points of $F$;
\item there exists a neighborhood $U$ of $\bar K_m$ in $[-2,2]$ and a constant $c>0$ such that, for any integer $n\geq 0$ and any $\bar y\in [-2,2]$ with $F^i(\bar y)\in U$ $(i=0,1,\dots,n-1)$, $\Vert dF^n(\bar y)\Vert\geq 3^nc$;
\item $\lim_{m(\mathrm{even})\rightarrow \infty} \tau(\bar{K}_{m})=\infty$.
\end{enumerate}
\end{proposition}
\begin{proof}
Define $\bar{K}_{m}=h({K}_{m})$,  where $h$ is the homeomorphism given  in Proposition \ref{prop_4_2}.
The assertion (i) is obvious since the critical points of $F$ are $\pm 1$ and $h^{-1}(\pm 1)=\pm 1/2\not\in K_m$.
Let $U'$ be the compact neighborhood of $K_m$ in $(-3/2,3/2)\setminus \{-1/2,1/2\}$ given in Proposition \ref{lem_4_1}\,(iii).
Then, $U=h(U')$ and the constant $c={c'}^2$ satisfy the conditions of the assertion (ii), where $0<c'<1$ is a constant with $c'\leq \Vert dh_x\Vert \leq {c'}^{-1}$ for any $x\in U'$.

Though the points $x=\pm 3/2$ are not contained in ${K}_{m}$, any neighborhood of these points intersects ${K}_{m}$ non-trivially for all sufficiently large even number $m$.
For a sufficiently small fixed $\delta>0$, $h$ is a diffeomorphism on $I_{\delta}=[-3/2+\delta, 3/2-\delta]$. 
We have the $\delta$-dependent constant 
$$c=\frac{\max \{h^{\prime}(x): x\in I_{\delta} \}}{\min\{h^{\prime}(x): x\in I_{\delta} \}}>0.$$
Note here that $h(x)$ is essentially quadratic near $x=\pm 3/2$.
Thus, one can show that
$$\tau(\bar{K}_{m})\geq \frac{\tau(K_{m})}{\max\{c, 4\}}$$
by using the argument in Palis-Takens \cite[\S 6.2, p.\ 119]{PT}.
This proves (iii).
\end{proof}

\section{Compatible foliations and Accompanying lemma}\label{C_foliations_Acc_lemma}
In this section, we will introduce a powerful scheme which is crucial in the
proof of Theorem \ref{main_a}. A method of estimating the thickness of compatible
stable and unstable foliations along a transverse curve to show that these
foliations have a tangency has been already presented by 
Newhouse,  Palis and Takens
\cite{N1,PT}. We also work along their line, but need a more elaborated method to
prove Theorem \ref{main_a}. Especially, Accompanying Lemma (Lemma \ref{lem_1_1}) given in
Subsection \ref{ac_lemma} is very useful to show that such a tangency is either contact-making
or breaking with respect to a certain one-parameter subfamily.

\subsection{Compatible foliations}\label{c_foliation}
Let $\mathcal{F}$ be a foliation consisting of smooth curves in the plane.
A smooth curve $\sigma$ in the plane is said to {\it cross} $\mathcal{F}$ {\it exactly} if each leaf of $\mathcal{F}$ intersects $\sigma$ transversely in a single point and any point of $\sigma$ is passed through by a leaf of $\mathcal{F}$.

\begin{figure}[hbtp]
\centering
\scalebox{0.6}{\includegraphics[clip]{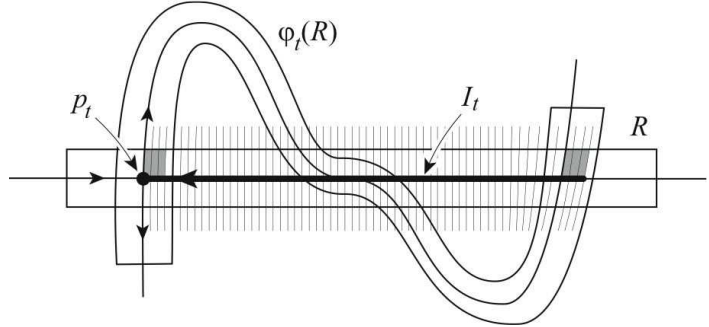}}
\caption{ This is the case where $\Lambda_t$ is given by $\bigcap_{n=-\infty}^\infty \varphi_t^n(R_0\cup R_1)$ for two curvilinear rectangle components $R_0, R_1$ of $R\cap \varphi_t(R)$ for a rectangle $R$ as above.
The shaded regions contain $\Lambda_t$, and the short vertical arcs represent leaves of $\mathcal{F}^u_t$.}
\label{fg_4_1}
\end{figure}

For some $\varepsilon>0$, let $\{\varphi_t\}_{-\varepsilon\leq t\leq \varepsilon}$ be a one-parameter family in $\mathrm{Diff}^3(\mathbb{R}^2)$ and $\{\Lambda_t\}$ a continuation of non-trivial basic sets of $\varphi_t$.
Suppose that $\{p_t\}$ is a continuation of saddle fixed points in $\Lambda_t$ and $\{I_t\}$ is a continuation of curves in $W^s_{\mathrm{loc}}(p_t)$ which are shortest among curves in $W^s_{\mathrm{loc}}(p_t)$ containing $\Lambda_t\cap W^s_{\mathrm{loc}}(p_t)$.
According to Lemma 4.1 in Kan-Ko\c{c}ak-Yorke \cite{KKY} based on results in Franks \cite{Fr}, there exists a continuation of foliations $\mathcal{F}^u_t$ in the plane satisfying the following conditions.
Such foliations are said to be {\it compatible with} $W_{\mathrm{loc}}^u(\Lambda_t)$, see Fig.\ \ref{fg_4_1}.
\begin{enumerate}[\rm (i)]
\item
Each leaf of $W^u_{\mathrm{loc}}(\Lambda_t)$ is a leaf of $\mathcal{F}^u_t$.
\item
$I_t$ crosses $\mathcal{F}^u_t$ exactly.
\item
Leaves of $\mathcal{F}^u_t$ are $C^3$-curves such that themselves, their directions, and their curvatures vary $C^1$ with respect to any transverse direction and $t$.
\end{enumerate}
Similarly, there exist foliations $\mathcal{G}^s_t$ {\it compatible with} $W^s_{\mathrm{loc}}(\Lambda_t)$.
Note that the definition of these foliations is slightly different from that of unstable and stable foliations in a usual sense.

\begin{figure}[hbtp]
\centering
\scalebox{0.6}{\includegraphics[clip]{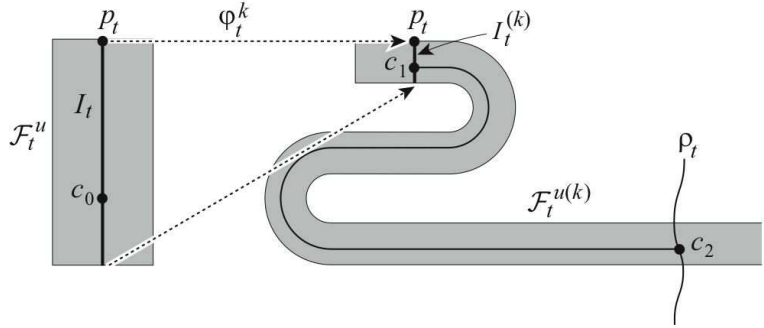}}
\caption{ The case in which $\varphi_t|W^s_{\mathrm{loc}}(p_t)$ is orientation-preserving.
For any point $c_0\in I_t$, $c_1=\varphi_t^k(c_0)\in I_t^{(k)}\subset I_t$ and $c_2=\pi_t^{u(k)}(c_1)=\hat\pi_t^{u(k)}(c_0)\in \rho_t$.}
\label{fg_4_2}
\end{figure}

Let $\rho_t$ be $C^1$-arcs in $\mathbb{R}^2$ $C^1$-varying with respect to $t$ and each of which contains a subarc crossing $\mathcal{F}^{u(k)}_t=\varphi^k_t(\mathcal{F}^u_t)$ exactly for some integer $k\geq 0$.
By the condition (iii), the projection $\pi_t^{u(k)}:I_t^{(k)}=\varphi_t^k(I_t)\longrightarrow  \rho_t$ along leaves of $\mathcal{F}^{u(k)}_t$ is a $C^1$-diffeomorphism $C^1$-depending on $t$.
Note that each leaf of $\varphi^k(W_{\mathrm{loc}}^u(\Lambda_t))$ meets $I^{(k)}_t$ transversely in a single point.
The restriction $\varphi_t^k|I_t:I_t\longrightarrow I_t^{(k)}$ is a $C^3$-diffeomorphism transforming $\Lambda_t\cap I_t$ onto $\Lambda_t\cap I_t^{(k)}$.
The $C^1$-diffeomorphisms $\hat\pi_t^{u(k)}$ defined by
\begin{equation}\label{pi_t}
\hat\pi_t^{u(k)}=\pi_t^{u(k)}\circ \varphi_t^k|I_t:I_t\longrightarrow \rho_t
\end{equation}
$C^1$-depend on $t$, see Fig.\ \ref{fg_4_2}.
Since $\varphi_t^k(W^u_{\mathrm{loc}}(\Lambda_t))\subset \varphi_t^k(W^u(\Lambda_t))=W^u(\Lambda_t)$, for each $c\in \Lambda_t\cap I_t$, $\hat\pi_t^{u(k)}(c)$ is contained in a leaf of $W^u(\Lambda_t)$.
Any leaf of $\mathcal{F}_t^{u(k)}$ meeting $\hat\pi_t^{u(k)}(\Lambda_t\cap I_t)$ non-trivially is called a $\Lambda_t$-{\it leaf}.

For any integer $l\geq 0$ and a continuation of foliations $\mathcal{G}_t^s$ compatible with $W^s_{\mathrm{loc}}(\Lambda_t)$, the $C^1$-diffeomorphism
\begin{equation}\label{pi_t2}
\hat \pi_t^{s(l)}=\pi^{s(l)}\circ \varphi^{-l}_t|J_t:J_t\longrightarrow \xi_t
\end{equation}
is defined similarly, where $J_t$ is a continuation of arcs in $W^u(p_t)$ crossing $\mathcal{G}_t^s$ exactly and $\xi_t$ are $C^1$-arcs $C^1$-depending on $t$ each of which contains a subarc crossing $\mathcal{G}_t^{s(l)}=\varphi_t^{-l}(\mathcal{G}_t^s)$ exactly.

From the definition of thickness, for any bridge $B$ in $I_t$ of the Cantor set $\Lambda_t\cap I_t$, $\tau(\Lambda_t\cap B)\geq \tau(\Lambda_t\cap I_t)$.
Since $\hat\pi_t^{u(k)}$ is a $C^1$-diffeomorphism, $\hat\pi_t^{u(k)}$ is well approximated by an affine map on any sufficiently short subinterval of $I_t$.
Also since the thickness of Cantor sets is invariant under affine transformations, 
for any $\delta >0$, if the bridge $B$ is sufficiently small, then
\begin{equation}\label{thick_B}
\tau(\hat\pi_t^{u(k)}(\Lambda_t\cap B))>\tau(\Lambda_t\cap I_t)-\delta.
\end{equation}

\subsection{Persistent contact-making and breaking tangencies}\label{p_tangency}
Let $\{\varphi_t\}$ be a one-parameter family  of  diffeomorphisms on the plane $\mathbb{R}^2$ with continuations of basic sets $\Lambda_{1,t}, \Lambda_{2,t}$ of $\varphi_t$ (possibly $\Lambda_{1,t}= \Lambda_{2,t}$).
Then, we say that $\{\varphi_t\}$ has a \textit{contact-making} (resp.\ \textit{contact-breaking}) tangency at $t=t_0$ if $\varphi_{t_0}$ has a quadratic tangency $r_{t_0}$ associated with $\Lambda_{1,t_0}$ and $\Lambda_{2,t_0}$ such that transverse points of $\varphi_t$ are generated (resp.\ annihilated) in a small neighborhood $\mathcal{N}(r_{t_0})$ of $r_{t_0}$ in $\mathbb{R}^2$.
More precisely, there exist continuations of curves $l^s_t\subset W^s(\Lambda_{1,t})\cap \mathcal{N}(r_{t_0})$, 
$l^u_t\subset W^u(\Lambda_{2,t})\cap \mathcal{N}(r_{t_0})$ and a sufficiently small $\delta >0$ such that 
\begin{itemize}
\item $l^u_t\cap l^s_t = \varnothing$ for $t<t_0$ (resp.\ $t>t_0$) with $\vert t-t_0\vert <\delta$, 
\item $l^u_{t_0}\cap l^s_{t_0}=\{r_{t_0}\}$, 
\item $l^u_t$ meets $l^s_t$ non-trivially and transversely for $t>t_0$ (resp.\ $t<t_0$) with $\vert t-t_0\vert <\delta$. 
\end{itemize}
The family $\{\varphi_t\}$ is said to exhibit \textit{persistent contact-making} (resp.\ \textit{persistent contact-breaking}) tangencies if each $\varphi_t$ has a contact-making (resp.\ contact-breaking) tangency $r_t$.

\begin{figure}[htb]
\centering
\scalebox{0.61}{\includegraphics[clip]{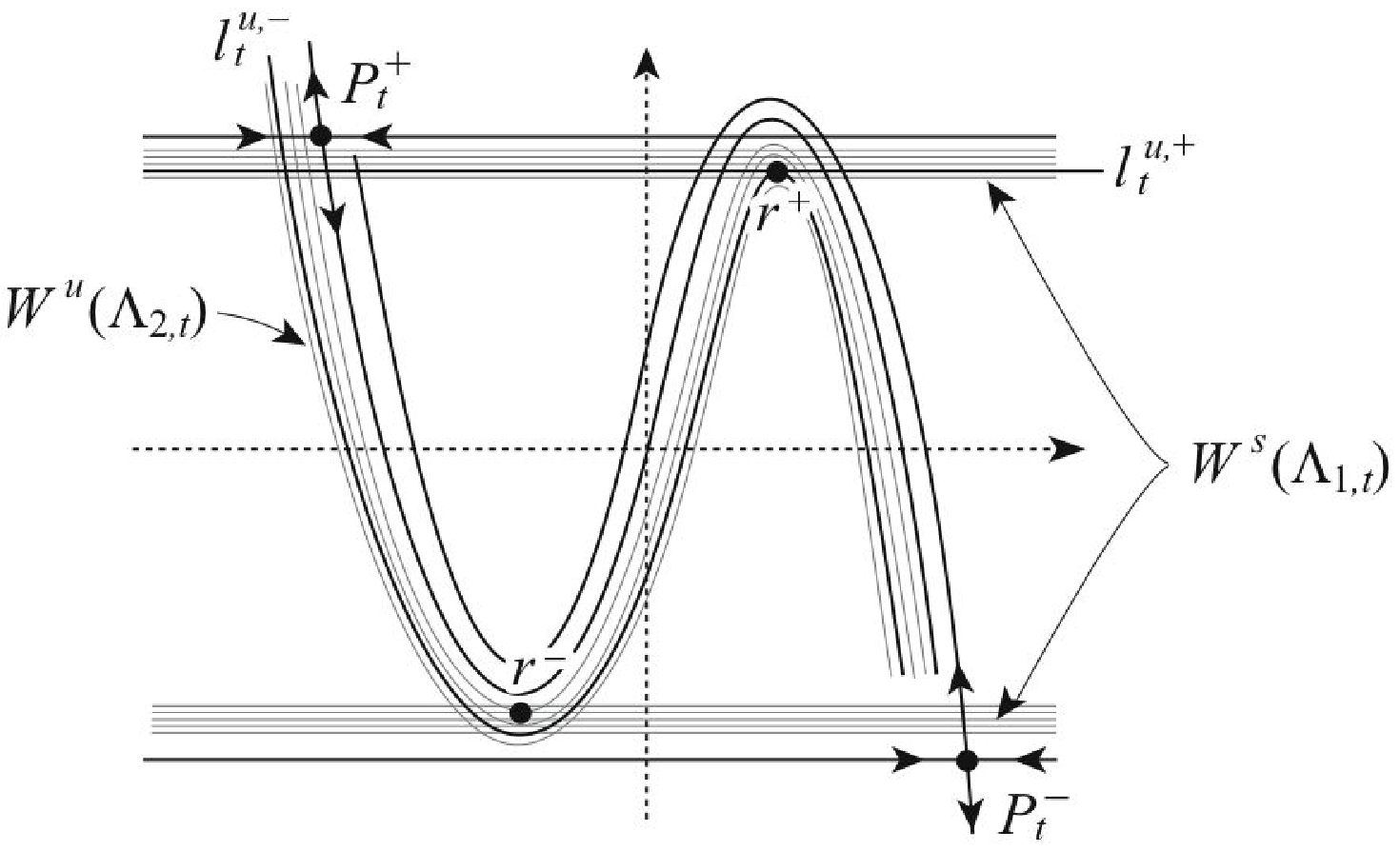}}
\caption{}
\label{fg_4_3}
\end{figure}

We consider the following situation for the above family $\{\varphi_t\}$.
Suppose that there exists a coordinate $\{(\bar x,\bar y):-5/2<\bar x<5/2,-5/2<\bar y<5/2\}$ on an open 
subset $U$ of $\mathbb{R}^2$ containing continuations $P^+_t,P^-_t$ of saddle periodic points of $\varphi_t$ 
with $P_t^\pm\approx (\mp 2,\pm 2)$.
We also assume that local stable manifolds $W^s_\mathrm{loc}(P_t^\pm)\subset U$ are well $C^1$-approximated by 
the horizontal lines $\bar y=\pm 2$, 
and local unstable manifolds $W^u_\mathrm{loc}(P_t^\pm)\subset U$ are contained in a sufficiently small 
neighborhood of the cubic curve $\{(\bar x,-\bar x^3+3\bar x):-5/2<\bar x<5/2\}$.
Then, we say that a subfamily $\{\varphi_t\}_{t_0<t<t_1}$ of $\{\varphi_t\}$ exhibits \emph{cubically 
related persistent contact-making and contact-breaking tangencies} if, for any $t\in (t_0,t_1)$, there 
exists a contact-making tangency 
$r_t^+\in l_t^{u,+}\cap l_t^{s,+}$ and a contact-breaking tangency $r_t^-\in l_t^{u,-}\cap l_t^{s,-}$ in $U$ 
with $r_t^\pm\approx (\pm 1,\pm 2)$, where $l_t^{s,\pm}$ (resp.\ $l_t^{u,\pm}$) are curves in 
$W^s(\Lambda_{1,t})\cap U$ (resp.\ $W^u(\Lambda_{2,t})\cap U$) arbitrarily $C^1$-close to 
$W^s_\mathrm{loc}(P_t^\pm)$ (resp.\ $W^u_\mathrm{loc}(P_t^\pm)$).
Figure \ref{fg_4_3} illustrates the cubically related situation, see also Fig.\ \ref{fg_5_1} in Section 
\ref{proofmainA}.

Here, we recall the standard method given by \cite{N1,PT} which guarantees the existence of persistent tangencies.
Suppose that $\{\varphi_t\}_{-\varepsilon\leq t\leq \varepsilon}$ is a one-parameter family in $\mathrm{Diff}^3(\mathbb{R}^2)$ with continuations of basic sets $\Lambda_t$ and $\Gamma_t$ admitting an unstable foliation $\mathcal{F}_t^u$ compatible with $\Lambda_t$ and a stable foliation $\mathcal{G}_t^s$ compatible with $\Gamma_t$.
We consider the situation where $\varphi_t$ has continuations of saddle fixed points $p_t\in \Lambda_t$ and $q_t\in \Gamma_t$ such that $W_{\mathrm{loc}}^s(p_t)$ (resp.\ $W_{\mathrm{loc}}^u(q_t)$) contains an arc $I_t$ (resp.\ $J_t$) crossing exactly $\mathcal{F}_t^u$ (resp.\ $\mathcal{G}_t^s$) and including $p_t$ (resp.\ $q_t$) as an end point.
Set $\mathcal{F}_t^{u(k)}=\varphi_t^k(\mathcal{F}_t^u)$ and $\mathcal{G}_t^{s(l)}=\varphi_t^{-l}(\mathcal{G}_t^s)$ for any integers $k,l\geq 0$.

\begin{figure}[hbtp]
\centering
\scalebox{0.6}{\includegraphics[clip]{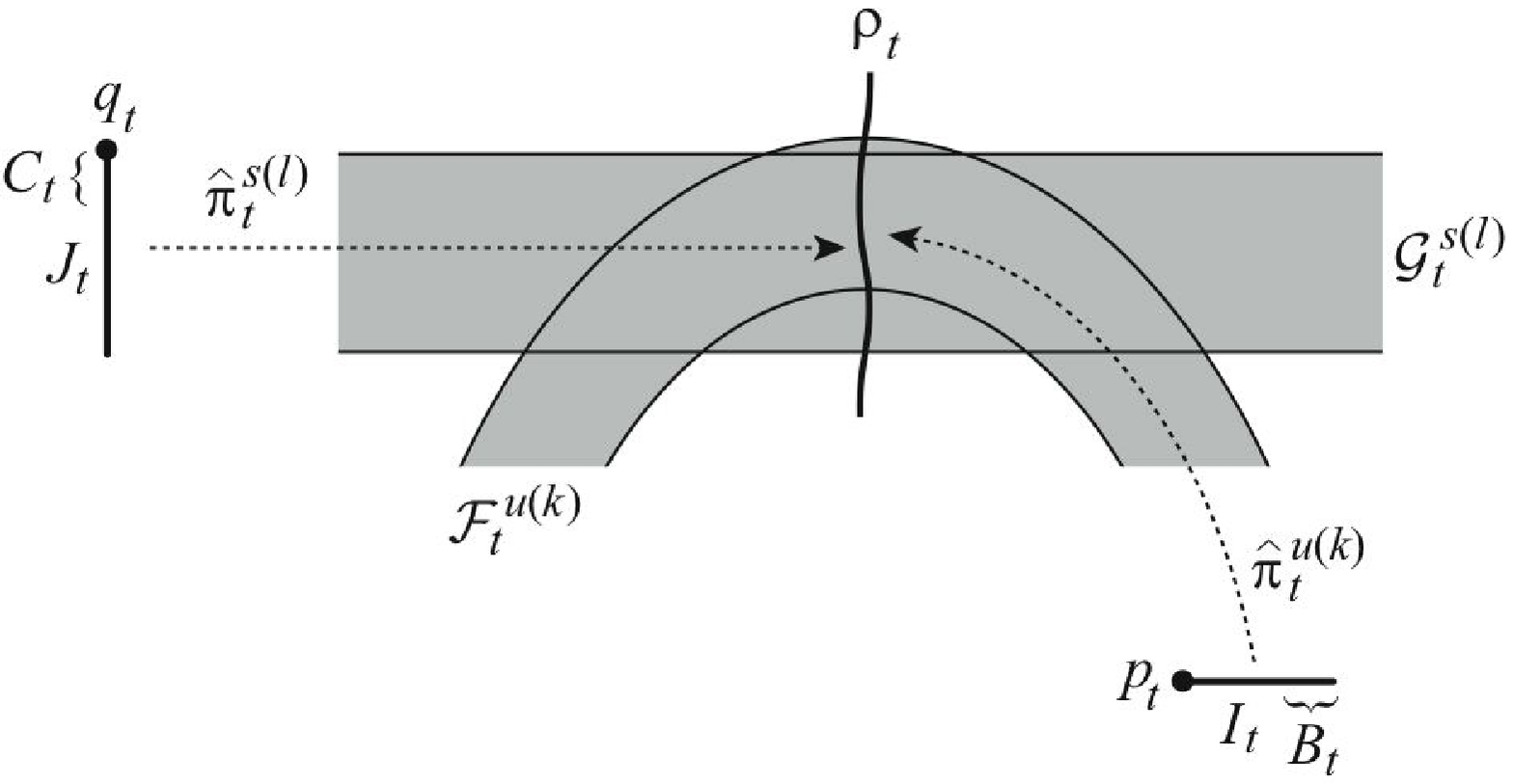}}
\caption{}
\label{fg_4_4}
\end{figure}

We suppose moreover that (i) the {\it highest} leaves of $\mathcal{F}_0^{u(k)}$ and $\mathcal{G}_0^{s(l)}$ 
have a quadratic tangency $r_0$ for some $k,l$, and (ii) for any $t$ sufficiently close to $0$, 
there exists a continuation of  $C^1$-curves $\rho_t$ in $\mathbb{R}^2$ 
containing two subarcs crossing $\mathcal{F}_t^{u(k)}$ and $\mathcal{G}_t^{s(l)}$ exactly and 
such that $\mathcal{F}_t^{u(k)}\cap \mathcal{G}_t^{s(l)}\cap \rho_t$ 
consists of tangencies of leaves in $\mathcal{F}_t^{u(k)}$ 
and $\mathcal{G}_t^{s(l)}$, see Fig.\ \ref{fg_4_4}.
Let $\hat\pi_t^{u(k)}:I_t\to \rho_t$ and $\hat\pi_t^{s(l)}:J_t\to \rho_t$ be the $C^1$-maps 
defined as (\ref{pi_t}) and (\ref{pi_t2}).

In our argument, the following assumption is crucial.
\begin{equation}\label{tau_tau}
\tau(\Lambda_t\cap I_t)\cdot \tau(\Gamma_t\cap J_t)>1.
\end{equation}
In Section \ref{proofmainA}, we will see that the assumption is satisfied in certain cases.

Let $B_0$ be a bridge of $\Lambda_0\cap I_0$ such that an end point of $B$ is mapped into the highest leave of $\mathcal{F}^{u(k)}_0$ by $\hat \pi^{u(k)}_0$.
A bridge $C_0$ of $\Gamma_0\cap J_0$ is given similarly.
By (\ref{thick_B}), if we choose $B_0,C_0$ sufficiently small, then $\tau(\hat\pi_0^{u(k)}(\Lambda_0\cap B_0))\cdot \tau(\hat\pi_0^{s(l)}(\Gamma_0\cap C_0))>1$ holds.
From the continuity of dynamically defined Cantor sets, we have
$$\tau(\hat\pi_t^{u(k)}(\Lambda_t\cap B_t))\cdot \tau(\hat\pi_t^{s(l)}(\Gamma_t\cap C_t))>1$$
for any $t$ sufficiently close to $0$, where $B_t$ and $C_t$ are the continuations of the bridges based at $B_0,C_0$ respectively.
Thus, Gap Lemma shows that $\hat\pi_t^{u(k)}(\Lambda_t\cap B_t)\cap \hat\pi_t^{s(l)}(\Gamma_t\cap C_t)\neq \emptyset$.
Then, any point $r_t\in \hat\pi_t^{u(k)}(\Lambda_t\cap B_t)\cap \hat\pi_t^{s(l)}(\Gamma_t\cap C_t)$ is a tangency of a $\Lambda_t$-leaf $l_t^u$ of $\mathcal{F}^{u(k)}_t$ and a $\Gamma_t$-leaf $l_t^s$ of $\mathcal{G}^{s(l)}_t$.
Note that $\{r_t\}$ is not a continuation on $t$, but one can suppose that each $r_t$ is contained in an arbitrarily small neighborhood of $r_0$ by taking the bridges $B_0,C_0$ sufficiently small.
Since $r_0$ is assumed to be a quadratic tangency, the leaves of $\mathcal{F}_0^{u(k)}$ and $\mathcal{G}_0^{s(l)}$ containing $r_0$ have the distinct curvatures at $r_0$.
From the continuity of the curvatures of leaves of compatible foliations and the closeness of $r_t$ to $r_0$, one can suppose that the curvatures of $l_t^u$ and $l_t^s$ at $r_t$ are mutually distinct.
This means that $r_t$ is a quadratic tangency.

Thus, we have shown that $W^u(\Lambda_t)$ and $W^s(\Gamma_t)$ have persistent quadratic tangencies $r_t$ under the assumption (\ref{tau_tau}).
But, we have \emph{not} shown yet that $r_t$ are contact-making (or breaking).
For the proof, we need Accompanying Lemma given in the next subsection.

\subsection{Accompanying Lemma}\label{ac_lemma}

\begin{figure}[hbtp]
\centering
\scalebox{0.65}{\includegraphics[clip]{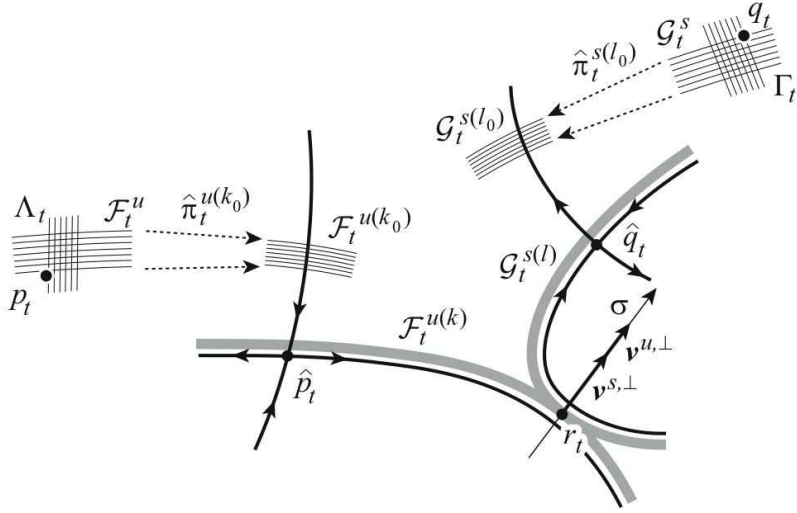}}
\caption{$\boldsymbol{v}^{u,\bot},\boldsymbol{v}^{s,\bot}$ represent the velocity vectors 
$d\eta_t^u/dt(t_0)$ and $d\eta_t^s/dt(t_0)$ respectively.}
\label{fg_4_5}
\end{figure}

We still work with the notation and situation of the previous subsections.
Suppose furthermore that there exist continuations of saddle fixed points $\hat p_t,\hat q_t$ of $\varphi_t$ other than $p_t,q_t$ such that $W^s(\hat p_t)\setminus \{\hat p_t\}$ has a subarc crossing $\mathcal{F}_t^{u(k_0)}=\varphi_t^{k_0}(\mathcal{F}^u_t)$ exactly and $W^u(\hat q_t)\setminus \{\hat q_t\}$ has a subarc crossing $\mathcal{G}_t^{s(l_0)}=\varphi_t^{-l_0}(\mathcal{G}^s_t)$ exactly for some integers $k_0,l_0\geq 0$.
If the integers $k,l$ given in Subsection \ref{p_tangency} are much greater than $k_0,l_0$, then the $\lambda$-lemma implies that $W^u_{\mathrm{loc}}(\hat p_t)$ (resp.\ $W^s_{\mathrm{loc}}(\hat q_t)$) is contained in a small neighborhood of $\mathcal{F}^{u(k)}_t$ (resp.\ $\mathcal{G}^{s(l)}_t$), see Fig.\ \ref{fg_4_5}.

\begin{figure}[hbtp]
\centering
\scalebox{0.6}{\includegraphics[clip]{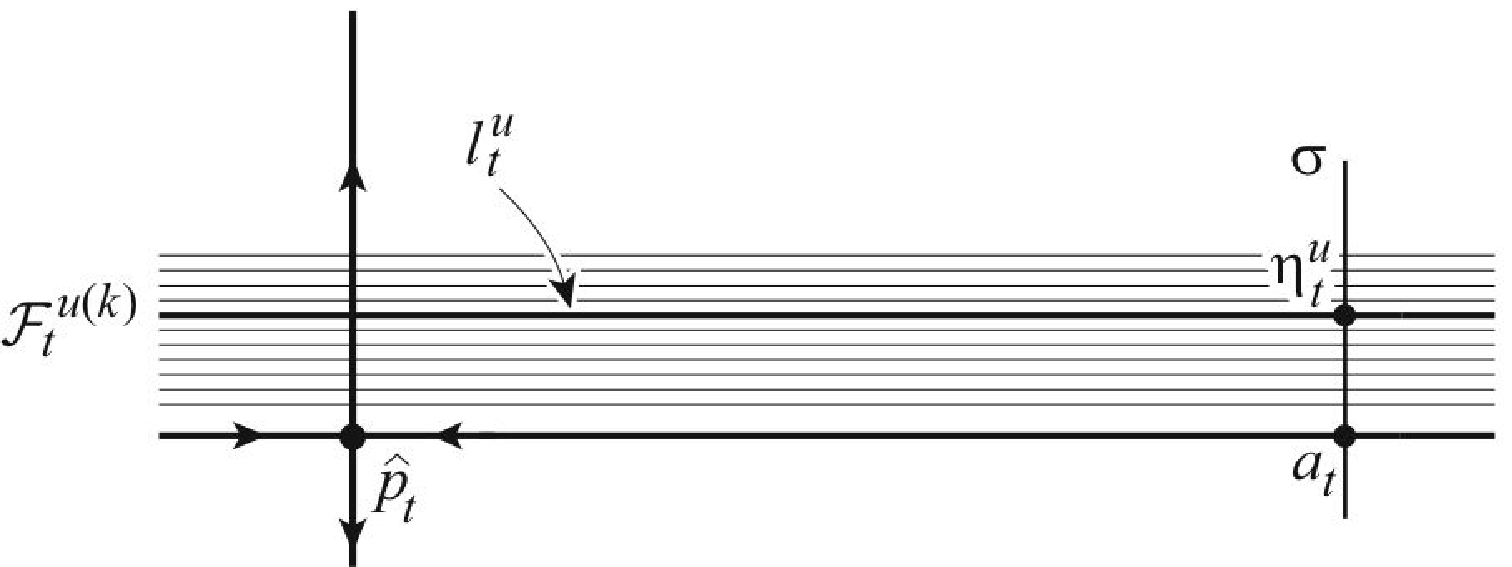}}
\caption{}
\label{fg_4_6}
\end{figure}

Here, we consider the case that the tangency $r_{t_0}$ of a $\Lambda_{t_0}$-leaf $l_{t_0}^u$ of $\mathcal{F}^{u(k)}_{t_0}$ and a $\Gamma_{t_0}$-leaf $l_{t_0}^s$ of $\mathcal{G}^{s(l)}_{t_0}$ is sufficiently close to both $W^u_{\mathrm{loc}}(\hat p_{t_0})$ and $W^s_{\mathrm{loc}}(\hat q_{t_0})$ for a $t_0\in (-\varepsilon,\varepsilon)$ .
Then, there exists a segment $\sigma$ meeting almost orthogonally all leaves of $\mathcal{F}^{u(k)}_{t_0}$, those of $\mathcal{G}^{s(l)}_{t_0}$, $W^u_{\mathrm{loc}}(\hat p_{t_0})$ and $W^s_{\mathrm{loc}}(\hat q_{t_0})$.
For any $t$ near $t_0$, let $\eta_{t}^u$ (resp.\ $\eta_{t}^s$) be the intersections of $\sigma$ and arcs $l_{t}^u$ in $\Lambda_{t}$-leaves of $\mathcal{F}^{u(k)}_{t}$ (resp.\ $l_{t}^s$ in $\Gamma_{t}$-leaves of $\mathcal{G}^{s(l)}_{t}$) which define a continuation starting at 
$l_{t_0}^u$ (resp.\ $l_{t_0}^s$), see Fig.\ \ref{fg_4_6}.
Note that $\eta_{t_0}^u=\eta_{t_0}^s=r_{t_0}$.
Propositions 1 and 2 in Pollicott \cite{Pol} based on arguments in \cite{LMM,Man,Ma} and others imply that 
$\{\eta_t^u\}$ and $\{\eta_t^s\}$ $C^1$-vary on $t$.
Fix a direction of $\sigma$ as shown in Fig.\ \ref{fg_4_5}, and suppose that $d\eta_t^u/dt(t_0)-d\eta_t^s/dt(t_0)$ is a non-zero vector.
If $d\eta_t^u/dt(t_0)-d\eta_t^s/dt(t_0)$ has the direction same as that of $\sigma$, then $r_{t_0}$ must be a contact-making tangency, and otherwise a contact-breaking tangency, for example see Lemma \ref{lem_5_2}.
The continuations $\{\eta_{t}^u\}$ and $\{\eta_{t}^s\}$ are away from their epicenters $\Lambda_{t}$ and $\Gamma_{t}$.
So, it may be difficult to estimate $d\eta_t^u/dt(t)$ and $d\eta_t^s/dt(t)$ by invoking these basic sets.
But, Accompanying Lemma as below suggests that $d\eta_t^u/dt(t)$ is well approximated by the velocity vector 
$da_t/dt(t_0)$ of the continuation $\{a_{t}\}$, where $a_{t}$ is the intersection point of $W^u_{\mathrm{loc}}(\hat p_{t})$ with $\sigma$.
A similar fact holds for $d\eta_t^s/dt(t)$.

\begin{lemma}[Accompanying Lemma]\label{lem_1_1}
For any $\delta_1>0$, there exists an integer $j_1>0$ and $\varepsilon_0>0$ such that  
any continuation $\{\eta_t^u\}$ of crossing points of $\sigma$ and $\Lambda_t$-leaves of $\mathcal{F}^{u(k_0+j)}_t$ satisfies
\begin{equation}\label{accomp}
\Bigl\Vert \frac{d  \eta_t^u}{d t}(t)-\frac{da_t}{dt}(t_0)\Bigr\Vert <\delta_1\end{equation}
if $j\geq j_1$ and $|t-t_0|<\varepsilon_0$.
\end{lemma}

The proof of this lemma is rather technical, so we defer it to Appendix.
The point $a_t$ is connected directly to $\hat p_t$ along $W^u_{\mathrm{loc}}(\hat p_t)$.
This implicitly suggests that $da_t/dt$ and hence $d  \eta_t^u/dt$ may be estimated from the behavior of $\hat p_t$.

\section{Proof of Theorem \ref{main_a}}\label{proofmainA}

\subsection{Outline of the proof}\label{outline_one}
We outline here the proof of Theorem \ref{main_a}.
Suppose that $\varphi\in \mathrm{Diff}^{\infty}(\mathbb{R}^2)$ and $\{\varphi_{\mu,\nu}\}\subset \mathrm{Diff}^{\infty}(\mathbb{R}^2)$ satisfy the conditions of Theorem \ref{main_a}.
That is, $\varphi=\varphi_{0,0}$ is a diffeomorphism with a dissipative and 
linearizable saddle fixed point $p$  and such that $W^u(p)$ and $W^s(p)$ has a cubic homoclinic tangency $q$ which unfolds generically with respect to $\{\varphi_{\mu,\nu}\}$.
Our proof consists of the three steps: 
\begin{enumerate}[(i)]
\item there is a renormalization at cubic homoclinic tangencies whose limit dynamics is the cubic map,
\item this map is conjugate to an N-map, and
\item the N-map has Cantor invariant sets with arbitrarily large thickness, 
\end{enumerate}
which is the cubic tangency version of Paris-Takens Theory \cite{PT} 
for one-parameter family of diffeomorphisms with a quadratic tangency unfolding generically.

Let $\{p_{\mu,\nu}\}$ be a continuation of saddle fixed points of $\varphi_{\mu,\nu}$ with $p_{0,0}=p$.
For a sufficiently large integer $k>0$, $\varphi_{\mu,\nu}^k$ has an invariant set (see Fig.\ \ref{fg_1_1}\,(b)) containing $p_{\mu,\nu}$ which defines a continuation of basic sets $\Lambda_{\mu,\nu}^{\mathrm{out}}$ of $\varphi_{\mu,\nu}$. 
From the continuity of dynamically defined Cantor set, there exists $\alpha>0$ with $\tau(\Lambda_{\mu,\nu}^{\mathrm{out}}\cap W^s(p_{\mu,\nu}))>\alpha$ for any $(\mu,\nu)$ near $(0,0)$, see Lemma \ref{lem_2_3}.

In Subsection \ref{renorm}, we renormalize $\varphi_{\mu,\nu}$ in a small neighborhood of $q$ by using a sequence of $C^4$-diffeomorphisms $\Phi_n$ and that of reparametrizations $\Theta_n$ so that the diffeomorphisms $\psi_{\bar\mu,\bar\nu,n}$ on $\mathbb{R}^2$ with $\psi_{\bar\mu,\bar\mu,n}(\bar{x},\bar{y})=\Phi_{n}^{-1}\circ\varphi_{\mu,\nu}^{N+n}\circ \Phi_{n}(\bar{x},\bar{y})$ $C^3$-converges to the endomorphism $\psi_{\bar\mu,\bar\nu}$ with $\psi_{\bar{\mu},\bar{\nu}}(\bar{x}, \bar{y})=(\bar{y},F_{\bar\mu,\bar\nu}(\bar y) )$ as $n\rightarrow \infty$, where $N$ is some fixed positive integer, $F_{\bar\mu,\bar\nu}$ is the cubic map with $F_{\bar\mu,\bar\nu}(\bar y)= - \bar{y}^{3}+\bar{\mu} \bar{y}+\bar{ \nu}$ and $(\bar\mu,\bar\nu)=\Theta_n^{-1}(\mu,\nu)$, see Lemma \ref{lem_3_1}.
We note that $\psi_{\bar\mu,\bar\nu}$ has two fixed points $P^\pm_{\bar\mu,\bar\nu}\approx (\mp 2,\pm 2)$ for any $(\bar\mu,\bar\nu)$ near $(3,0)$.
In particular, when $(\bar\mu,\bar\nu)=(3,0)$, $W^u(P^+_{3,0})$ is equal to $W^u(P^-_{3,0})$ and it has quadratic tangencies with $W^s(P^+_{3,0})$ and $W^s(P^-_{3,0})$, see Fig.\ \ref{fg_5_2}.
Moreover, for any even number $m\geq 6$, each $\psi_{\bar\mu,\bar\nu,n}$ has a basic set $\Lambda_{\bar\mu,\bar\nu,n}^m$ and a saddle periodic point $Q^m_{\bar\mu,\bar\nu,n}$ converging uniformly to a certain basic set $\Lambda_{\bar\mu,\bar\nu}^m$ and a saddle periodic point $Q_{\bar\mu,\bar\nu}^m$ of $\psi_{\bar\mu,\bar\nu}$ respectively as $n\rightarrow \infty$, see Fig.\ \ref{fg_5_3}.

To see that the inequality condition corresponding to (\ref{tau_tau}) holds, 
we need to show that $\tau(\Lambda_{\bar\mu,\bar\nu,n}^m\cap W^u(Q^m_{\bar\mu,\bar\nu,n}))>1/\alpha$.
In fact, this inequality is derived from the limit inequality $\tau(\Lambda_{3,0}^m\cap W^u(Q^m_{3,0}))>1/\alpha$ 
when $n$ is sufficiently large and $(\bar\mu,\bar\nu)$ is sufficiently close to $(3,0)$.
This is accomplished by the two facts that (i) the function $F_{3,0}$ on $[-2,2]$ 
is conjugate to the N-map $S$ on $[-3/2,3/2]$ defined in \ref{affine} via 
the homeomorphism $h:[-3/2,3/2]\longrightarrow [-2,2]$ with $h(x)=2\sin(\pi x/3)$ (Proposition \ref{prop_4_2}) and (ii) the dynamically defined Cantor set $K_m$ in $[-3/2,3/2]$ for $S$ corresponding to $\Lambda^m_{3,0}$ satisfies $\lim_{m(\mathrm{even})\rightarrow \infty}\tau(K_m)=\infty$ (Lemma \ref{lem_4_1}\,(iv)).
Thus, for any sufficiently large even number $m$, we have the inequality
$$\tau(\Lambda_{\mu,\nu}^{\mathrm{out}}\cap W^s(p_{\mu,\nu}))\cdot \tau(\Lambda_{\bar\mu,\bar\nu,n}^m\cap W^u(Q^m_{\bar\mu,\bar\nu,n}))>1$$
corresponding to (\ref{tau_tau}).

\begin{figure}[htb]
\centering
\scalebox{0.61}{\includegraphics[clip]{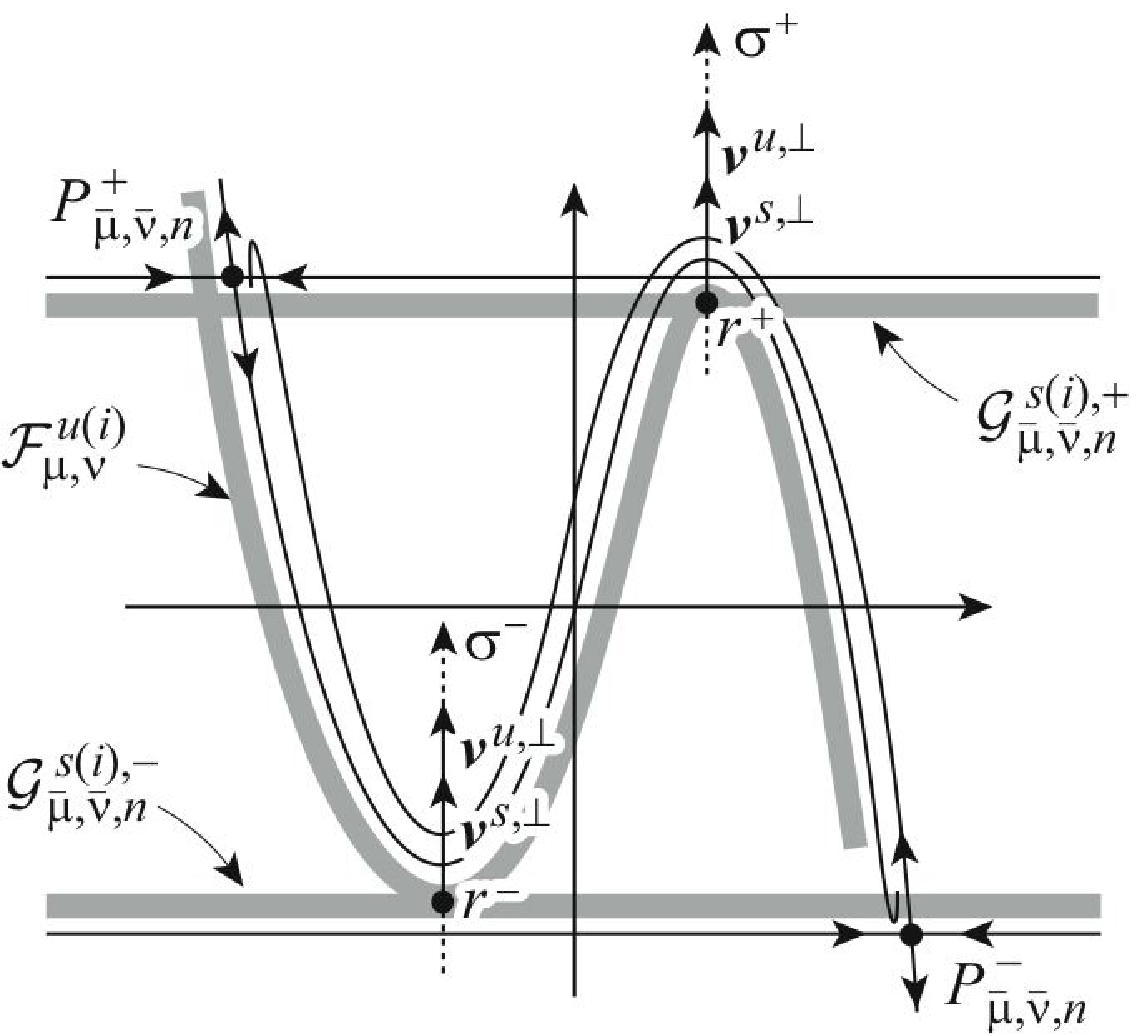}}
\caption{ The upper $\boldsymbol{v}^{u,\bot}$, $\boldsymbol{v}^{s,\bot}$ represent $d  \eta^{u,+}_t/dt(t)$ and $d  \eta^{s,+}_t/dt(t)$ respectively. The lower $\boldsymbol{v}^{u,\bot}$, $\boldsymbol{v}^{s,\bot}$ represent $d  \eta^{u,-}_t/dt(t)$ and $d  \eta^{s,-}_t/dt(t)$.
The tangencies $r^\pm_{\bar c(t),n}$ are indicated shortly by $r^\pm$.}
\label{fg_5_1}
\end{figure}

In Subsection \ref{Hyperbolic}, we will show that, for an unstable foliation $\mathcal{F}^u_{\mu,\nu}$ compatible with $\Lambda_{\mu,\nu}^{\mathrm{out}}$, $W^s(P^+_{\bar\mu,\bar\nu,n})\setminus \{P^+_{\bar\mu,\bar\nu,n}\}$ has a subarc crossing $\mathcal{F}^{u(k_0)}_{\mu,\nu}=\psi_{\bar\mu,\bar\nu,n}^{k_0}(\mathcal{F}^u_{\mu,\nu})$ exactly, where the $\Phi_n^{-1}$-image of $\mathcal{F}^u_{\mu,n}$ is also denoted by the same notation.
Thus, if $i$ is sufficiently greater than $k_0$, a small neighborhood of $\mathcal{F}^{u(i)}_{\mu,\nu}$ in the $\bar x\bar y$-plane contains $W^u_{\mathrm{loc}}(P^+_{\bar\mu,\bar\nu,n})$, see Fig.\ \ref{fg_5_1}.
On the other hand, one can choose a stable foliation $\mathcal{G}^s_{\bar\mu,\bar\nu,n}$ compatible with $\Lambda^m_{\bar\mu,\bar\nu,n}$ so that both $W^u(P^+_{\bar\mu,\bar\nu,n})\setminus \{P^+_{\bar\mu,\bar\nu,n}\}$ and $W^u(P^-_{\bar\mu,\bar\nu,n})\setminus \{P^-_{\bar\mu,\bar\nu,n}\}$ have subarcs crossing $\mathcal{G}^s_{\bar\mu,\bar\nu,n}$ exactly.
Thus, for any sufficiently large integer, we have two foliations $\mathcal{G}^{s(i),\pm}_{\bar\mu,\bar\nu,n}$ obtained
by shortening the leaves of $\psi^{-2i}(\mathcal{G}^s_{\bar\mu,\bar\nu,n})$ 
as illustrated in  Fig.\ \ref{fg_5_1}.

Let us here give the Table \ref{Table1} which may be helpful to read the remaining parts of this section.
The above basic sets/foliations and related saddle points in this section
are listed  as pairs in the first column, 
and 
corresponding pairs of notations given in Sections \ref{C_foliations_Acc_lemma} 
are listed in the second column 
of Table \ref{Table1}.
\begin{table}[hbt]
{\renewcommand\arraystretch{1.2}
\begin{tabular}{c|c}
Subsections \ref{outline_one}, \ref{Existence}& Subsection \ref{ac_lemma} \\ \hline
$(\Lambda_{\mu,\nu}^{\mathrm{out}},p_{\mu,\nu})$&$(\Lambda_t,p_t)$ \\ 
$(\Lambda_{\bar\mu,\bar\nu,n}^m,Q^m_{\bar\mu,\bar\nu,n})$& $(\Gamma_t,q_t)$\\ 
$(\mathcal{F}^{u(i)}_{\mu,\nu},P^+_{\bar\mu,\bar\nu,n})$& $(\mathcal{F}^{u(k)}_t,\hat p_t)$\\ 
$(\mathcal{G}^{s(i),\pm}_{\bar\mu,\bar\nu,n},P^\pm_{\bar\mu,\bar\nu,n})$& $(\mathcal{G}^{s(l)}_t,\hat q_t)$\\
\hline
\end{tabular}
}
\bigskip
\caption{}\label{Table1}
\end{table}

For any $n$ not less than some integer $n_0>0$, there exists $(\bar \mu_1,\bar \nu_1)$ near $(3,0)$ such that the highest (resp.\ lowest) leaves of $\mathcal{F}^{u(i)}_{\mu_1,\nu_1}$ and $\mathcal{G}^{s(i),+}_{\bar\mu_1,\bar\nu_1,n}$ (resp.\ $\mathcal{G}^{s(i),-}_{\bar\mu_1,\bar\nu_1,n})$ have a quadratic tangency.
Here, we use the fact that our family $\{\varphi_{\mu,\nu}\}$ has \emph{two} parameters essentially to obtain these \emph{two} quadratic tangencies simultaneously.
As was seen in Subsection \ref{p_tangency}, for any $(\bar\mu,\bar\nu)$ in a sufficiently small open neighborhood $\mathcal{N}_n$ of $(\bar\mu_1,\bar\nu_1)$, $\mathcal{F}^{u(i)}_{\mu,\nu}$ and $\mathcal{G}^{s(i),+}_{\bar\mu,\bar\nu,n}$ (resp.\ $\mathcal{G}^{s(i),-}_{\bar\mu,\bar\nu,n}$) have a quadratic tangency $r^+_{\bar\mu,\bar\nu,n}$ (resp.\ $r^-_{\bar\mu,\bar\nu,n}$).
From the definition of $\Theta_n$, we know that $\Theta_n(\mathcal{N}_n)$ converges to the origin $\mathbf{0}$ of the $\mu\nu$-plane as $n\rightarrow \infty$.
Thus, $\mathcal{O}=\bigcup_{n\geq n_0}\Theta_n(\mathcal{N}_n)$
is an open set with $\mathrm{Cl}(\mathcal{O})\ni \mathbf{0}$.

For any $(\bar\mu,\bar\nu)\in \mathcal{N}_n$, consider a regular curve $c:(-\delta,\delta)\longrightarrow \Theta_n(\mathcal{N}_n)$ such that $\bar c(t)=\Theta_n^{-1}\circ c(t)=(\bar\mu,\bar \nu+t)$ for a small $\delta>0$.
For \emph{any} $t_0\in (-\delta,\delta)$ and $\bullet:=\pm$, let $\{\eta_{t}^{u,\bullet}\}$, $\{\eta_{t}^{s,\bullet}\}$ be continuations  based at the quadratic tangency $r^\bullet_{\bar c(t_0),n}$ defined as in Subsection \ref{ac_lemma} with respect to the vertical segment $\sigma^\bullet$ passing though $r^\bullet_{\bar c(t_0),n}$.
By using Accompanying Lemma (Lemma \ref{lem_1_1}), we show in Subsection \ref{Existence} that the velocity vectors of $\eta_t^{u,\bullet}$ and $\eta_t^{s,\bullet}$ satisfy $d\eta_t^{u,\bullet}/dt(t_0)-d
\eta_t^{s,\bullet}/dt(t_0)\approx (0,9/10)$.
This implies that $r^+_{\bar c(t_0),n}$ is a contact-making tangency and $r^-_{\bar c(t_0),n}$ is a contact-breaking tangency with respect to the one-parameter family $\{\psi_{\bar c(t),n}\}$ and hence to $\{\varphi_{c(t)}\}$.
This ends the outline of proof
which will be given in
the next subsections.

\subsection{Basic sets with large stable thickness}\label{Hyperbolic} 
Now, we begin full-fledged discussions to show Theorem \ref{main_a} from this subsection.

Let $\{\varphi_{\mu,\nu}\}$ be any two-parameter family of diffeomorphisms on the plane satisfying the conditions in Setting \ref{set}.
For any $n\geq 0$ and $(\bar\mu,\bar\nu)\in \bar\Sigma=[0,4]\times [-1,1]$, let $\psi_{\bar\mu,\bar\nu,n}$ be the diffeomorphism on the plane defined by
\begin{equation}\label{eqn_3}
\psi_{\bar\mu,\bar\mu,n}(\bar{x},\bar{y})=\Phi_{n}^{-1}\circ\varphi_{\Theta_{n}(\bar\mu,\bar\nu)}^{N+n}\circ \Phi_{n}(\bar{x},\bar{y}).
\end{equation}
By Lemma \ref{lem_3_1}, the sequence $\{\psi_{\bar{\mu}, \bar{\nu},n}\}_n$ $C^3$-converges as $(\bar x,\bar y,\bar \mu,\bar \nu)$-functions to the endomorphism $\psi_{\bar\mu,\bar \nu}$ on $\mathbb{R}^2$ defined by $\psi_{\bar\mu,\bar\nu}(\bar x,\bar y)=(\bar y,F_{\bar \mu,\bar \nu}(\bar y))$, where
$$F_{\bar{\mu},\bar{\nu}}(\bar{y})= -\bar{y}^{3}+\bar{\mu}  \bar{y}+ \bar{\nu}.$$

\begin{proposition} \label{prop_5_1}
For any even number $m\geq 6$, there exists a neighborhood $\mathcal{N}$ of $(3,0)$ in $\mathrm{Int}(\bar \Sigma)$ and an integer $n_0>0$ satisfying the following conditions.
\begin{enumerate}[\rm (i)]
\item 
For any $(\bar{\mu},\bar{\nu})\in \mathcal{N}$ and any integer $n\geq n_0$, $\psi_{\bar{\mu},\bar{\nu},n}$ has two $2$-periodic points $P^{\pm}_{\bar{\mu},\bar{\nu}, n}$, an $m$-periodic point $Q^{m}_{\bar{\mu},\bar{\nu}, n}$ and 
a Cantor invariant set $\Lambda^{m}_{\bar{\mu},\bar{\nu}, n}$ with $\{P^\pm_{\bar \mu,\bar \nu,n}\}\rightarrow (\pm 2,\mp 2)$, $\{Q^m_{\bar \mu,\bar \nu,n}\}\rightarrow (\bar{q}(m),  F(\bar{q}(m)))$ and $\{\Lambda_{\bar\mu,\bar\nu,n}^m\}\rightarrow \Lambda^{m}_{3,0}=\{ (\bar{y}, F(\bar{y})) : \bar{y}\in   \bar{K}_{m} \}$ as $(\bar \mu,\bar \nu)\rightarrow (3,0)$ and $n\rightarrow \infty$, where $F=F_{3,0}$ and $\bar{q}=\bar{q}(m)$ is the $m$-periodic point of $F$ given in Proposition  \ref{prop_4_3}.
Moreover, 
$\Lambda^{m}_{\bar{\mu},\bar{\nu}, n}$ is a basic set containing $Q^{m}_{\bar{\mu},\bar{\nu}, n}$.
\item
For any $n\geq n_0$, the stable thickness $\tau_{\bar\mu,\bar\nu,n}^m$ of the Cantor set $\Lambda^{m}_{\bar{\mu},\bar \nu, n}\cap W^u(Q^m_{\bar{\mu},\bar \nu,n})$ is greater than an arbitrarily given constant if we take $m$ sufficiently large and $\mathcal{N}$ sufficiently small.
\end{enumerate}
\end{proposition}

\begin{proof} 
The points $P^{\pm}_{3,0}=(\mp 2, \pm 2)$ are 2-periodic saddle points for $\psi_{3,0}$.
Consider the arcs $\tilde{I}_{i}=\{(\bar{y}, F(\bar{y})): \bar{y}\in h(I_{i})\}$ $(i=0,1,\cdots,m-1)$ contained in the cubic curve $l=\mathrm{Im}(\psi_{3,0})$ as in Fig.\ \ref{fg_5_2}, where $h$ is the homeomorphism given in the proof of Proposition \ref{prop_4_2} and $I_{i}$ are the intervals given in the proof of Lemma \ref{lem_4_1}.
Note that the point $Q^{m}_{3,0}=( \bar{q}(m), F(\bar{q}(m)))$ is an an $m$-periodic of $\psi_{3,0}$ contained in the boundary of $\tilde{I}_{m-1}$.
\begin{figure}[hbt]
\begin{center}
\includegraphics{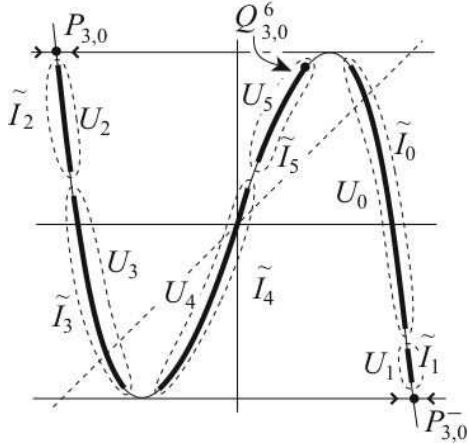}
\caption{The case of $m=6$.}
\label{fg_5_2}
\end{center}
\end{figure}

Now, we show the existence of a basic set $\Lambda^m_{\bar{\mu},\bar \nu,n}$ satisfying the conditions as above by using arguments similar to those in Proposition 1 in \cite[\S 6.3]{PT}.
Take mutually disjoint and sufficiently thin neighborhoods $U_{i}$ of $\tilde{I}_{i}$ $(i=0,1,\cdots, m-1)$ such that, for any horizontal line $\gamma$ with $\gamma \cap U_i\neq \emptyset$, the intersection $\gamma \cap U_i$ consists of a single short segment.
For each point $p$ of $\mathcal{U}_m=\bigcup_{i=0}^{m-1}U_i$, we take the unstable cone so that it does not contains the horizontal direction.
Moreover, when the point $p$ is in the cubic curve $l$, we take the unstable cone so that it contains the tangent direction of $l$.
By Proposition \ref{prop_4_3}\,(ii), the set $\Lambda^m_{3,0}$ as above is an expanding hyperbolic set for $\psi_{3,0}|l$.
It follows that there exists an integer $n(m)>0$ and a constant $C>1$ such that, for any $p\in \mathcal{U}_m$ with $\psi_{3,0}^i(p)\in \mathcal{U}_m$ for $i=1,\dots,n(m)$, $\Vert d\psi_{3,0}^{n(m)}(\boldsymbol{v})\Vert \geq C\Vert \boldsymbol{v}\Vert$ for any $\boldsymbol{v}$ in the unstable cone at $p$.
This implies that, for all sufficiently large integer $n>0$ and any $(\bar\mu,\bar\nu)$ sufficiently close to $(3,0)$, the vectors in the unstable cones are eventually growing with respect to $\psi_{\bar\mu,\bar\nu,n}$ if they do not leave $\mathcal{U}_m$.
The corresponding fact for the stable cones is evident.
From these facts, we have an hyperbolic invariant set $\Lambda^m_{\bar{\mu},\bar \nu,n}$ for $\psi_{\bar\mu,\bar\nu,n}$ in $\mathcal{U}_m$ with $\{\Lambda_{\bar\mu,\bar\nu,n}^m\}\rightarrow \Lambda^{m}_{3,0}$ and containing an $m$-periodic point $Q^m_{\bar\mu,\bar\nu,n}$ which converges to $Q^m_{3,0}$ as $(\bar\mu,\bar\nu)\rightarrow (3,0)$, $n\rightarrow \infty$.
One can prove that the $\Lambda_{\bar\mu,\bar\nu,n}^m$ is a basic set by using the Markov partition for $\psi_{\bar\mu,\bar\nu,n}$ consisting of $m$ boxes which are near $\mathcal{U}_m$, see \cite[Appendix 2]{PT} for details.

The assertion (ii) is verified by using Proposition \ref{prop_4_3}\,(iii) together with the continuity of thickness (Theorem 1 in \cite[\S 4.3]{PT}).
\end{proof}

By Proposition \ref{prop_5_1}, one can choose an even number $m\geq 6$, an open neighborhood $\mathcal{N}$ of $(3,0)$ and an integer $n_0>0$ so that, for any $(\bar \mu,\bar \nu)\in \mathcal{N}$, $n\geq n_0$ and $i\geq 0$,
\begin{equation}\label{eqn_2}
\tau^m_{\bar\mu,\bar\nu,n}\geq \frac{2}{\alpha},
\end{equation}
where $\alpha>0$ is the constant given in Lemma \ref{lem_2_3}.
The even number $m$ satisfying (\ref{eqn_2}) is fixed throughout the remainder of this section.
However, we may retake the neighborhood $\mathcal{N}$ of $(3,0)$ by a smaller one and the integer $n_0$ by a larger one as needed.

\begin{figure}[hbt]
\begin{center}
\scalebox{0.85}{\includegraphics[clip]{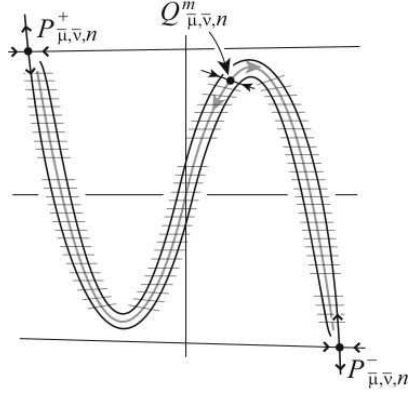}}
\caption{ The gray curve represents $W^u_{\mathrm{loc}}(Q_{\bar\mu,\bar\nu,n}^m)$.
The short segments are leaves of $W^s_{\mathrm{loc}}(\Lambda_{\bar\mu,\bar\nu,n}^m)$.}
\label{fg_5_3}
\end{center}
\end{figure}

\begin{figure}[hbt]
\begin{center}
\scalebox{0.85}{\includegraphics[clip]{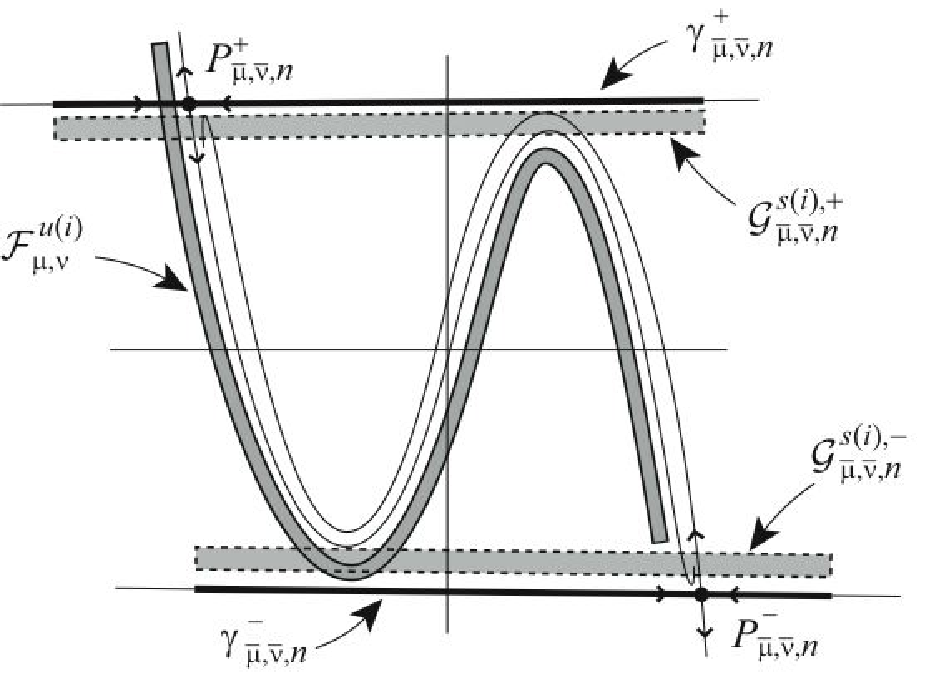}}
\caption{}
\label{fg_5_4}
\end{center}
\end{figure}

For any sufficiently large integer $n>0$, there exists a foliation $\mathcal{G}^s_{\bar\mu,\bar\nu,n}$ compatible with $W^s_{\mathrm{loc}}(\Lambda^m_{\bar\mu,\bar\nu,n})$ such that both $W^u(P^+_{\bar\mu,\bar\nu,n})$ and $W^u(P^-_{\bar\mu,\bar\nu,n})$ have subarcs crossing $\mathcal{G}^s_{\bar\mu,\bar\nu,n}$ exactly, see Fig.\ \ref{fg_5_3}.
For any sufficiently large integer $i>0$, one can have foliations 
$\mathcal{G}^{s(i),\pm}_{\bar\mu,\bar\nu,n}$ obtained by shortening the leaves of 
$\psi^{-2i}_{\bar\mu,\bar\nu,n}(\mathcal{G}^s_{\bar\mu,\bar\nu,n})$ 
so that each leaf of $\mathcal{G}^{s(i),\pm}_{\bar\mu,\bar\nu,n}$ is contained in a small neighborhood of 
$\gamma^\pm_{\bar\mu,\bar\nu,n}$ respectively, where $\gamma^\pm_{\bar\mu,\bar\nu,n}$ are arcs in 
$W^s(P^\pm_{\bar\mu,\bar\nu,n})$ 
such that the end points of $\gamma^+_{\bar\mu,\bar\nu,n}$ are close to $(-5/2,2),(2,2)$ and those of 
$\gamma^-_{\bar\mu,\bar\nu,n}$ are close to $(-2,-2), (5/2,-2)$, see Fig.\ \ref{fg_5_4}.

In any case that $\bar \mu,\bar \nu,n$ and $\mu,\nu$ appear simultaneously, for just simplicity,
$\Theta_n(\bar \mu,\bar \nu)$ will be denoted by $(\mu,\nu)$. 
For a subset $A$ of the $xy$-plane, the inverse image $\Phi_n^{-1}(A)$ in the $\bar x\bar y$-plane 
will be also denoted by $A$.

Let $\Lambda_{\mu,\nu}^{\mathrm{out}}$ be the basic invariant set of $\varphi_{\mu,\nu}$, and $\mathcal{O}_1$ the neighborhood of $(0,0)$ in the $\mu\nu$-plane given in Lemma \ref{lem_2_3}.
We may assume that $\Theta_n(\mathcal{N})\subset \mathcal{O}_1$ for any $n\geq n_0$.
Let $\mathcal{F}^u_{\mu,\nu}$ be a foliation compatible with $W^u_{\mathrm{loc}}(\Lambda_{\mu,\nu}^{\mathrm{out}})$.
As in the proof of Proposition 1 in \cite[\S 6.4]{PT}, there exists a compact arc $\sigma^s_{\bar\mu,\bar\nu,n}$ in $W^s(P^+_{\bar\mu,\bar\nu,n})$ containing $P^+_{\bar\mu,\bar\nu,n}$ and converging as $n\rightarrow \infty$ to an arc in $W^s(p_{\mu,\nu})$ which contains at least one fundamental domain of $W^s(p_{\mu,\nu})$.
Moreover, an argument as in \cite[p.\ 129, Remark 1]{PT} shows that, for all sufficiently large $n$ and some $j>0$, $\varphi^{-(n+j)}_{\mu,\nu}(\sigma^s_{\bar\mu,\bar\nu,n})$ meets $\mathcal{F}_{\mu,\nu}^u$ non-trivially and transversely.
From this fact together with the $\lambda$-lemma, one can have a foliation $\mathcal{F}^{u(i)}_{\mu,\nu}$ obtained by shortening the leaves of $\psi^{2i}_{\bar\mu,\bar\nu,n}(\mathcal{F}^u_{\mu,\nu})$ such that each leaf of $\mathcal{F}^{u(i)}_{\mu,\nu}$ is contained in a small neighborhood of a compact arc in $W^u(P^+_{\bar\mu,\bar\nu,n})$ passing through $P^+_{\bar\mu,\bar\nu,n}$ and ending at a point near $P^-_{\bar\mu,\bar\nu,n}$ if $i$ is an integer sufficiently greater than $j$.
Note that $W^s(P^+_{\bar\mu,\bar\nu,n})$ has a subarc crossing $\mathcal{F}^{u(i)}_{\mu,\nu}$ exactly, 
see Fig.\ \ref{fg_5_4} again.

\subsection{Existence of persistent contact-making and breaking tangencies}\label{Existence}
In this section,  we will benefit from Accompanying Lemma in Subsection \ref{ac_lemma}.
Table \ref{Table1} in Subsection \ref{outline_one} 
may be helpful for the reader 
to grasp the correspondence between our settings here and general settings in Subsection \ref{p_tangency} and \ref{ac_lemma}.

Recall that $P^\pm_{\bar\mu,\bar\nu}=(\bar x_1^\pm,\bar y_1^\pm)$ are fixed points of the endomorphism $\psi^2_{\bar\mu,\bar\nu}$ defined by
$$\psi^2_{\bar\mu,\bar\nu}(\bar x,\bar y)=(F_{\bar\mu,\bar\nu}(\bar y),F_{\bar\mu,\bar\nu}^2(\bar y)).$$
Thus,
$$\bar y_1^\pm=-(-(\bar y_1^\pm)^3+\bar \mu \bar y_1^\pm+\bar \nu)^3+\bar \mu(-(\bar y_1^\pm)^3+\bar \mu \bar y_1^\pm+\bar \nu)+\bar \nu.$$
By regarding $\bar y_1^\pm$ as functions of $(\bar\mu,\bar\nu)$ satisfying the above equation and applying the Implicit Function Theorem to $\bar y_1^\pm$, one can show at $\bar \mu=3,\bar \nu=0,\bar y_1^\pm=\pm 2$ that
\begin{equation}\label{y_1}
\frac{\partial \bar y_1^\pm}{\partial \bar \mu}=\pm\frac{1}{4}\quad\mbox{and}\quad
\frac{\partial \bar y_1^\pm}{\partial \bar \nu}=\frac{1}{10}.
\end{equation}

The critical points of $F_{\bar\mu,\bar\nu}(\bar y)$ are $c^\pm_{\bar\mu,\bar\nu}= \pm\sqrt{\bar \mu/3}$.
In the cases when $\bar \nu=0$ and $\bar \mu=3$, we have respectively
$$F_{\bar\mu,0}(c^\pm_{\bar\mu,0})=\pm \frac{2}{3\sqrt{3}}\bar\mu^{3/2}\quad\mbox{and}\quad F_{3,\bar\nu}(c^\pm_{3,\bar\nu})=\bar\nu\pm 2.$$
These imply that
\begin{equation}\label{y_2}
\left.\frac{d}{d \bar \mu} F_{\bar\mu,0}(c^\pm_{\bar\mu,0})\right|_{\bar \mu=3}=\pm 1\quad\mbox{and}\quad \left.\frac{d}{d \bar \nu} F_{3,\bar\nu}(c^\pm_{3,\bar\nu})\right|_{\bar \nu=0}=1.
\end{equation}

\begin{figure}[hbt]
\begin{center}
\scalebox{0.8}{\includegraphics[clip]{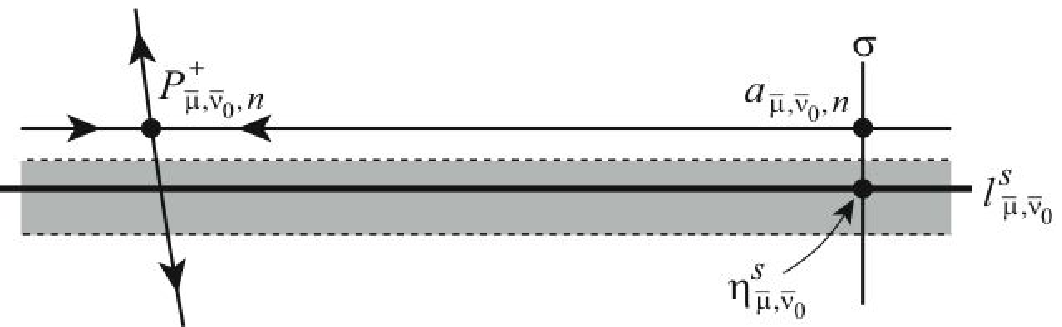}}
\caption{}
\label{fg_5_5}
\end{center}
\end{figure}

For a fixed $(\bar\mu_0,\bar\nu_0)$ near $(3,0)$, let $r_{\bar\mu_0,\bar\nu_0}$ be a point in a $\Lambda^m_{\bar\mu_0,\bar\nu_0,n}$-leaf $l_{\bar\mu_0,\bar\nu_0}^s$ of $\mathcal{G}^{s(i),+}_{\bar\mu_0,\bar\nu_0,n}$.
Recall that, as was defined in Subsection \ref{c_foliation}, $l_{\bar\mu_0,\bar\nu_0}^s$ being a $\Lambda^m_{\bar\mu_0,\bar\nu_0,n}$-leaf means that $l_{\bar\mu_0,\bar\nu_0}^s$ is contained in a leaf of $W^s(\Lambda^m_{\bar\mu_0,\bar\nu_0,n})$.
Suppose that $\sigma$ is a short segment passing through $l_{\bar\mu_0,\bar \nu_0}^s$ at $r_{\bar\mu_0,\bar\nu_0}$ and parallel to the $\bar y$-axis.
For any $\bar\mu$ near $\bar\mu_0$, let $\eta_{\bar\mu,\bar\nu_0}^s$ be the intersection point of $\sigma$ and a $\bar \mu$-continuation $l_{\bar\mu,\bar\nu_0}^s\subset W^s(\Lambda^m_{\bar \mu,\bar\nu_0,n})$ based at $l_{\bar\mu_0,\bar\nu_0}^s$, see Fig.\ \ref{fg_5_5}.
In particular, $\eta^s_{\bar \mu_0,\bar \nu_0}=r_{\bar \mu_0,\bar \nu_0}$.
Then, the velocity vector $\boldsymbol{v}_{\bar\mu_0}^{s,\bot}(r_{\bar\mu_0,\bar\nu_0})$ of $\eta^s_{\bar\mu,\bar\nu_0}$ at $r_{\bar\mu_0,\bar\nu_0}$ is defined by
$$\boldsymbol{v}_{\bar\mu}^{s,\bot}(r_{\bar\mu_0,\bar\nu_0})=\frac{d\eta^s_{\bar\mu,\bar\nu_0}}{d\bar\mu}(\bar\mu_0).$$
The velocity vector $\boldsymbol{v}^{s,\bot}_{\bar\nu}(r_{\bar\mu_0,\bar\nu_0})$ is defined similarly.
Let $\{a_{\bar\mu,\bar\nu_0,n}\}_n$ be a sequence of the intersection points of $\sigma$ and $W^s(P^+_{\bar \mu,\bar\nu_0,n})$, which $C^3$-converges to a point $a_{\bar\mu,\bar\nu_0}\in W^s(P^+_{\bar\mu,\bar\nu_0})\cap \sigma$ as $n\rightarrow \infty$.
Since $W^s(P^+_{\bar\mu,\bar\nu_0})$ is the horizontal line passing through $P^+_{\bar\mu,\bar\nu_0}$, the $\bar y$-entry of the coordinate of $a_{\bar \mu,\bar\nu_0}$ is $\bar y_1^+(\bar\mu,\bar\nu_0)$.
It follows from this fact and (\ref{y_1}) that
$$\boldsymbol{v}^{s,\bot}_{\bar\mu}(r_{\bar\mu_0,\bar\nu_0})\approx \frac{da_{\bar\mu,\bar\nu_0,n}}{d\bar \mu}(\bar\mu_0)\approx \frac{da_{\bar\mu,\bar\nu_0}}{d\bar \mu}(\bar\mu_0)=\Bigl(0,\frac{\partial\bar y_1^+}{\partial \bar\mu}(\bar\mu_0,\bar\nu_0)\Bigr)\approx \Bigl(0,\frac{1}{4}\Bigr),$$
where the first approximation is given by Accompanying Lemma (Lemma \ref{lem_1_1}).

Similarly, for any point $r$ in a $\Lambda_{\bar\mu,\bar\nu,n}^m$-leaf of $\mathcal{G}^{s(i),+}_{\bar\mu,\bar\nu,n}$ (resp.\ $\mathcal{G}^{s(i),-}_{\bar\mu,\bar\nu,n}$), Accompanying Lemma together with (\ref{y_1}) implies
\begin{equation}\label{eqn_sbot}
\boldsymbol{v}^{s,\bot}_{\bar\mu}(r)\approx \Bigl(0,\frac{1}{4}\Bigr)\ \bigl(\mbox{resp.}\ \Bigl(0,-\frac{1}{4}\Bigr)\bigr)\quad\mbox{and}\quad \boldsymbol{v}^{s,\bot}_{\bar\nu}(r)\approx \Bigl(0,\frac{1}{10}\Bigr).
\end{equation}

\begin{figure}[hbtp]
\centering
\scalebox{0.7}{\includegraphics[clip]{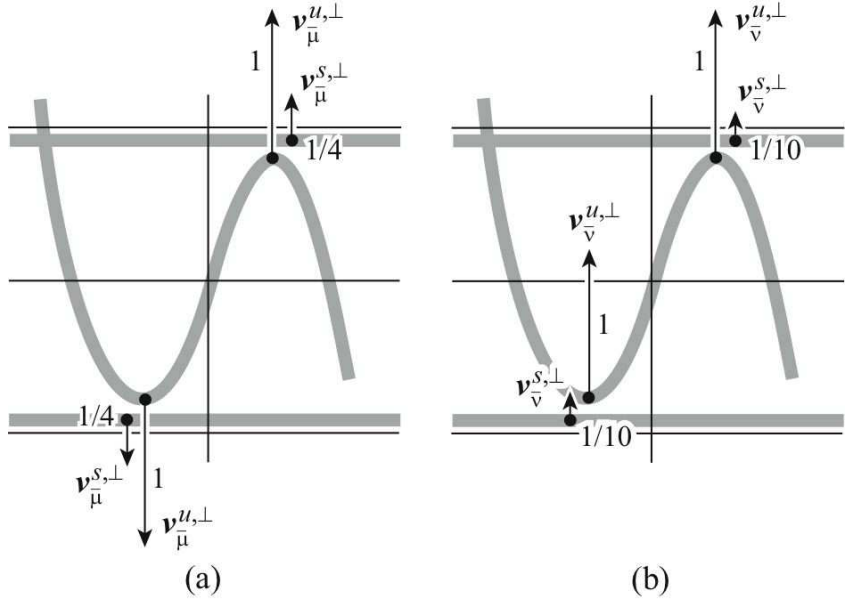}}
\caption{ The numbers $1,1/4,1/10$ associated to the vectors are approximate values of their lengths.}
\label{fg_5_6}
\end{figure}

For any point $r_{\bar\mu_0,\bar\nu_0}$ of $\Lambda_{\mu_0,\nu_0}^{\mathrm{out}}$-leaf $l_{\bar\mu_0,\bar\nu_0}^u$ of $\mathcal{F}^{u(i)}_{\mu_0,\nu_0}$, the velocity vector $\boldsymbol{v}^{u,\bot}_{\bar\mu}(r_{\bar\mu_0,\bar\nu_0})$ (resp.\ $\boldsymbol{v}^{u,\bot}_{\bar\nu}(r_{\bar\mu_0,\bar\nu_0})$) can be defined in the same manner as above by using the intersection point $\eta_{\bar\mu,\bar\nu_0}^u$ (resp.\  $\eta_{\bar\mu_0,\bar\nu}^u$) of a short vertical segment $\sigma$ passing through $r_{\bar\mu_0,\bar\nu_0}$ and a $\bar\mu$-continuation $l_{\bar\mu,\bar\nu_0}^u$ (resp.\ a $\bar\nu$-continuation $l_{\bar\mu_0,\bar\nu}^u$) based at $l_{\bar\mu_0,\bar\nu_0}^u$.
Let $b_{\bar\mu,\bar\nu,n}^+$ (resp.\ $b_{\bar\mu,\bar\nu,n}^-$) be the maximal (resp.\ minimal) point of the highest (resp.\ lowest) leaf of $\mathcal{F}^{u(i)}_{\mu,\nu}$.
Let $r$ be any point in $\Lambda_{\mu,\nu}^{\mathrm{out}}$-leaf of $\mathcal{F}^{u(i)}_{\mu,\nu}$ near $b_{\bar\mu,\bar\nu,n}^+$ (resp.\ $b_{\bar\mu,\bar\nu,n}^-$).
Again, by Accompanying Lemma together with (\ref{y_2}),
\begin{equation}\label{eqn_ubot}
\boldsymbol{v}^{u,\bot}_{\bar\mu}(r), \boldsymbol{v}^{u,\bot}_{\bar\nu}(r)\approx (0,1)\quad\Bigl(\mbox{resp.}\ \boldsymbol{v}^{u,\bot}_{\bar\mu}(r)\approx (0,-1),\ \boldsymbol{v}^{u,\bot}_{\bar\nu}(r)\approx (0,1)\Bigr).
\end{equation}
Figure \ref{fg_5_6} illustrates the approximations (\ref{eqn_sbot}) and (\ref{eqn_ubot}).

\begin{lemma}\label{lem_5_2}
For a sufficiently small $\varepsilon>0$, suppose that $\bar c:(-\varepsilon,\varepsilon)\longrightarrow \mathbb{R}^2$ is the parametrization of the curve in the $\bar\mu\bar\nu$-plane such that $\bar c(t)=(\bar\mu,t+\bar\nu)$.
Let $r$ be a tangency of a $\Lambda_{c(t_0)}^{\mathrm{out}}$-leaf of $\mathcal{F}^{u(i)}_{c(t_0)}$ and a $\Lambda_{\bar c(t_0),n}^m$-leaf of $\mathcal{G}^{s(i),\pm}_{\bar c(t_0),n}$ for some $t_0\in (-\varepsilon,\varepsilon)$.
Then, the following {\rm (i)} and {\rm (ii)} hold.
\begin{enumerate}[\rm (i)]
\item
If the position of $r$ is sufficiently close to $b_{\bar c(t_0),n}^+$, then $r$ is a contact-making tangency with respect to $\{\psi_{\bar c(t),n}\}$.
\item
If the position of $r$ is sufficiently close to $b_{\bar c(t_0),n}^-$, then $r$ is a contact-breaking tangency with respect to $\{\psi_{\bar c(t),n}\}$.
\end{enumerate}
\end{lemma}
\begin{proof}
(i)
Suppose that $r$ is sufficiently close to $b_{\bar c(t_0),n}^+$.
If necessary reparametrizing $c$, one can suppose that $t_0=0$ and hence $c(t_0)=c(0)=(\bar\mu,\bar\nu)$.
By (\ref{eqn_sbot}) and (\ref{eqn_ubot}), $\boldsymbol{v}^{u,\bot}_{\bar\nu}(r)-\boldsymbol{v}^{s,\bot}_{\bar\nu}(r)\approx (0,9/10)$.

For any sufficiently small $\Delta t>0$, let $\hat\sigma$ be a short vertical segment in the $\bar x\bar y$-plane such that, for the $t$-continuations $\{\hat\eta^u_{\bar c(t)}\}$, $\{\hat\eta^s_{\bar c(t)}\}$ with $\hat\eta^u_{\bar c(t)}\in \hat \sigma\cap l_{\bar c(t)}^u$, $\hat\eta^s_{\bar c(t)}\in \hat \sigma\cap l_{\bar c(t)}^s$, $\mathrm{dist}(\hat\eta^u_{\bar c(\Delta t)},\hat\eta^s_{\bar c(\Delta t)})$ has the maximum value among all short vertical segments near $\sigma$.
Since leaves of $\mathcal{F}^{u(i)}_{c(t)}$ and $\mathcal{G}^{s(i),+}_{\bar c(t),n}$ vary $C^1$ with respect to $t$ and the $\sigma$-direction, the move from $l^u_{\bar c(t)}$ to $l^u_{\bar c(t+\Delta t)}$ and that from $l^s_{\bar c(t)}$ to $l^s_{\bar c(t+\Delta t)}$ are well approximated by parallel translations in a small neighborhood of $r$.
This implies that $\hat{\boldsymbol{v}}^{u,\bot}_{\bar\nu}(\hat r^u)-\hat{\boldsymbol{v}}^{s,\bot}_{\bar\nu}(\hat r^u)$ is also well approximated by $(0,9/10)$, where $\hat{\boldsymbol{v}}^{u,\bot}_{\bar\nu}(\hat r^u)$, $\hat{\boldsymbol{v}}^{s,\bot}_{\bar\nu}(\hat r^s)$ are the velocity vectors of $\{\hat\eta^u_{\bar c(t)}\}$ and $\{\hat\eta^s_{\bar c(t)}\}$ at $\hat r^u=\hat\eta^u_{\bar c(0)}$ and $\hat r^s=\hat\eta^s_{\bar c(0)}$ respectively.
Here, we note that $\mathrm{dist}(\sigma,\hat\sigma)=O(\Delta t)$ and hence $\mathrm{dist}(\hat r^u,\hat r^s)=O(\Delta t^2)$.
Thus, the position of $\hat \eta_{\bar c(\Delta t)}^u$ is higher than that of $\hat \eta_{\bar c(\Delta t)}^s$ by $9\Delta t/10$ 
approximately, see Fig.\ \ref{fg_5_7}\,(a).
\begin{figure}[hbtp]
\centering
\scalebox{0.6}{\includegraphics[clip]{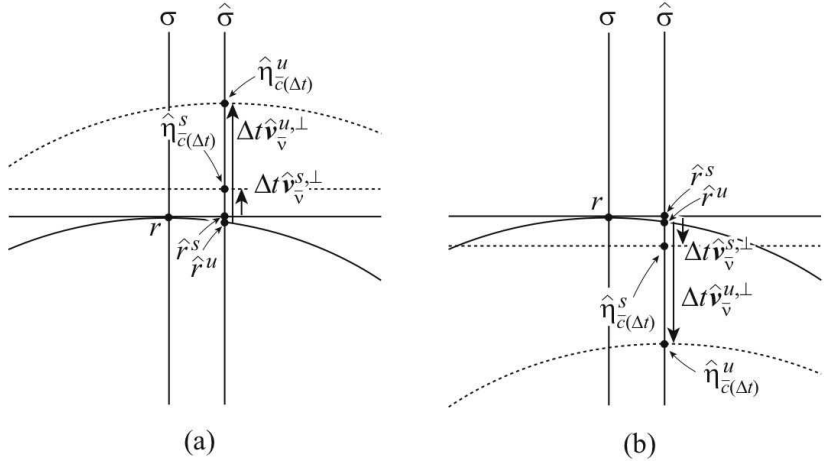}}
\caption{(a) The case of $\Delta t>0$.
(b) The case of $\Delta t<0$.}
\label{fg_5_7}
\end{figure}

As was seen in Subsection \ref{p_tangency}, $r$ is a quadratic tangency and hence we may assume that the curvature of $l^u_{\bar c(0)}$ at $r$ is greater than that of $l^s_{\bar c(0)}$.
From the continuity of the curvatures of leaves of compatible foliations, one can suppose that the infimum of the curvatures of $l_{\bar c(\Delta t)}^u$ in a small neighborhood $\mathcal{N}(r)$ of $r$ is greater than the supremum of those of $l_{\bar c(\Delta t)}^s$ for any sufficiently small $\Delta t>0$.
This shows that $l_{\bar c(\Delta t)}^u$ meets $l_{\bar c(\Delta t)}^s$ transversely in $\mathcal{N}(r)$.

Now, we consider the case in which $\Delta t<0$.
By the argument similar to that as above, it is shown that $l_{\bar c(\Delta t)}^u\cap l_{\bar c(\Delta t)}^s\cap \mathcal{N}(r)=\emptyset$ when $\vert\Delta t\vert$ is sufficiently small, see Fig.\ \ref{fg_5_7}\,(b).
Here, a short vertical segment $\hat \sigma$ is taken so that $\mathrm{dist}(\hat\eta^u_{\bar c(\Delta t)},\hat\eta^s_{\bar c(\Delta t)})$ has the minimum value among all short vertical segments near $\sigma$.

It follows that $r$ is a contact-making tangency with respect to $\{\psi_{\bar c(t),n}\}$.

(ii)
Next, we suppose that $r$ is close to $b_{\bar c(t_0),n}^-$.
The situation in the case of $\Delta t>0$ (resp.\ $\Delta t<0$) is shown by turning Fig.\ \ref{fg_5_7}\,(b) (resp.\ Fig.\ \ref{fg_5_7}\,(a)) upside down.
This proves that $r$ is a contact-breaking tangency with respect to $\{\psi_{\bar c(t),n}\}$.
\end{proof}

\subsection{Proof of Theorem \ref{main_a}}

With the notation and definitions as above, we prove Theorem \ref{main_a}.

For a fixed $\bar\mu_0<3$ with $(\bar\mu_0,0)\in \mathcal{N}$, one can choose $n_0$ and $i_0$ so that $b_{\bar\mu_0,0,n}^+$ is lower than $\mathcal{G}^{s(i),+}_{\bar\mu_0,0,n}$ and $b_{\bar\mu_0,0,n}^-$ is higher than $\mathcal{G}^{s(i),-}_{\bar\mu_0,0,n}$ if $n\geq n_0$ and $i\geq i_0$.
The relations (\ref{eqn_sbot}) and (\ref{eqn_ubot}) guarantee that one can have the first value $\bar\nu_0$ by increasing $\bar\nu$ from zero such that the highest leaf of $\mathcal{F}^{u(i)}_{\mu_0,\nu_0}$ touches that of $\mathcal{G}^{s(i),+}_{\bar\mu_0,\bar\nu_0,n}$ near $b_{\bar\mu_0,\bar\nu_0,n}^+$.

Consider a continuous function $\bar f(\bar \mu)$ defined on $\bar\mu\geq \bar \mu_0$ with $\bar f(\bar\mu_0)=\bar\nu_0$ such that the highest leaves of $\mathcal{G}^{s(i),+}_{\bar\mu,\bar f(\bar\mu),n}$ and $\mathcal{F}^{u(i)}_{\mu,f(\mu)}$ keep the tangency.

\begin{lemma}\label{lem_5_3}
The function $\bar f$ is strictly decreasing.
\end{lemma}
\begin{proof}
By (\ref{eqn_sbot}) and (\ref{eqn_ubot}), for any sufficiently small $\Delta \bar \mu,\Delta \bar \nu$, the height of the highest leaf of $\mathcal{F}^{u(i)}_{\mu+\Delta \mu,\nu+\Delta \nu}$ (resp.\ $\mathcal{G}^{s(i),+}_{\bar\mu+\Delta\bar\mu,\bar \nu+\Delta\bar\nu,n}$) in a small neighborhood of $b^+_{\bar\mu,\bar\nu,n}$ relative to that of $\mathcal{F}^{u(i)}_{\mu,\nu}$ (resp.\ $\mathcal{G}^{s(i),+}_{\bar\mu,\bar \nu,n}$) is $\Delta \bar\mu+\Delta\bar\nu$ (resp.\ $\Delta\bar\mu/4+\Delta\bar\nu/10$) approximately.
Then, from the definition of $\bar f$,
$$\Delta\bar\mu+\Delta\bar f(\bar \mu)\approx \frac{1}{4}\Delta \bar\mu+\frac{1}{10}\Delta \bar f(\bar\mu),$$
where $\Delta\bar f(\bar\mu)=\bar f(\bar\mu+\Delta\bar\mu)-\bar f(\bar\mu)$.
It follows that $\Delta \bar f(\bar\mu)\approx -5\Delta\bar\mu/6<0$ for any $\Delta \bar\mu>0$.
Thus, $\bar f$ is strictly decreasing.
\end{proof}

\begin{proof}[Proof of Theorem \ref{main_a}]
By Lemma \ref{lem_5_3}, the level of $b_{\bar\mu,\bar f(\bar \mu),n}^-$ is also strictly decreasing.
More strictly, the proof of the lemma implies 
that the difference of the level is approximated by $-\Delta\bar\mu-5\Delta\bar\mu/6=-11\Delta\bar\mu/6$.
Thus, we have $(\bar\mu_1,\bar\nu_1)=(\bar\mu_1,\bar f(\bar\mu_1))\in \mathcal{N}$ such that 
not only the highest leaves of $\mathcal{G}^{s(i),+}_{\bar\mu_1,\bar\nu_1,n}$ 
and $\mathcal{F}^{u(i)}_{\mu_1,\nu_1}$ 
but also the lowest leaves of $\mathcal{G}^{s(i),-}_{\bar\mu_1,\bar\nu_1,n}$ and $\mathcal{F}^{u(i)}_{\mu_1,\nu_1}$ admit tangencies, see Fig.\ \ref{fg_5_8}.
\begin{figure}[hbt]
\begin{center}
\scalebox{0.85}{\includegraphics[clip]{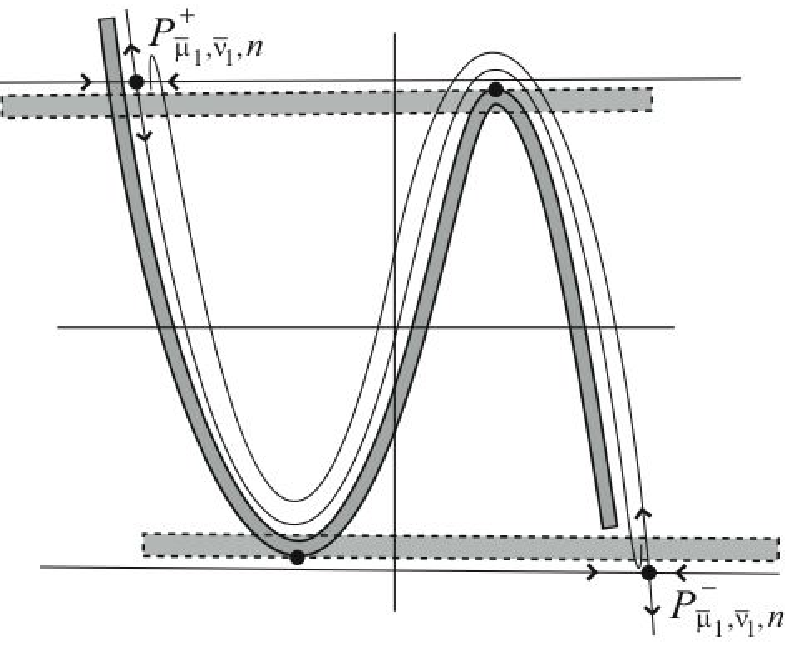}}
\caption{}
\label{fg_5_8}
\end{center}
\end{figure}

By Lemma \ref{lem_2_3} and (\ref{eqn_2}),
$$
\tau(\Lambda_{\mu_1,\nu_1}^{\mathrm{out}}\cap W^s(p_{\mu_1,\nu_1}))\cdot \tau^m_{\bar\mu_1,\bar\nu_1,n}>1.$$
As was seen in Subsection \ref{p_tangency}, the standard theory of Newhouse \cite{N1} and Palis-Takens \cite{PT} implies the existence of an open neighborhood $\mathcal{N}_n\subset \mathcal{N}$ of $(\bar\mu_1,\bar\nu_1)$ such that, for any $(\bar\mu,\bar\nu)\in \mathcal{N}_n$, $\mathcal{G}^{s(i),+}_{\bar\mu,\bar\nu,n}$ and $\mathcal{F}^{u(i)}_{\mu,\nu}$ have respectively a $\Lambda^m_{\bar\mu,\bar\nu,n}$-leaf and a $\Lambda^{\mathrm{out}}_{\mu,\nu}$-leaf admitting a quadratic tangency $\hat b_{\bar\mu,\bar \nu,n}^+$ near $b_{\bar\mu,\bar \nu,n}^+$, and $\mathcal{G}^{s(i),-}_{\bar\mu,\bar\nu,n}$ and $\mathcal{F}^{u(i)}_{\mu,\nu}$ also have respectively a $\Lambda^m_{\bar\mu,\bar\nu,n}$-leaf  and a $\Lambda^{\mathrm{out}}_{\mu,\nu}$-leaf admitting a quadratic tangency $\hat b_{\bar\mu,\bar \nu,n}^-$ near $b_{\bar\mu,\bar \nu,n}^-$.
For any $(\bar\mu,\bar\nu)\in \mathcal{N}_n$, consider the vertical curve $\bar c:(-\varepsilon,\varepsilon)\longrightarrow \mathcal{N}_n$ defined by $\bar c(t)=(\bar\mu,t+\bar\nu)$ for a small $\varepsilon>0$.
By Lemma \ref{lem_5_2}, for any $t\in (-\varepsilon,\varepsilon)$, $\hat b_{\bar c(t),n}^+$ is a contact-making tangency and $\hat b_{\bar c(t),n}^-$ is a contact-breaking tangency with respect to $\psi_{\bar c(t),n}$.

By Lemma \ref{lem_3_1}\,(ii), the union $\mathcal{O}=\bigcup_{n=n_0}^\infty \Theta_n(\mathcal{N}_n)$ is an open set in the $\mu\nu$-plane with $\mathrm{Cl}(\mathcal{O})\ni (0,0)$.
For any $(\mu,\nu)\in \Theta_n(\mathcal{N}_n)$, let $c:(-\varepsilon,\varepsilon)\longrightarrow \mathcal{O}$ be the curve defined by $c(t)=\Theta_n(\bar c(t))$.
By the relation (\ref{eqn_3}), tangencies of $\psi_{\bar c(t),n}$ are transformed to those of $\varphi_{c(t)}$.
Thus, the one-parameter family $\{\varphi_{c(t)}\}$ exhibits both persistent contact-making and contact-braking tangencies.
This completes the proof.
\end{proof}

\section{Proof of Theorem \ref{main_b}}\label{proofB}

In this section, we prove the second main theorem.
Our proof is based on results in Wang-Young \cite{WY}.

\subsection{Strange attractors and SRB measures}\label{strange_attractor}

First of all, we recall the definitions of strange attractors and SRB measures.

\begin{definition}
An invariant set $\Omega$ of a $2$-dimensional diffeomorphism is called a \emph{strange attractor} if it satisfies the following conditions.
\begin{itemize}
\item
There exists a saddle point $p\in \Omega$ such that the unstable manifold $W^u(p)$ has dimension $1$ and $\mathrm{Cl}(W^u(p))=\Omega$.
\item
There exists an open neighborhood $U$ of $\Omega$ such that $\{f^n(U)\}_{n=1}^\infty$ is a decreasing sequence with $\Omega=\bigcap_{n=1}^\infty f^n(U)$.
\item
There exists a point $z_0\in \Omega$ whose positive orbit is dense in $\Omega$ and a non-zero vector $v_0\in T_{z_0}(\mathbb{R}^2)$ with
\begin{equation}\label{Lyapunov_exponent}
\Vert d\varphi^n_{z_0}(v_0)\Vert \geq e^{cn}\Vert v_0\Vert
\end{equation}
for any integer $n\geq 0$ and some constant $c>0$.
\end{itemize}
\end{definition}

We say that $\varphi$ has a \emph{positive Lyapunov exponent} at $z_0$ if there exists a non-zero vector $v_0\in T_{z_0}(\mathbb{R}^2)$ satisfying the inequality (\ref{Lyapunov_exponent}).
A $\varphi$-invariant, Borel probability measure $\boldsymbol{\mu}$ with $\mathrm{supp}(\boldsymbol{\mu})=\Omega$ is called an \emph{SRB measure} supporting $\Omega$ if $\varphi$ has a positive Lyapunov exponent $\boldsymbol{\mu}$-almost everywhere and the conditional measures of $\boldsymbol{\mu}$ on unstable manifolds are absolutely continuous with respect to the Lebesgue measure on these manifolds.

\subsection{Wang-Young's Theorems for the existence of strange attractors}

For closed intervals $I,J$, consider a two-parameter family of $C^4$-diffeomorphisms $T_{\alpha,\beta}:[-1,1]\times I\longrightarrow \mathbb{R}^2$ $C^4$-depending on the parameters $(\alpha, \beta)\in J\times (0,1]$.
Furthermore, we suppose that $T_{\alpha,\beta}$ has the form
\begin{equation}\label{T_ab}
T_{\alpha,\beta}:(x,y)\longmapsto \bigl(\beta u(x,y,\alpha,\beta),F(x,y,\alpha)+\beta v(x,y,\alpha,\beta)\bigr),
\end{equation}
where the $C^3$-norms of $u$, $v$ are uniformly bounded as maps of $x,y,\alpha$ for all $\beta\in (0,1]$.
Let $F_{\alpha}:I\rightarrow \mathbb{R}$ be the function defined by $F_{\alpha}(y)=F(0,y,\alpha)$ and $C_\alpha$ the critical point set $\{y\in I;F_{\alpha}'(y)=0\}$ of $F_\alpha$.
In the case that $F_{\alpha}(y)=1-\alpha y^{2}$, the family $\{T_{\alpha,\beta}\}$ is well known as a quadratic H\'enon-like family, see \cite{BC,MV} and \cite[\S 6.3]{PT}.
When  $\{F_{\alpha}\}$ is a family of cubic polynomials, we say that the family $\{T_{\alpha,\beta}\}$ is \textit{cubic polynomial-like}.

Now, we consider the case in which $F_\alpha (I)\subset I$ and 
$C_\alpha$ consists of continuations $c_\alpha^{(1)},\dots,c_\alpha^{(n)}$.
Then, the complement $I\setminus C_\alpha$ consists of the continuations of the intervals $I_\alpha^{(0)},I_\alpha^{(1)},\dots,I_{\alpha}^{(n)}$.
We say that $p\in I$ has the {\it same symbolic dynamics} with respect to $F_\alpha$ as that of $q\in I$ with respect to $F_{\alpha_*}$ if they have the same itinerary.
That is, for any integer $k\geq 0$, $F^k_\alpha (p)\in I_\alpha^{(i)}$ (resp.\ $F^k_\alpha (p)=c_\alpha^{(j)}$) if and only if $F^k_{\alpha_*}(q)\in I_{\alpha_*}^{(i)}$ (resp.\ $F^k_{\alpha_*} (q)=c_{\alpha_*}^{(j)}$).

\begin{definition}\label{def_6_1}
We say that the family $\{T_{\alpha,\beta}\}$ satisfies the \emph{Wang-Young condition} if the following three conditions {\rm (a)-(c)} hold for some $\alpha_*\in J$.
These conditions correspond to the Steps I, II and IV in \cite[\S 1.1]{WY} respectively (see also \cite[\S 1.5]{WY}).
(The two parameter family $T_{\alpha,\beta}$ as above satisfies the condition corresponding to the Step III in \cite{WY} by itself.)
\begin{enumerate}[\rm (a)]
\item
The image $F_{\alpha_*}(I)$ is contained in the interior $\mathrm{Int} (I)$, and $F_{\alpha_*}$ satisfies the {\it Misiurewicz conditions} on $I$, that is, (i) $F_{\alpha_*}''(y)\neq 0$ for all $y\in C_{\alpha_*}$, (ii) $F_{\alpha_*}$ has negative Schwarzian derivative on $I\setminus C_{\alpha_*}$, (iii) there exists no element $y\in I$ with $F_{\alpha_*}^m(y)=y$ and $|(F_{\alpha_*}^m)'(y)|\leq 1$ for any $m>0$ and (iv) $\inf_{m>0}\mathrm{dist}(F_{\alpha_*}^m(y),C_{\alpha_*})>0$ for all $y\in C_{\alpha_*}$.
\item
The one-parameter family $\{F_\alpha\}$ satisfies the {\it transversality condition} at $\alpha=\alpha_*$, that is,
$$
\frac{d}{d\alpha}p(\alpha_*)\neq \frac{d}{d\alpha}F_\alpha(y(\alpha))\Bigr|_{\alpha=\alpha_*}
$$
for any continuation of critical points $y(\alpha)\in C_\alpha$, where $p(\alpha)$ is a unique point of $I$ which has the same symbolic dynamics with respect to $F_\alpha$ as that of $F_{\alpha_*}(y(\alpha_*))$ with respect to $F_{\alpha_*}$.
\item
The function $F_{\alpha_*}$ satisfies the {\it non-degeneracy condition}, that is,
$$
\frac{\partial}{\partial x} F(0,y,\alpha_*)\neq 0 
\quad \mbox{whenever}\quad F^{\prime}_{\alpha_*}(y)=0.
$$
\end{enumerate}
\end{definition}

Main theorems in Wang-Young \cite{WY} are summarized as follows.

\begin{theorem}[\cite{WY}]\label{thm_6_2}
If a two-parameter family $T_{\alpha,\beta}$ as above satisfies the Wang-Young condition, then there exists a positive measure set $\Delta\subset J\times (0,1]$ satisfying the following conditions.
\begin{enumerate}[\rm (i)]
\item
For all sufficiently small $\beta>0$, $\Delta_\beta:=\{\alpha\in J;\,(\alpha,\beta)\in \Delta\}$ has positive $1$-dimensional Lebesgue measure.
\item
For any $(\alpha,\beta)\in \Delta$, $T_{\alpha,\beta}$ exhibits  a strange attractor $\Omega_{\alpha,\beta}$ supported by an SRB measure.
\end{enumerate}
\end{theorem}

\begin{remark}
In fact, Theorems 1.1-1.7 in \cite{WY} imply that the strange attractors $\Omega_{\alpha,\beta}$ have extra properties other than those in Subsection \ref{strange_attractor}.
\end{remark}

\subsection{Outline of the proof of Theorem \ref{main_b}}\label{outline2}

We work with the notation and terminologies given in Subsection \ref{renorm}, but consider parameters $(\bar \mu,\bar \nu)$ only contained in $\bar \Sigma_0=J\times [0,1]$, where $J=[0,3]$.

Let $n_0>0$ be an integer such that Lemma \ref{lem_3_1} holds for any $n\geq n_0$.
Conjugating again the conjugated diffeomorphism $\psi_{\bar\mu,\bar\nu,n}=\Phi_{n}^{-1}\circ\varphi_{\Theta_{n}(\bar\mu,\bar\nu)}^{N+n}\circ \Phi_{n}(\bar{x},\bar{y})$ in (\ref{renormal}) for any $n\geq n_0$, we have a two-parameter family of diffeomorphisms $\tilde \psi_{\bar \mu,\bar \nu,n}$ on $\mathbb{R}^2$ satisfying the following conditions.
\begin{itemize}
\item
There exists a cubic polynomial-like two-parameter family $\{T_{\bar\mu,\beta}\}$ and a constant $0<\xi<1$ independent of $\bar\mu$ such that $T_{\bar\mu,\xi^n}=\tilde \psi_{\bar\mu,\xi^n,n}$ for any $n\geq n_0$, see Equation (\ref{T=psi}).
\item
The family $\{T_{\bar\mu,\beta}\}$ satisfies the Wang-Young condition, see Lemma \ref{lem_6_6}.
\end{itemize}
By Theorem \ref{thm_6_2}, there exists an integer $n_1\geq n_0$ such that, for any $n\geq n_1$, the 1-dimensional Lebesgue measure of $\Delta_{\xi^n}$ is positive, and for any $(\mu,\nu)$ corresponding to an element of $\Delta_{\xi^n}$, $\varphi_{\mu,\nu}$ has a strange attractor supported by an SRB measure.
However, since the 2-dimensional Lebesgue measure of the union $\bigcup_{n\geq n_1}\Delta_{\xi^n}$ in the $\bar\mu\beta$-plane is zero, this fact does not imply Theorem \ref{main_b} directly. 
To avoid this difficulty, we will apply the \emph{sliding trick}.
That is, for any $s\in [1,2]$, we define an $s$-dependant two-parameter family $\{T_{s;\bar\mu,\beta}\}$ satisfying the Wang-Young condition and such that $T_{s;\bar\mu,\xi^n}=\tilde \psi_{\bar\mu,s\xi^n,n}$.
In particular, there exists an $s$-dependant integer $n_1(s)\geq n_0$ such that, for any $n\geq n_1(s)$, there is a subset $\Delta_{s,\xi^n}$ of $J=[0,3]$ with positive 1-dimensional Lebesgue measure for any element $\bar\mu$ of which $\tilde \psi_{\bar\mu,s\beta_n.n}$ has a strange attractor as above.
Thus, for any element $(\mu,\nu)$ of the union $\mathcal{Z}=\bigcup_{s\in [1,2]}\bigcup_{n=n_1(s)}^\infty\Theta_n(\Delta_{s,\xi^n}\times \{s\xi^n\})$, $\varphi_{\mu,\nu}$ has a strange attractor supported by an SRB measure.
Moreover, Fubini's Theorem shows that, for any open neighborhood $U$ of the origin in the $\mu\nu$-plane, $\mathcal{Z}\cap U$ has a positive measure.

\subsection{Reforming to Wang-Young's settings}\label{Reformation}

Consider the two-parameter family $\{\psi_{\bar\mu,\bar\nu,n}\}$ of $C^4$-diffeomorphisms defined in (\ref{eqn_3}) converging in $C^3$-topology to the endomorphisms given by 
$$\psi_{{\bar \mu}, {\bar \nu}}(\bar{x}, \bar{y})=(\bar{y},-\bar{y}^{3}+{\bar\mu} \bar{y}+ {\bar\nu})$$
for $(\bar \mu,\bar \nu)\in \bar \Sigma=[0,4]\times [-1,1]$.

\begin{proposition} \label{prop_6_4} 
The map $\psi_{\bar\mu,\bar\nu}$ is smoothly conjugate to 
the endomorphism given by 
$\tilde{\psi}_{\bar{\mu}, \bar{\nu}}(\tilde{x},  \tilde{y})=(0,  -\tilde{y}^{3}+\bar{\mu} \tilde{y}+ \bar{\nu}+\tilde{x})$. 
\end{proposition}
\begin{proof}
The diffeomorphisms $f_{\bar\mu,\bar\nu}$ of $\mathbb{R}^2$ defined by
$$
(\tilde{x}, \tilde{y})=f_{\bar\mu, \bar\nu}(\bar x,\bar y)=
(\bar x^{3}-\bar\mu\bar x -\bar\nu+\bar y, \bar x)
$$
accomplish our desired conjugation.
\end{proof}

For any $n\geq n_0$, let $\tilde \psi_{\bar \mu,\bar \nu,n}$ be the diffeomorphism of $\mathbb{R}^2$ obtained by conjugating $\psi_{\bar\mu,\bar\nu,n}=\Phi_{n}^{-1}\circ\varphi_{\Theta_{n}(\bar\mu,\bar\nu)}^{N+n}\circ \Phi_{n}$ by $f_{\bar \mu,\bar \nu}$.
Then, the new map has the form
$$
\tilde \psi_{\bar \mu,\bar \nu,n}(\tilde x,\tilde y)=\bigl(\tilde H_{n;1}(\bar \mu,\bar \nu,\tilde x,\tilde y), -\tilde y^3+\bar\mu \tilde y+\bar \nu +\tilde x+\tilde H_{n;2}(\bar \mu,\bar \nu,\tilde x,\tilde y)\bigr),
$$
where $\tilde H_{n;1}, \tilde H_{n;2}$ are $C^{4}$-functions.
The equations in (\ref{norm}) imply that, for any compact subset $K$ of $\mathbb{R}^2$,
\begin{align*}
\Vert \tilde H_{n;1}\Vert_{C^3,\bar \Sigma_0\times K}&=O(\sigma_{\bar \Sigma_0}^{-n/2})+O((\lambda\sigma)_{\bar \Sigma_0}^n),\\
\Vert \tilde H_{n;2}\Vert_{C^3,\bar \Sigma_0\times K}&=O(\sigma_{\bar \Sigma_0}^{-n/2})+O((\lambda\sigma)_{\bar \Sigma_0}^n).
\end{align*}
Thus, we have the following lemma.

\begin{lemma}\label{lem_6_5}
There exists a constant $\xi$ with $0<\xi<1$ independent of $(\bar\mu,\bar\nu)$ in $\bar \Sigma_0$ satisfying
$$\Vert \tilde H_{n;1}\Vert_{C^3,\bar \Sigma_0\times K}=O(\xi^n)\quad\mbox{and}\quad \Vert \tilde H_{n;2}\Vert_{C^3,\bar \Sigma_0\times K}=O(\xi^n).$$
\end{lemma}
\begin{proof}
The constant $\xi=\max\{\sigma_{\bar \Sigma_0}^{-1/2}, (\lambda\sigma)_{\bar \Sigma_0}\}$ satisfies the above equations.
\end{proof}

Fix an element $s$ of $[1,2]$, and define the $C^{4}$-functions $u_{s;n},v_{s;n}$ by
$$
u_{s;n}(\bar \mu,\tilde x,\tilde y)=\frac{1}{\xi^n}\tilde H_{n;1}(\bar \mu,s\xi^n,\tilde x,\tilde y),\quad
v_{s;n}(\bar \mu,\tilde x,\tilde y)=s+\frac{1}{\xi^n}\tilde H_{n;2}(\bar \mu,s\xi^n,\tilde x,\tilde y).
$$
Then, Lemma \ref{lem_6_5} implies
$$
\sup_n\{\Vert u_{s;n}\Vert_{C^3,J\times K}\}<\infty,\quad
\sup_n\{\Vert v_{s;n}\Vert_{C^3,J\times K}\}<\infty.
$$
Thus, for any closed interval $I$, one can define an $s$-dependent two-parameter family of $C^{4}$-diffeomorphisms $T_{s;\bar\mu,\beta}:[-1,1]\times I\longrightarrow \mathbb{R}^2$ such that
\begin{equation}\label{cubic diffeomorphism}
T_{s;\bar{\mu},\beta} : 
(\tilde{x},\tilde{y})\longmapsto 
\bigl(\beta u(\bar{\mu},\beta,\tilde{x},\tilde{y}),  -\tilde{y}^{3}+\bar{\mu} \tilde{y }+  \tilde{x}+\beta v(\bar{\mu},\beta,\tilde{x},\tilde{y}) \bigr)
\end{equation}
for $(\bar \mu,\beta)\in J\times (0,1]$, where $u,v$ are functions satisfying the following two conditions.
\begin{itemize}
\item
For any $n\geq n_0$,
$$u(\bar \mu,\xi^n,\tilde x,\tilde y)=u_{s;n}(\bar \mu,\tilde x,\tilde y),\quad
v(\bar \mu,\xi^n,\tilde x,\tilde y)=v_{s;n}(\bar \mu,\tilde x,\tilde y).$$
\item
The $C^3$-norms of $u,v$ as $(\bar \mu,\tilde x,\tilde y)$-functions on any compact subset of $J\times \mathbb{R}^2$ are bounded by a constant independent of $\beta\in (0,1]$.
\end{itemize}
Note that the former condition implies that
\begin{equation}\label{T=psi}
T_{s;\bar \mu,\xi^n}=\tilde \psi_{\bar \mu,s\xi^n,n}\quad\mbox{on}\quad [-1,1]\times I.
\end{equation}

Here, we set
$$
F_{\bar{\mu}}(\tilde{y})= -\tilde{y}^{3}+\bar{\mu} \tilde{y},
$$
and show that the family $\{T_{s;\bar\mu,\beta}\}$ satisfies the Wang-Young conditions as follows.

\begin{lemma}\label{lem_6_6}
There exists $\bar \mu_*\in (3\sqrt{3}/2,3)$ and a closed interval $I$
such that 
\begin{enumerate}[\rm (a)]
\item \label{W-Y condition 1}
$F_{\bar{\mu}_{\ast}}(I)\subset \mathrm{Int} (I)$ and $F_{\bar{\mu}_{\ast}}$ satisfies the Misiurewicz conditions on $I$.
\item \label{W-Y condition 2}
$\{ F_{\bar{\mu}} \}$ satisfies the transversality condition at $\bar{\mu}=\bar{\mu}_{\ast}$.
\item  \label{W-Y condition 4}
The function $F$ with $F(\tilde x,\tilde y,\bar \mu)=-\tilde y^3+\bar \mu\tilde y+\tilde x$ satisfies the non-degeneracy condition at $\bar\mu=\bar\mu_*$.
\end{enumerate}
\end{lemma}
\begin{proof}
(a)
The positive critical point $c(\bar\mu)$ of $F_{\bar\mu}$ is $\sqrt{\bar\mu/3}$.
Since $F^2_{\bar\mu}(c(\bar\mu))+\sqrt{\bar\mu}=-2+\sqrt{3}<0$ (resp.\ $3^{3/4}2^{-1/2}>0$) when $\bar\mu=3$ (resp.\ $\bar\mu=3\sqrt{3}/2$), by the Intermediate Value Theorem, there exists $\bar\mu_*\in (3\sqrt{3}/2,3)$ such that $F^2_{\bar\mu_*}(c(\bar\mu_*))=-\sqrt{\bar\mu_*}$ and hence $F^3_{\bar\mu_*}(c(\bar\mu_*))=0$.

The positive fixed point $e(\bar\mu_*)$ of $F_{\bar\mu_*}$ is $\sqrt{\bar\mu_*-1}$, and hence $F_{\bar\mu_*}^2(-c(\bar\mu_*))=\sqrt{\bar\mu_*}>e(\bar\mu_*).$
Note that both $F_{\bar\mu_*}^2(c(\bar\mu_*))=-\sqrt{\bar\mu_*}$, $F_{\bar\mu_*}(-c(\bar\mu_*))=-(2\bar\mu_*/3)\sqrt{\bar\mu_*/3}$ belong to $(-\infty,-c(\bar\mu_*)]$.
Since $F_{\bar\mu_*}(F^2_{\bar\mu_*}(c(\bar\mu_*)))=0<e(\bar\mu_*)<F_{\bar\mu_*}(F_{\bar\mu_*}(-c(\bar\mu_*)))$ and $F_{\bar\mu_*}|_{(-\infty,-c(\bar\mu_*)]}$ is strictly decreasing, there exists an $r'$ with
$$F_{\bar\mu_*}(-c(\bar\mu_*))<r'<F^2_{\bar\mu_*}(c(\bar\mu_*))\quad\mbox{and}\quad F_{\bar\mu_*}(r')=e(\bar\mu_*),$$
see Fig.\ \ref{fg_6_1}.
\begin{figure}[hbt]
\begin{center}
\scalebox{0.9}{\includegraphics[clip]{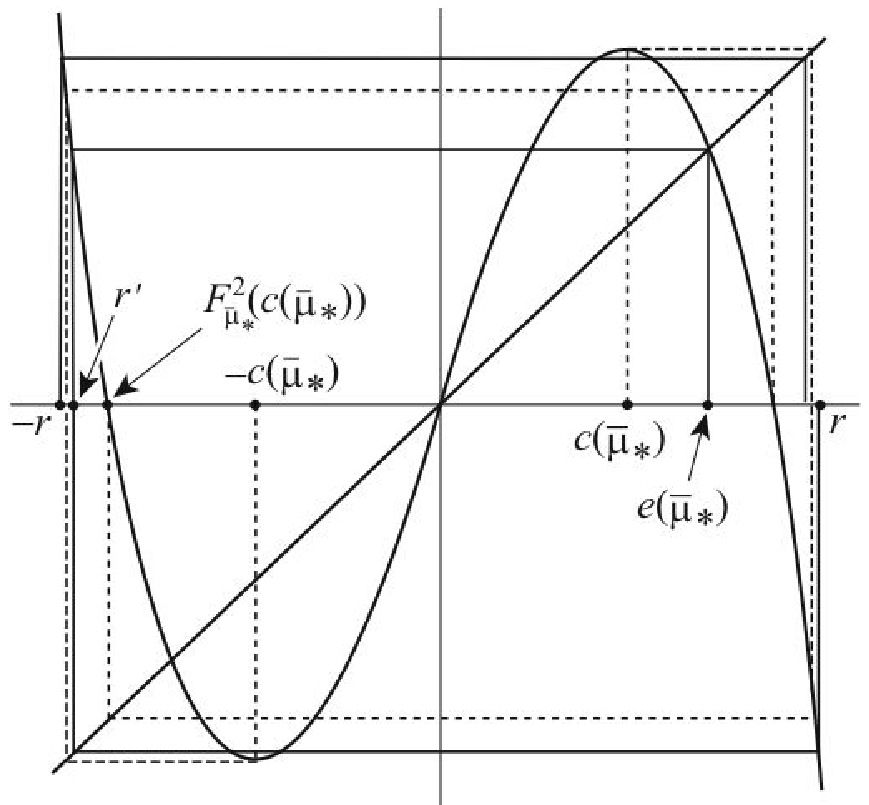}}
\caption{}
\label{fg_6_1}
\end{center}
\end{figure}
Since $r'<F_{\bar\mu_*}(F_{\bar\mu_*}(c(\bar\mu_*)))$, we have an $r>F_{\bar\mu_*}(c(\bar\mu_*))$ with $F_{\bar\mu_*}(r)=r'$.
In particular, $F^2_{\bar\mu_*}(r)=e(\bar\mu_*)$.
Since $F_{\bar\mu_*}$ is an odd function, the inequalities as above imply
$$-r<F_{\bar\mu_*}(-c(\bar\mu_*))<F_{\bar\mu_*}(r)<-c(\bar\mu_*)<c(\bar\mu_*)<F_{\bar\mu_*}(-r)<F_{\bar\mu_*}(c(\bar\mu_*))<r.$$
From this, we know that the closed interval $I=[-r,r]$ satisfies $F_{\bar\mu_*}(I)\subset \mathrm{Int}(I)$.

The direct calculation of the Schwarzian derivative of $F_{\bar{\mu}}$ shows
$$
SF_{\bar{\mu}}(\tilde{y})=\frac{F^{\prime\prime\prime}_{\bar{\mu}}(\tilde{y})}{F^{\prime}_{\bar{\mu}}(\tilde{y})}-\frac{3}{2}\left(\frac{F^{\prime\prime}_{\bar{\mu}}(\tilde{y})}{F^{\prime}_{\bar{\mu}}(\tilde{y})}\right)^{2} 
=-\frac{6\tilde{y}^{2}+\bar{\mu}}{(-3\tilde{y}^{2}+\bar{\mu})^{2}}<0
$$
for any $\tilde {y}\neq \pm c(\bar{\mu})$. 
Since $F^{3}_{\bar{\mu}_{\ast}}(-c(\bar{\mu}_{\ast}))=0$, any critical points of $F_{\bar{\mu}_{\ast}}$ can not return to a small neighborhood of the critical point set $C_{\bar\mu_*}=\{-c(\bar \mu_*), c(\bar \mu_*)\}$.
A well-known result by Singer on functions with negative Schwarzian derivative (see for example \cite[p.\ 155, Theorem 6.1]{DV}) implies that, if there existed a periodic point $q$ such that $F^{n}_{\bar{\mu}_{\ast}}(q)=q$ and 
$\vert (F^{n}_{\bar{\mu}_{\ast}})^{\prime}(q)\vert \leq 1$, then the basin of $q$ would contain either a critical point of $F_{\bar \mu_*}$ or a boundary point of $I$.
However, in our case, all the critical points $\pm c(\bar{\mu}_{\ast})$ and boundary points $\pm r$ are mapped  eventually to the repelling fixed points $\pm e(\bar\mu_*)$.
This implies that, for every periodic point $q$ with  $F^{n}_{\bar{\mu}_{\ast}}(q)=q$, $\vert (F^{n}_{\bar{\mu}_{\ast}})^{\prime}(q)\vert > 1$.  
Thus, we have proved that $F_{\bar{\mu}_{\ast}}$ satisfies the Misiurewicz conditions.

(b)
Since $p(\bar{\mu})$ has the same itinerary with respect to $F_{\bar{\mu}}$ as that of $F_{\bar{\mu}_{\ast}}(c(\bar{\mu}_{\ast}))$ with respect to $F_{\bar{\mu}_*}$, it follows that $p(\bar\mu)\in (c(\bar\mu),r]$, $F_{\bar\mu}(p(\bar\mu))\in [-r,-c(\bar\mu))$ and $F^k(p(\bar\mu))\in (-c(\bar\mu),c(\bar\mu))$ for any integer $k\geq 2$.
This implies that
$$F_{\bar\mu}^2(p(\bar\mu))=F_{\bar\mu}(p(\bar\mu))\bigl(-F_{\bar\mu}(p(\bar\mu))^2+\bar\mu\bigr)=0$$
for any $\bar\mu\approx \bar\mu_*$.
Since $F_{\bar\mu}(p(\bar\mu))<0$, it follows that $F_{\bar{\mu}}(p(\bar{\mu}))=-\sqrt{\bar{\mu}}$, that is,
$$
-(p(\bar{\mu}))^{3}+\bar\mu p(\bar{\mu})= -\sqrt{\bar{\mu}}.
$$
Then,
$$
\frac{d}{d\bar{\mu}}  p(\bar{\mu}) =\frac{2p(\bar{\mu})+1/\sqrt{\bar{\mu}}}{6(p(\bar{\mu}))^{2}-2\bar\mu}.
$$
Since 
$$p(\bar{\mu}_{\ast}) =F_{\bar{\mu}_{\ast}}(c(\bar{\mu}_{\ast}))=\frac{2\bar{\mu}_{\ast}}{3}\sqrt{\frac{\bar{\mu}_{\ast}}{3}},$$
we have
$$
\frac{d}{d\bar{\mu}} p(\bar{\mu}_{\ast})
=\frac{4\sqrt{3} \bar{\mu}_{\ast}^{2}+9}{2\bar\mu_*\sqrt{\bar{\mu}_{\ast}}(4\bar{\mu}_{\ast}^{2}-9)}.
$$
Since the function $h(t)=(4\sqrt{3} t^{2}+9)/ (2t\sqrt{t}(4t^{2}-9))$ is strictly decreasing on $3\sqrt{3}/{2}\leq t\leq  3$,
$$\frac{d}{d\bar \mu}p(\bar \mu_*)=h(\bar \mu_*)<h(3\sqrt{3}/2)<0.4.$$
On the other hand, $(d/d\bar{\mu}) F_{\bar{\mu}}(c(\bar{\mu})) =\sqrt{\bar{\mu}/{3}}$ is a strictly increasing function on $3\sqrt{3}/{2}\leq \bar{\mu}\leq  3$,
$$
\frac{d}{d\bar{\mu}}F_{\bar{\mu}} ( c(\bar \mu))\Bigr|_{\bar \mu=\bar \mu_*}>\frac{d}{d\bar{\mu}} F_{\bar{\mu}} (c(\bar \mu))\Bigr|_{\bar \mu=3\sqrt{3}/2} >0.9.
$$
This completes the  proof of (\ref{W-Y condition 2}).  

The assertion (c) is immediately derived from the form of $F(\tilde x,\tilde y,\bar\mu)$.
\end{proof}

\subsection{Proof of Theorem \ref{main_b}}

We still work with the notation as above.
In particular, $I$ is the closed interval given in Lemma \ref{lem_6_6}\,(a), and $J=[0,3]$.

\begin{proof}[Proof of Theorem \ref{main_b}]
By Lemma \ref{lem_6_6}, $\{T_{s;\bar{\mu},\beta}\}$ admits the Wang-Young condition for any $s\in [1,2]$.
Hence, there exists an integer $n_1=n_1(s)\geq n_0$ such that, for any $n\geq n_1$, $2\xi^n\leq 1$ and $\{T_{s;\bar \mu,\xi^n}\}_{0\leq \bar \mu \leq 3}$ satisfies the conclusions of Theorem \ref{thm_6_2}.
In particular, there exists a subset $\Delta_{s,\xi^n}$ of $J$ with positive 1-dimensional Lebesgue measure such that, for any $\bar \mu \in \Delta_{s,\xi^n}$, $T_{s;\bar \mu,\xi^n}=\tilde \psi_{\bar \mu,s\xi^n,n}$ has a cubic polynomial-like strange attractor $\Omega_{\bar \mu,s\xi^n,n}$ with an SRB measure, where the $n_1$ is chosen so that $\Omega_{\bar \mu,s\xi^n,n}$ is contained in a sufficiently small rectangle neighborhood $B$ of $\{0\}\times I$ in $\mathbb{R}^2$.
Consider the subset
$$\mathcal{Z}=\bigcup_{s\in [1,2]}\bigcup_{n=n_1(s)}^\infty \hat \Delta_{s,n}$$
of the $\mu\nu$-plane, where $\hat \Delta_{s,n}=\Theta_n(\Delta_{s,\xi^n}\times \{s\xi^n\})$.
Since $\tilde\psi_{\bar \mu,s\xi^n,n}$ is smoothly conjugate to $\varphi_{\Theta_{n}(\bar\mu,s\xi^n)}^{N+n}$, for any $(\mu,\nu)=\Theta_{n}(\bar\mu,s\xi^n)\in \mathcal{Z}$, $\varphi_{\mu,\nu}^{N+n}$ also has a strange attractor of the same type as that of $\Omega_{\bar \mu,s\xi^n,n}$.

Let $U$ be any open neighborhood of the origin in the $\mu\nu$-plane.
It remains to show that $\mathcal{Z}\cap U$ has positive Lebesgue measure.
Fix an integer $n_2\geq n_0$ with $\Theta_n(\bar \Sigma)\subset U$ for any $n\geq n_2$.
Note that, for any $n\geq \max\{n_1(s),n_2\}$,
\begin{equation}\label{meas}
\mathrm{meas}^{(1)}_{\,n}(\mathcal{Z}\cap\Gamma_{s,n})\geq \mathrm{meas}^{(1)}_{\,n}(\hat \Delta_{s,n})>0,
\end{equation}
where $\Gamma_{s,n}$ denotes the smooth curve $\Theta_n(J\times \{s\xi^n\})$ in $U$ and $\mathrm{meas}^{(1)}_{\,n}(\cdot)$ does the 1-dimensional Lebesgue measure on $\Gamma_{s,n}$, see Fig.\ \ref{fg_6_2}.
\begin{figure}[hbt]
\begin{center}
\scalebox{0.85}{\includegraphics[clip]{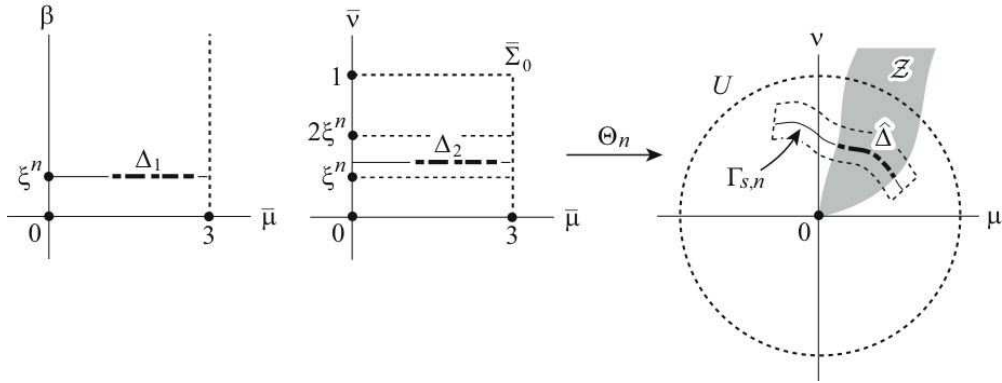}}
\caption{ $\Delta_1,\Delta_2,\hat \Delta$ represent respectively $\Delta_{s,\xi^n}\times \{\xi^n\}$, $\Delta_{s,\xi^n}\times \{s\xi^n\}$ and $\hat \Delta_{s,n}$.}
\label{fg_6_2}
\end{center}
\end{figure}

If the 2-dimensional Lebesgue measure $\mathrm{meas}^{(2)}(\mathcal{Z}\cap U)$ of $\mathcal{Z}\cap U$ were zero, then we would also have $\mathrm{meas}^{(2)}(\mathcal{Z}\cap R_n)=0$ for any $n\geq n_2$, where $R_n$ is the curvilinear rectangle $\bigcup_{1\leq s\leq 2}\Gamma_{s,n}$ in $U$.
By Fubini's Theorem, there would exist a zero-measure subset $E_n$ of $[1,2]$ such that $\mathrm{meas}^{(1)}_{\,n}(\mathcal{Z}\cap\Gamma_{s,n})=0$ for any $s\in [1,2]\setminus E_n$.
This implies that, for any element $s_0$ of the non-empty subset $[1,2]\setminus \bigcup_{n=n_2}^\infty E_n$ of $[1,2]$, $\mathrm{meas}^{(1)}_{\,n}(\mathcal{Z}\cap\Gamma_{s_0,n})=0$ if $n\geq n_2$, which contradicts (\ref{meas}).
Thus, $\mathrm{meas}^{(2)}(\mathcal{Z}\cap U)>0$.
This completes the proof of Theorem \ref{main_b}.
\end{proof}

\appendix

\section{Proof of Accompanying Lemma}

We will here prove Accompanying Lemma (Lemma \ref{lem_1_1}) stated in Subsection \ref{ac_lemma}.
Figure\ \ref{fg_4_5} in Subsection \ref{ac_lemma} illustrates the situation under which the lemma is considered.

First of all, we present an outline of the proof of this lemma.
Recall that $\{\varphi_t\}_{-\varepsilon \leq t\leq \varepsilon}$ is a one-parameter family in $\mathrm{Diff}^3(\mathbb{R}^2)$ and $\{\Lambda_t, p_t\}$ is a continuation of pairs of basic sets and saddle fixed points of $\varphi_t$.
Moreover, there exists a one-parameter family of another saddle fixed points $\hat p_t$ such that $W^s(\hat p_t)\setminus \{\hat p_t\}$ contains an arc $\alpha_t$ crossing the foliation $\mathcal{F}_t^{u(k_0)}=\varphi_t^{k_0}(\mathcal{F}_t^u)$ exactly for some integer $k_0>0$.
We need to see that how any continuation $\{\eta_t^u\}$ of crossing points of $\sigma$ and $\Lambda_t$-leaves of $\mathcal{F}_t^{u(k_0+j)}$ is represented.
The $\lambda$-lemma shows that, for all sufficiently large integer $j>0$, $\sigma^j_t=\varphi_t^{-j}(\sigma)$ contains an arc $\alpha_t^j$ $C^3$-converging to $\alpha_t$ as $j\rightarrow \infty$ and crossing $\mathcal{F}^{u(k_0)}_t$ exactly.
For such a large $j$, the $C^1$-diffeomorphisms $\hat \pi_t^{u(k_0);j}:I_t\longrightarrow \alpha_t^j$ can be defined as (\ref{pi_t}).
Then, we know that $\eta_t^u$ is an element of $\mathcal{X}_t=\varphi^j_t(\hat \pi_t^{u(k_0);j}(\Lambda_t\cap I_t))$.

Lemma \ref{lem_a_1} shows that, for any continuation $\{w_t\}$ with $w_t\in\Lambda_t\cap I_t$, the norm $\Vert  dw_t/dt\Vert$ of 
the velocity vector $dw_t/dt$ is bounded by a uniform constant.
This fact implies that, for any continuation of points $r_t^{\langle j;j\rangle}=\hat \pi_t^{u(k_0);j}(w_t)$ in $\hat \pi_t^{u(k_0);j}(\Lambda_t\cap I_t)$, the norm of its velocity vector is also uniformly bounded.
Lemma \ref{lem_a_2} shows that the differential $d\varphi^j_t$ reduces the velocity vector $dr_t^{\langle j;j\rangle}/dt-da_t^j/dt$ of $r_t^{\langle j;j\rangle}$ relative to $a_t^{j}=\varphi^{-j}_t(a_t)$ so that the norm of the resultant vector is arbitrarily small if $j$ is sufficiently large.
The actual proof is done by using the arc-length parametrizations $\beta^{(i)}_t$ of $\sigma^i_t$ based at $a_t^i=\varphi^{-i}_t(a_t)$ for $i=0,1,\dots,j$.
Accompanying Lemma is derived immediately from this lemma.

\begin{lemma}\label{lem_a_1}
There exists a constant $K>0$ independent of $t$ and the choice of a continuation $w_t\in \Lambda_t\cap I_t$ such that
$$
\Bigl\Vert \frac{dw_t}{dt}\Bigr\Vert\leq K.
$$
\end{lemma}
\begin{proof}
The space $C^0(\Lambda_0\cap I_0,\mathbb{R}^2)$ of continuous maps from $\Lambda_0\cap I_0$ to $\mathbb{R}^2$  
with the supremum norm is a Banach space.
By Proposition 1 in \cite{Pol}, the map $\theta:[-\varepsilon,\varepsilon]\rightarrow 
C^0(\Lambda_0\cap I_0,\mathbb{R}^2)$ defined by $(\theta (t))(w_0)=w_t$ is of $C^2$-class.
In particular, $d\theta /dt(t)$ is contained in $C^0(\Lambda_0\cap I_0,\mathbb{R}^2)$ and 
continuous on $t$.
Since  $[-\varepsilon,\varepsilon]\times (\Lambda_0\cap I_0)$ is compact and $(d\theta /dt(t))(w_0)=dw_t/dt(t)$, 
there exists a constant $K$ satisfying the above inequality.
\end{proof}

\begin{figure}[hbtp]
\centering
\scalebox{0.65}{\includegraphics[clip]{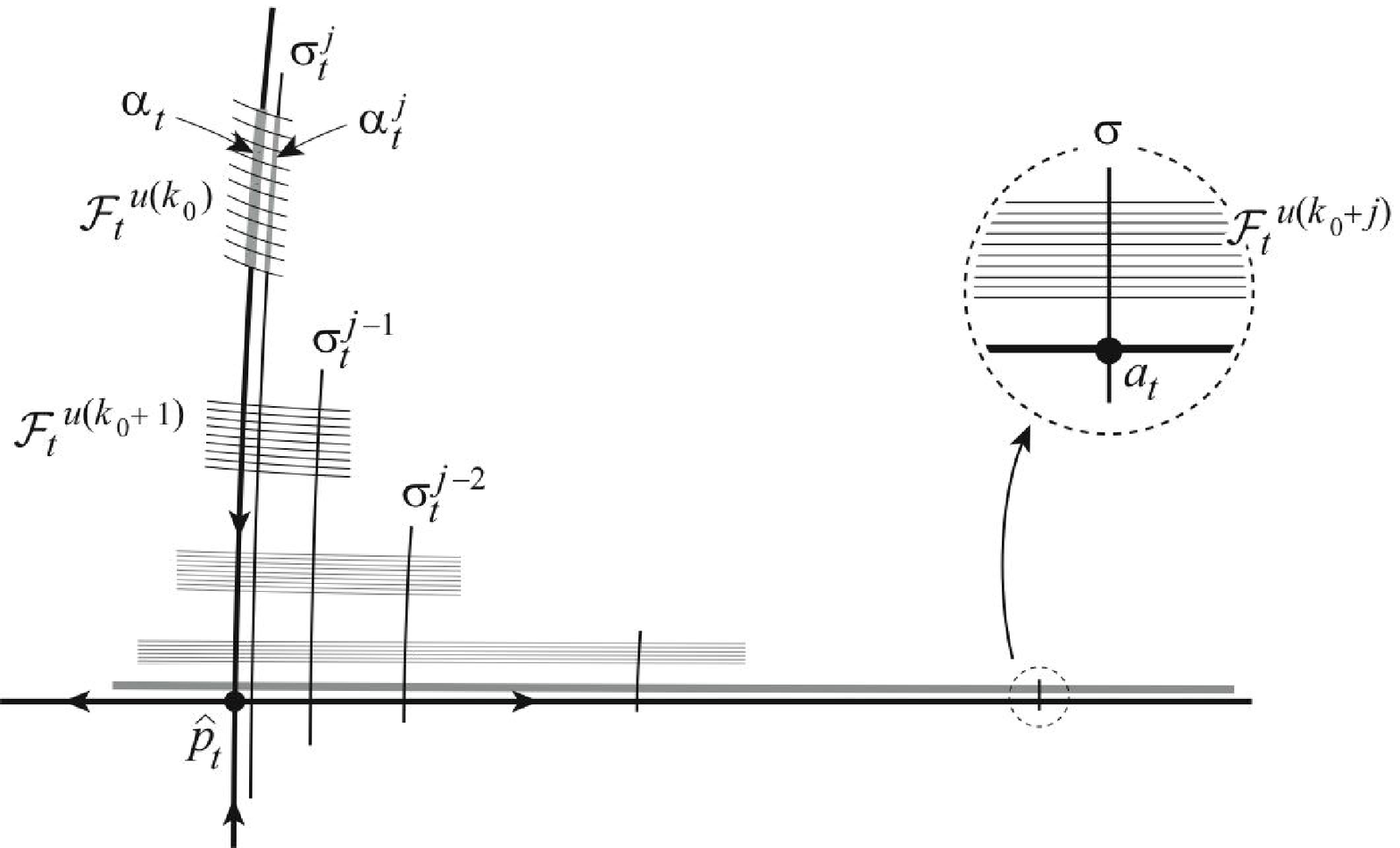}}
\caption{}
\label{fg_a_1}
\end{figure}

Recall that, for a $t_0\in (-\varepsilon,\varepsilon)$, $\sigma$ is a short segment in $\mathbb{R}^2$
 passing through $W^u(\hat p_{t_0})\setminus \{\hat p_{t_0}\}$ at a point $a_{t_0}$ almost orthogonally.
For any $t$ near $t_0$, $\sigma$ also meets $W^u(\hat p_t)$ almost orthogonally at a point 
$a_t$ close to $a_{t_0}$, see Fig.\ \ref{fg_a_1}.
Consider the segments $\sigma^j_t=\varphi_t^{-j}(\sigma)$ for $j=1,2,\ldots$.
Since $\hat p_t$ is a saddle point, the $\lambda$-lemma shows that, for any continuation of compact arcs $\tau_t$ in $W^s(\hat p_t)$, all but finitely many $\sigma^j_t$ have compact subarcs uniformly $C^3$-close to $\tau_t$.
The uniformity of the $C^3$-closeness is independent of $t\in [-\varepsilon,\varepsilon]$.

For $j=0,1,2,\ldots$, let $\beta^{(j)}_t:[u^j_{t,-},u^j_{t,+}]\longrightarrow \sigma^j_t$ be the arc-length parametrization of $\sigma^j_t$ with $\beta^{(j)}_t(0)=a_t^j\in W^u(\hat p_t)$, where $\sigma^0_t=\sigma$.
For $j\geq 1$, the composition
$$\zeta^{(j)}_t=(\beta^{(j-1)}_t)^{-1}\circ \varphi_t\circ \beta^{(j)}_t:[u^j_{t,-},u^j_{t,+}]\longrightarrow [u^{j-1}_{t,-},u^{j-1}_{t,+}]$$
is a $C^3$-diffeomorphism $C^3$-depending on $t$.
The 2-variable function $\zeta^{(j)}$ is defined by $\zeta^{(j)}(t,u)=\zeta^{(j)}_t(u)$.
Since $\zeta^{(j)}(t,0)=0$ for any $t$, $(\partial \zeta^{(j)}/\partial t)(t,0)=0$.
From this fact together with the uniform closeness on $\sigma^j_t$ as above, for any $\delta >0$, 
there exists a constant  $u_0>0$ and an integer $j_1>0$ such that 
\begin{equation}\label{eqn_4n}
\Bigl|\frac{\partial\zeta^{(j)}}{\partial t}(t,u)\Bigr|<\delta
\end{equation}
if $t\in [-\varepsilon,\varepsilon]$, $u\in [-u_0,u_0]$ and $j>j_1$. 
One can retake if necessary the $u_0$ so
that the inequality (\ref{eqn_4n}) holds also for the remaining $j = 1,\cdots, j_1$.

Since $\hat p_t$ is a saddle fixed point $C^3$-depending on $t$, one can have a constant $0<\lambda_0<1$ independent of all $t$ near $t_0$ such that
\begin{equation}\label{eqn_41}
\Bigl|\frac{d\zeta^{(j)}_t}{du}(u)\Bigr|<\lambda_0
\end{equation}
for all $j=1,2,\ldots$ and $u\in [-u_0,u_0]\cap [u^j_{t,-},u^j_{t,+}]$ if necessary replacing $u_0$ by a smaller positive number.
Usually, $u_0$ is less than $|u^j_{t,\pm}|$ if $j$ is sufficiently large, but greater if $j$ is small.

For all sufficiently large $j$, $\sigma^j_t$ contains an arc $\alpha_t^j$ crossing $\mathcal{F}_t^{u(k_0)}$ exactly which $C^3$-converges to $\alpha_t$ as $j\rightarrow \infty$ for any $t\in [-\varepsilon, \varepsilon]$.
Let $\hat\pi^{u(k_0);j}_t:I_t\longrightarrow \alpha_t^j$ be the $C^1$-diffeomorphism defined as (\ref{pi_t}), which $C^1$-converges to $\hat\pi^{u(k_0)}_t:I_t\longrightarrow \alpha_t$ as a $C^1$-embedding of $I_t$ to $\mathbb{R}^2$.
If necessary replacing $k_0$ by some greater integer, we may assume that $\hat\pi^{u(k_0);j}_t(I_t)\subset \beta_t^{(j)}([-u_0,u_0])$ for all sufficiently large $j$.

For any continuation of points $w_t\in \Lambda_t\cap I_t$, set
$$r^{\langle j; j\rangle}_t=\hat\pi^{u(k_0);j}_t(w_t)\in \alpha^j_t\quad\mbox{and}\quad u^{\langle j; j\rangle}_t=(\beta_t^{(j)})^{-1}(r_t^{\langle j;j\rangle})\in [-u_0,u_0].$$
Furthermore, we set $r^{\langle j-i;j\rangle}_t=\varphi^{i}_t(r^{\langle j;j\rangle}_t)$ and $u^{\langle j-i;j\rangle}_t=(\beta_t^{(j-i)})^{-1}(r^{\langle j-i;j\rangle}_t)$ for $i=1,2,\ldots,j$.
That is, $u_t^{\langle j-i;j\rangle}=\zeta_t^{(j-i+1)}\circ\cdots\circ \zeta_t^{(j)}(u^{\langle j;j\rangle}_t)$.
By Lemma \ref{lem_a_1} together with the $C^1$-convergence of $\hat\pi^{u(k_0);j}_t$ to $\hat\pi^{u(k_0)}_t$, there exists a constant $K_0=K_0(k_0)>0$ independent of $t$ and the choice of the continuation $\{w_t\}$ such that
\begin{equation}\label{eqn_5n}
\Bigl| \frac{d u^{\langle j;j\rangle}_t}{d t}\Bigr| <K_0.
\end{equation}

By using this inequality, we will prove the following reducing lemma.

\begin{lemma}\label{lem_a_2}
For any $\delta_0>0$, there exists an integer $j_0>0$ independent of $t$ and the choice of a continuation $\{w_t\}$ with $w_t\in \Lambda_t\cap I_t$ such that, for any $j\geq j_0$,
\begin{equation}\label{u_t}
\Bigl| \frac{d u^{\langle 0;j\rangle}_t}{d t}\Bigr| <\delta_0.
\end{equation}
\end{lemma}
\begin{proof}
Choose $\delta>0$ so that $2\delta/(1-\lambda_0)<\delta_0$.
Since $u_t^{\langle j-1;j\rangle}=\zeta^{(j)}(t,u_t^{\langle j;j\rangle})$,
$$\frac{d u^{\langle j-1;j\rangle}_t}{d t}=\frac{\partial \zeta^{(j)}}{\partial t}(t,u_t^{\langle j;j\rangle})+\frac{\partial \zeta^{(j)}}{\partial u}(t,u_t^{\langle j;j\rangle})\frac{d u^{\langle j;j\rangle}_t}{d t}.$$
By (\ref{eqn_4n}), (\ref{eqn_41}) and (\ref{eqn_5n}),
$$
\Bigl|\frac{d u^{\langle j-1;j\rangle}_t}{d t}\Bigr|<\delta+\lambda_0 K_0.
$$
By repeating the same argument, we have
$$
\Bigl|\frac{d u^{\langle 0;j\rangle}_t}{d t}\Bigr|<\delta+\lambda_0\delta+\lambda_0\delta+\cdots +\lambda_0^{j-1}\delta+\lambda^j_0 K_0<\frac{\delta}{1-\lambda_0}+\lambda_0^jK_0.
$$
Thus, our desired inequality holds for all sufficiently large $j$.
\end{proof}

Now, we are ready to prove Accompanying Lemma.

\begin{proof}[Proof of Lemma \ref{lem_1_1}]
Note that $\mathcal{X}_t=\varphi^j_t(\hat \pi_t^{u(k_0);j}(\Lambda_t\cap I_t))$ is the set of crossing points $\eta_t^u$ of $\sigma$ 
and $\Lambda_t$-leaves of $\mathcal{F}_t^{u(k_0+j)}$.
Hence, one can choose $w_t\in \Lambda_t\cap I_t$ so that $\beta_t^{(0)}(u_t^{\langle 0;j\rangle})=\eta_t^u$.
Since $u_t\approx 0$ and $t\approx t_0$, one can suppose that $(\partial \beta^{(0)}/\partial t)(t,u_t)$ is arbitrarily close to the tangent vector $(\partial \beta^{(0)}/\partial t)(t_0,0)=da_t/dt|_{t=t_0}$.
Since moreover
$$
\frac{d\eta_t^u}{dt}(t)=\frac{\partial \beta^{(0)}}{\partial t}(t,u_t^{\langle 0;j\rangle})+
(d\beta^{(0)}_t)_{u_t^{\langle 0;j\rangle}}\Bigl(\frac{du_t^{\langle 0;j\rangle}}{dt}\Bigr).
$$
and $\Vert (d\beta^{(0)}_t)_{u_t}\Vert=1$, Lemma \ref{lem_a_2}
implies the inequality (\ref{accomp}).
This completes the proof.
\end{proof}


\begin{thebibliography}{100}


\bibitem{BC}
M. Benedicks and L. Carleson, 
The dynamics of the H\'enon map, 
 \textit{Ann.\ of Math.} (2) \textbf{133} (1991), no.\ 1, 73--169

\bibitem{BDV}
C. Bonatti, L.\,J. D\'{\i}az and F. Vuillemin, 
Cubic tangencies and hyperbolic diffeomorphisms, 
\textit{Bol.\ Soc.\ Brasil.\ Mat.\ (N.S.)} 
\textbf{29} (1998), no.\ 1, 99--144

\bibitem{Ca}
M. Carvalho,
First homoclinic tangencies in the boundary of Anosov diffeomorphisms,  
\textit{Discrete Contin.\ Dynam.\ Systems} \textbf{4} (1998), no.\ 4, 765--782

\bibitem{Dav}
G.\,J. Davis, Infinitely many co-existing sinks from degenerate homoclinic tangencies.
\textit{Trans.\ Amer.\ Math.\ Soc.} \textbf{323} (1991), no.\ 2, 727--748



\bibitem{DR}
L.\,J. D\'{\i}az and J. Rocha, 
Heterodimensional cycles, partial hyperbolicity and limit dynamics, 
\textit{Fundamenta Math.}, 
\textbf{174}, (2002)
127--187



\bibitem{E}
H. Enrich,
A heteroclinic bifurcation of Anosov diffeomorphisms, 
\textit{Ergodic Theory Dynam.\ Systems} \textbf{18} (1998), no.\ 3, 567--608


\bibitem{Fr}
J. Franks, 
Differentiably $\Omega$-stable diffeomorphisms, 
\textit{Topology} \textbf{11} (1972), 107--113






\bibitem{G93}
S.\,V. Gonchenko, L.\,P.  Shilnikov and D.\,V. Turaev,
On models with non-rough Poincar\'{e} homoclinic curves,
\textit{Physica D},  \textbf{62}, Nos.\ 1-4, (1993) 1--14



\bibitem{G93b}
S.\,V. Gonchenko, L.\,P.  Shilnikov and D.\,V. Turaev,
On the existence of Newhouse regions near systems with
a structurally unstable homoclinic Poincar\'e curve (multidimensional case) 
\textit{Russian Acad. Sci. Dokl. Math.}
\textbf{47} 268--273

\bibitem{G96}
S.\,V. Gonchenko,  L.\,P. Shilnikov and D.\,V. Turaev,
 Dynamical phenomena in systems
with structurally unstable Poincarr\'e homoclinic orbits,
\textit{Interdisc. J. CHAOS}  \textbf{6}(1)  (1996)15--31.

\bibitem{G97}
 S.\,V. Gonchenko,    L.\,P. Shilnikov and  D.\,V. Turaev,
On Newhouse domains of two-dimensional diffeomorphisms
which are close to a diffeomorphism with a structurally unstable heteroclinic cycle, \textit{Proc. Steklov Inst. Math.}
\textbf{216} (1997) 70--118


\bibitem{G99}
 S.\,V. Gonchenko,   L.\,P. Shilnikov and  D.\,V. Turaev,
Homoclinic tangencies of an arbitrary order in Newhouse domains, 
\textit{Itogi Nauki i Tekhniki, Seriya Sovremennaya Matematika iEe Prilozheniya} \textbf{67} (1999)
69--128, 
 \textit{J. Math. Sci.} \textbf{105}  (2001) 1738--1778 (English transl.)
 
 
\bibitem{G02}
S.\,V. Gonchenko, L.\,P.  Shilnikov and O.\,V. Stenkin,  
On Newhouse regions with infinite many stable and
unstable invariant tori,
\textit{Progress in Nonlinear Science: Proc.\ Int.\ Conf.\ Dedicated to A.\,A.  Andronov}, Vol I:
Mathematical Problems of Nonlinear Dynamics University of Nizhny Novgorod, 
(2002) 80--102

\bibitem{G07}
S.\,V. Gonchenko,  L.\,P. Shilnikov and  D.\,V. Turaev, 
Homoclinic tangencies of arbitrarily high orders in
conservative and dissipative two-dimensional maps,
\textit{Nonlinearity} \textbf{20} (2007) 241--275

\bibitem{Kal}
V.\,Y. Kaloshin,
Generic diffeomorphisms with superexponential
growth of number of periodic orbits,
\textit{Commun. Math. Phys.} \textbf{211}(2000) 253--271 

\bibitem{KKY}
I. Kan,  H. Ko\c{c}ak and J.\,A. Yorke, 
Antimonotonicity: concurrent creation and annihilation of periodic orbits.
\textit{Ann.\ of Math.} (2) \textbf{136} (1992), no.\ 2, 219--252


\bibitem{LS}
J.\,S.\,W. Lamb and O.\,V. Stenkin,
Newhouse regions for reversible systems with infinitely
many stable, unstable and elliptic periodic orbits,
\textit{Nonlinearity} \textbf{17} (2004) 1217--1244

\bibitem{Li}
M.-C. Li, 
Point bifurcations and bubbles for a cubic family,  
\textit{J. Difference Equ.\ Appl.} \textbf{9} (2003), no.\ 6, 553--558


\bibitem{LMM}
R. de la Llave, J.\,M. Marco and R. Moriy\'{o}n, Canonical perturbation theory of Anosov systems and regularity results for the 
Liv\v sic cohomology equationm,  \textit{Ann.\ of Math.} (2) \textbf{123} (1986), no.\ 3, 537--611



\bibitem{Man}
R. M\~{a}n\'{e}, The Hausdorff dimension of horseshoes of diffeomorphisms of surfaces, \textit{Bol.\ Soc.\ Brasil.\ Mat.\ (N.S.)}
\textbf{20} (1990), no.\ 2, 1--24


\bibitem{Ma}
J. Mather, Differentiable dynamical systems Appendix, \textit{Bull.\ Amer.\ Math.\ Soc.} 
\textbf{73} (1967), 747--817





\bibitem{DV}
W. de Melo and S. Van Strien, 
One-Dimensional Dynamics (Ergebnisse Der Mathematik Und Ihrer Grenzgebiete, 3 Folge Vol.\ 25) Springer-Verlag (1993)



\bibitem{MT}
J. Milnor and W.\,P. Thurston,
On iterated maps of the interval, 
\textit{Dynamical systems (College Park, MD, 1986--87)}, 465--563, 
Lecture Notes in Math., 1342, Springer, Berlin, 1988


\bibitem{MTr}
J. Milnor and C. Tresser,
On entropy and monotonicity for real cubic maps. With an appendix by Adrien Douady and Pierrette Sentenac, \textit{Comm.\ Math.\ Phys.} \textbf{209} (2000), no.\ 1, 123--178




\bibitem{MV}
L. Mora and M. Viana,
Abundance of strange attractors, \textit{Acta Math.} \textbf{171} (1993), 
no.\ 1, 1--71



 \bibitem{N1}
S.\,E. Newhouse, 
The abundance of wild hyperbolic sets and non-smooth stable sets for diffeomorphisms, 
\textit{Publ.\ Math.\ I.H.\'E.S.} \textbf{50} (1979), 101--151
 
\bibitem{N2}
S.\,E. Newhouse, 
New phenomena associated with homoclinic tangencies, 
\textit{Ergodic Theory Dynam.\ Systems} \textbf{24} (2004), no.\ 5, 1725--1738



\bibitem{PT}
J. Palis and F. Takens, 
Hyperbolicity and sensitive chaotic dynamics at homoclinic bifurcations, Fractal dimensions and infinitely many attractors. Cambridge Studies in Advanced Mathematics 35, Cambridge University Press, Cambridge, 1993


\bibitem{Pol}
M. Pollicott, 
Stability of mixing rates for Axiom A attractors, 
\textit{Nonlinearity} \textbf{16} (2003), 567-578


 
\bibitem{R} 
C. Robinson, 
Bifurcation to infinitely many sinks, 
\textit{Comm.\ Math.\ Phys.} \textbf{90} (1983), no.\ 3, 433--459

\bibitem{R2} 
C. Robinson, Dynamical systems, Stability, symbolic dynamics, and chaos, 
Second ed.\ Studies in Advanced Mathematics, CRC Press, Boca Raton, FL, 1999



\bibitem{S}
S. Sternberg,
On the structure of local homeomorphisms of euclidean $n$-space, II, 
\textit{Amer.\ J. Math.} \textbf{80} (1958), 623--631



\bibitem{T0}
F. Takens, 
Partially hyperbolic fixed points,  
\textit{Topology} \textbf{10} (1971), 133--147




\bibitem{TY}
L. Tedeschini-Lalli and J.\,A. Yorke,
How often do simple dynamical processes have infinitely many coexisting sinks? 
\textit{Comm.\ Math.\ Phys.} \textbf{106} (1986), no.\ 4, 635--657




\bibitem{WY}
Q. Wang and L.-S. Young, 
Strange attractors with one direction of instability,  
\textit{Comm.\ Math.\ Phys.}
\textbf{218}
(2001), no.\ 1, 1--97



\end{thebibliography}
\end{document}